\documentclass{asl}

\usepackage{tikz-cd,tikz,tikz-3dplot}
\usetikzlibrary{shapes, arrows, calc, positioning, shapes.geometric,%
	decorations.pathreplacing, decorations.markings}
\usepackage{adjustbox}
\usepackage{caption}
\usetikzlibrary{mindmap}
\usepackage{color}
\usepackage[final]{hyperref}
\usepackage{graphicx}
\usepackage{tcolorbox}
\usepackage{stmaryrd}
\usepackage[shortcuts]{extdash}
\usepackage{faktor}

\usepackage[utf8]{inputenc}
\usepackage[T1]{fontenc}
\usepackage{CJKutf8}

\usepackage{mathrsfs}
\usepackage{mathabx}
\usepackage{accents}
\usepackage{dsfont}
\usepackage{bussproofs}

\EnableBpAbbreviations%
\numberwithin{equation}{section}

\newtheorem{theorem}{Theorem}[section]
\newtheorem*{theorem*}{Theorem}
\newtheorem{proposition}[theorem]{Proposition}
\newtheorem*{proposition*}{Proposition}
\newtheorem{lemma}[theorem]{Lemma}
\newtheorem*{lemma*}{Lemma}
\newtheorem{corollary}[theorem]{Corollary}
\newtheorem*{corollary*}{Corollary}
\newtheorem{remark}[theorem]{Remark}
\newtheorem*{remark*}{Remark}
\newtheorem{question}[theorem]{Question}
\newtheorem*{question*}{Question}
\newtheorem{mathrule}[theorem]{Rule}

\theoremstyle{definition}
\newtheorem{definition}[theorem]{Definition}
\newtheorem*{definition*}{Definition}

\newtheorem*{example*}{Example}

\newtheorem*{theorem:BiInterpretability-result}%
	{Theorem~\ref{theorem:BiInterpretability-result}}
\newtheorem*{theorem:CThEqCoh-Correspondence}%
	{Theorem~\ref{theorem:CThEqCoh-Correspondence}}	
\newtheorem*{theorem:cth0}%
  {Theorem~\ref{theorem:cth0}}%

\newcommand{\One}{\mathds{1}}
\newcommand{\ep}{\varepsilon}
\makeatletter
\newcommand*{\doublerightarrow}[2]{\mathrel{
  \settowidth{\@tempdima}{$\scriptstyle#1$}
  \settowidth{\@tempdimb}{$\scriptstyle#2$}
  \ifdim\@tempdimb>\@tempdima \@tempdima=\@tempdimb\fi
  \mathop{\vcenter{
    \offinterlineskip\ialign{\hbox to\dimexpr\@tempdima+1em{##}\cr
    \rightarrowfill\cr\noalign{\kern.5ex}
    \rightarrowfill\cr}}}\limits^{\!#1}_{\!#2}}}
\makeatother
\def\fCenter{\ \vdash\ }
\renewcommand{\leq}{\leqslant}
\newcommand\Eff{\mathfrak{F}}
\newcommand\Gee{\mathfrak{G}}
\newcommand\sfF{\mathscr{F}}
\newcommand\sfG{\mathscr{G}}
\newcommand\Eta{\mathfrak{H}}
\providecommand{\C}{}
\renewcommand{\C}{\mathscr{C}}
\newcommand{\M}{\mathscr{M}}
\newcommand{\Cl}{\mathscr{L}}

\newcommand{\T}{\mathscr{T}}
\newcommand{\1}{\mathbf{1}}
\newcommand{\2}{\mathbf{2}}
\newcommand{\Mod}{\normalfont{\text{Mod}}}
\newcommand\X{\mathscr{X}}
\newcommand{\h}{\mathsf{h}}
\newcommand{\CTh}{\mathsf{CTh}}
\newcommand{\CThEq}{\mathsf{CThEq}}
\newcommand{\Th}{\mathsf{Th}}
\newcommand{\ThEq}{\mathsf{ThEq}}
\newcommand{\Coh}{\mathsf{Coh}}
\newcommand{\Bool}{\mathsf{Bool}}
\DeclareMathOperator\Sub{Sub}
\newcommand{\Exact}{\mathsf{ExactCoh}}
\newcommand{\Hom}{\normalfont{\text{Hom}}}
\newcommand{\Pretopos}{\mathsf{Pretopos}}
\newcommand\Set{{
    \mathsf{Set}
}}
\newcommand{\eqdef}{\stackrel{\normalfont{\text{def}}}{=}}
\newcommand{\equivdef}{\stackrel{\normalfont{\text{def}}}{\equiv}}

\DeclareMathOperator{\Dom}{Dom}
\DeclareMathOperator{\dom}{dom}
\DeclareMathOperator{\ob}{ob}
\newcommand{\id}{\mathsf{id}}
\newcommand\eq{{\normalfont{\text{eq}}}}
\newcommand{\ex}{{ex}}
\def\implies{\DOTSB\;\Longrightarrow\;}

\newcommand\SRel[1]{#1{\operatorname{-Rel}}}
\newcommand\SSort[1]{#1{\operatorname{-Sort}}}
\newcommand\SFunc[1]{#1{\operatorname{-Func}}}
\newcommand\SSub[1]{#1{\operatorname{-Sub}}}
\newcommand\SForm[1]{#1{\operatorname{-Form}}}
\newcommand{\formbracket}[1]{\left[{#1}\right]}
\newcommand{\fb}[1]{\formbracket{#1}}
\newcommand{\doublebracket}[1]{\llbracket{#1}\rrbracket}
\newcommand{\fbb}[1]{\doublebracket{#1}}

\newcommand{\ic}[1]{\doublebracket{#1}}
\newcommand{\bwrap}[1]{{\left\{{#1}\right\}}}
\newcommand{\pwrap}[1]{{\left(#1\right)}}
\newcommand{\ul}[1]{\underline{#1}}
\newcommand\ol[1]{\overline{#1}}
\newcommand{\minisub}[2]{#1_{\scalebox{0.65}{\(#2\)}}}
\newcommand{\hw}[1]{\widehat{#1}}
\DeclareMathAccent{\wtilde}{\mathord}{largesymbols}{"65}
\newcommand{\utilde}[1]{{\underaccent{\wtilde}{#1}}}
\newcommand\axiom[2][]{{%
\ifx&#1&%
  \text{#2}%
\else%
  \text{#2}\ensuremath{\pwrap{#1}}%
\fi%
}}%
\newcommand{\pvec}[1]{\vec{#1}\mkern2mu\vphantom{#1}}
\newcommand{\mvec}[1]{\vec{#1}\mkern-2mu\vphantom{#1}}

\newcommand\restr[2]{{
   \left.\kern-\nulldelimiterspace 
    #1 
    \vphantom{\big| } 
    \right|_{#2} 
}}

\newenvironment{bprooftree}%
  {\leavevmode\hbox\bgroup}
  {\DisplayProof\egroup}
\makeatletter
\renewcommand*\subjclass[2][2020]{%
  \def\@subjclass{#2}%
  \@ifundefined{subjclassname@#1}{%
    \ClassWarning{\@classname}{Unknown edition (#1) of Mathematics
      Subject Classification; using '2020'.}%
  }{%
    \@xp\let\@xp\subjclassname\csname subjclassname@#1\endcsname
  }%
}
\@namedef{subjclassname@1991}{%
  \textup{1991} Mathematics Subject Classification}
\@namedef{subjclassname@2000}{%
  \textup{2000} Mathematics Subject Classification}
\@namedef{subjclassname@2010}{%
  \textup{2010} Mathematics Subject Classification}
\@namedef{subjclassname@2020}{%
  \textup{2020} Mathematics Subject Classification}
\@xp\let\@xp\subjclassname\csname subjclassname@2020\endcsname
\makeatother

\title{%
	Bicategories, Biequivalence, and Bi\-/Interpretability%
}

\author{Anthony D'Arienzo}
\revauthor{D'Arienzo, Anthony}
\twoaddress{%
	Department of Mathematics\\
	The University of Illinois\\ 
	Urbana, IL 61801 USA%
}{%
	Department of Mathematics\\
	Princeton University\\
	Princeton, NJ 08544
}{3}
\curraddr{%
	Department of Mathematics\\
	The University of Illinois\\
	Urbana, IL 61801 USA%
}
\email{apd6@illinois.edu}
\genaddr{Correspondence address}{%
	Anthony P. D'Arienzo\\
	Department of Mathematics\\
	The University of Illinois\\
	Urbana, IL 61801 USA%
}

\author{Vinny Pagano}
\revauthor{Pagano, Vinny}
\address{%
	Department of Mathematics\\
	Princeton University\\
	Princeton, NJ 08544
}
\email{vpagano@alumni.princeton.edu}

\author{Ian M.J. McInnis}
\revauthor{McInnis, Ian}
\address{%
	Department of Mathematics\\
	Princeton University\\
	Princeton, NJ 08544
}%
\email{imj.mcinnis@gmail.com}


\begin{document}
\begin{abstract}
	We make explicit the correspondence between syntax and syntactic categories
	for coherent first\-/order logic, providing a categorical characterization of
	bi\-/interpretability. This is done by creating a biequivalence between a
	bicategory of coherent theories and the (strict) bicategory of coherent
	categories. While the biequivalence concerns the stronger
	equality\-/preserving bi\-/interpretability, we use it to obtain a necessary
	and sufficient condition for two theories to be bi\-/interpretable in
	general, by relating the exact completions of their syntactic categories.
	These results extend analogously to familiar fragments of first\-/order
	logic, thereby clarifying the long\-/intuited relation between logical syntax
	and syntactic categories.
\end{abstract}

\maketitle%

\tableofcontents%

\section{Introduction} 

We formalize and answer the following questions. (1) What kind of categorical
structure do (coherent) predicate theories form? (2) How do
bi\-/interpretability (see~\cite{hodges_1993}) and other notions of equivalence
of theories fit into this categorical structure? (3) How does this structure
interact with the syntactic category and internal logic operations
of~\cite{MakkaiReyes1977}? These questions have been explored in the
literature~\cite{%
  AHLBRANDT198663,%
  Friedman2014WhenBI,%
  harnik2011,%
  MakkaiReyes1977,%
  mceldowney_2020,%
  tsementzis2015syntactic,%
  Visser2006ProlegomenaTT,%
  Visser_categoriesof%
}; however there have been no comprehensive, unifying answers yet.

For example, there are five recent works~\cite{%
    halvorson-logic,%
    harnik2011,%
    mceldowney_2020,%
    tsementzis2015syntactic,%
    washington2018%
} that compare bi\-/interpretability, \emph{Morita equivalence}
(see~\cite{barrett_halvorson_2016,johnstone2002sketches,tsementzis2015syntactic}),
and the classifying pretopos of a theory (the pretopos completion of its
syntactic category; see \cite{johnstone2002sketches,MakkaiReyes1977}).
In~\cite{tsementzis2015syntactic}, Tsementzis showed that two theories \(T\)
and \(T'\) have equivalent classifying pretopoi if and only if they are Morita
equivalent in the sense of~\cite{barrett_halvorson_2016}. Tran-Hoang\footnote{%
	Formerly named McEldowney.
} in~\cite{mceldowney_2020} proved that Morita equivalence in the sense
of~\cite{barrett_halvorson_2016} is a strictly coarser relation than
bi\-/interpretability, and the two notions coincide if and only if coproduct
Morita extensions may be eliminated. This shows that, in general, classifying
pretopoi do not classify bi\-/interpretability; a smaller `classifying' category
is needed for a categorical analogue of bi\-/interpretability. However,
reconciling this observation with Harnik's work in~\cite{harnik2011} is subtle.
Harnik makes the conjecture: \emph{%
	for first\-/order theories \(T\) and \(T'\) in finite languages, the most
	general notion of interpretation is a logical functor \(\C(T) \to
	\C({T'}^\eq)\), so \(T\) and \(T'\) are bi\-/interpretable if and only if
	\(\C(T^\eq)\) and \(\C({T'}^\eq)\) are equivalent categories, where \(\C(T)\)
	denotes the syntactic category of \(T\) and \(T^\eq\) denotes Shelah's
	elimination of imaginaries construction (see discussion
	following~\cite[Definition 6.2]{harnik2011}).
}%
Indeed Harnik proves, for \emph{proper} theories, that \(\C(T^\eq)\) is
equivalent to the classifying pretopos of \(T\)~\cite[Theorem 5.3]{harnik2011},
whereas in general \(\C(T^\eq)\) is at least the exact completion of \(\C(T)\)
as in~\cite{MakkaiReyes1977}. This implies an important consequence of Harnik's
conjecture: Morita equivalence and bi\-/interpretability coincide for proper
theories. This consequence was proven by Washington in~\cite{washington2018} 
and further explored in~\cite{halvorson-logic}; however, Washington's proof 
does not use syntactic categories or Shelah's completion, leaving Harnik's 
conjecture unsettled. On the other hand, observe that \(T^\eq\) freely adjoins 
only quotients to \(T\). Since coproducts are not adjoined, this hints at why 
coproduct Morita extensions must be eliminated in order to recover 
bi\-/interpretability. Proving some form of Harnik's conjecture would yield a 
categorical characterization of bi\-/interpretability, similar to Tsementzis'
characterization of Morita equivalence (in the sense
of~\cite{barrett_halvorson_2016}) in terms of classifying pretopoi.
Furthermore, it would explain the \emph{proper}
hypothesis in Washington's theorem: proper theories are those theories \(T\)
for which the exact completion \(\C(T^\eq)\) is also a pretopos.

There is also a converse direction to this rationale. While interpretations
could be identified with logical functors \(\C(T) \to \C({T'}^\eq)\), we
may also ask \emph{what kind of syntactic operation corresponds to a logical
functor \(\C(T) \to \C(T')\)?} By composing such a functor with the exact
completion \(\C(T') \to \C({T'}^\eq)\), we could identify these functors
with a special class of interpretations. The motivation for Harnik's hypothesis
suggests that this class should correspond to \emph{equality\-/preserving
interpretations} (called \emph{simple interpretations} in~\cite{harnik2011}). In turn
these interpretations pick out the stricter notion of
\emph{equality\-/preserving bi\-/interpretability}. In set theory and model
theory, non\-/equality\-/preserving bi\-/interpretability is often
considered~\cite{%
  AHLBRANDT198663,%
  halvorson-logic,%
  hodges_1993,%
  Robinson1952-ROBU-4,%
  roque_freire_hamkins_2021%
}, whereas the equality\-/preserving case is considered in a handful of other
model theory publications~\cite{%
  Button-Walsh2018,%
  Friedman2014WhenBI,%
  harnik2011,%
  Marker2002,%
  Tarski1953-TARUT-3%
}. Therefore understanding how syntactic categories interact with
equality\-/preserving bi\-/interpretability puts these works and their
applications in context.

We answer the three posed questions by formalizing a purely syntactic theory
of interpretations and proving Harnik's conjecture.
Following~\cite{halvorson-logic,washington2018}, we call these interpretations
\emph{translations}. Similar to~\cite{Visser_categoriesof}, we show that these
translations fit into a bicategory \(\CTh_0\) whose 0\-/cells are (coherent)
theories, 1\-/cells are translations, and 2\-/cells are \emph{t\-/maps}
(see~\cite{halvorson-logic,washington2018}, the \emph{i\-/maps}
of~\cite{Visser_categoriesof}, or the \emph{homotopies} in~\cite{AHLBRANDT198663}). We
account for the equality\-/preserving situation by showing that the
sub\-/bicategory \(\CThEq\) spanned by equality\-/preserving translations is
biequivalent to the familiar 2\-/category \(\Coh\) of coherent categories, 
coherent functors, and natural transformations. We prove the following main theorems.

\begin{theorem:cth0}
  \(\CTh_0\) is a well\-/defined bicategory.
\end{theorem:cth0}
\begin{theorem:CThEqCoh-Correspondence}
  The syntactic category and internal logic operations of~\cite{MakkaiReyes1977}
  extend to a biequivalence \(\C: \CThEq \to \Coh\) and \(\T: \Coh \to \CThEq\).
  In particular, two coherent theories \(T_1\) and \(T_2\) are
  bi\-/interpretable via equality\-/preserving translations if and only if
  \(\C(T_1)\) and \(\C(T_2)\) are equivalent categories.
\end{theorem:CThEqCoh-Correspondence}
\begin{theorem:BiInterpretability-result}
  Two coherent theories \(T_1\) and \(T_2\) are bi\-/interpretable if and only
  if \(\C(T_1^\eq)\) and \(\C(T_2^\eq)\) are equivalent categories.
\end{theorem:BiInterpretability-result}

These have consequences in categorical logic and model theory beyond
bi\-/interpretability. In Section 4 we show that the units of the biequivalence
between \(\CThEq\) and \(\Coh\) are the \emph{canonical interpretations}
of~\cite{MakkaiReyes1977}. In Section 6 we show that our theory of translations
readily extends to classical and intuitionistic first\-/order theories, and we
sketch a characterization of Morleyization as a bi\-/adjunction using the
biequivalence of Theorem~\ref{theorem:CThEqCoh-Correspondence}. We also discuss
the directly categorical aspects of the collection of coherent theories, such
as the existence of limits and colimits. Section 6 ends with a discussion on
Morita equivalence and its relation to Harnik's conjecture.

Our theory of translations builds upon the works of others. Our starting point
is Szczerba's work~\cite{Szczerba1977}, which argues that it is natural to
allow translations to assign a single variable in one language to finite
strings of variables in another language. Szczerba's notion of translation was
honed by van\ Benthem and Pearce~\cite{van1984mathematical}, who also hint
about the connection between translations and functors. Their results were
incorporated into the work of Visser and collaborators
in~\cite{Visser_categoriesof}, in which they provide a sharp definition for a
(strict) 2\-/category of theories, called \(\mathsf{INT}^{\mathsf{iso}}\).
However, the study of \(\mathsf{INT}^{\mathsf{iso}}\) is insufficient for
answering the three posed questions or for proving our three main theorems.
This is for three key reasons. (1) Visser et al.\ only consider \emph{purely
relational} signatures. Therefore a more extensive framework is needed to study
theories in signatures with function symbols. This restriction was a decisive
simplification: adding function symbols makes the collection of theories into a
\emph{weak} 2\-/category (bicategory) due to the appearance of nontrivial
unitors. (2) The 2\-/cells of \(\mathsf{INT}^{\mathsf{iso}}\) are all
invertible, so a biequivalence with \(\Coh\)---or any strict 2\-/category of
syntactic categories---is impossible; this prevents a direct relation between
translations and functors. (3) Visser et al.\ do not compare
\(\mathsf{INT}^{\mathsf{iso}}\) to \(\Coh\). A central thesis
of~\cite{MakkaiReyes1977} is that (coherent) theories may be faithfully
replaced with their syntactic categories. Finding an analogue for
bi\-/interpretability in categorical logic is tantamount to making this
comparison.

The categorical formalism developed by~\cite{MakkaiReyes1977} and refined by,
e.g.,~\cite{maclane-moerdijk} and~\cite{johnstone2002sketches} is not
sufficient for our goals either. It is unclear what is preserved when moving
between syntactic and categorical realms using the syntactic category
operation \(\C\) and internal logic operation \(\T\). Among the results proven
in~\cite{MakkaiReyes1977} are two natural bijections between models of a
coherent theory and coherent functors: \(\Mod(T,D) \cong \Coh(\C(T),D)\) and
\(\Coh(C,D) \cong \Mod(\T(C),D)\). These bijections show that the categories
of models of \(T\) and \(\T\C(T)\) are equivalent, i.e., \(T\) and \(\T\C(T)\)
are \emph{categorically equivalent} (see~\cite{barrett_halvorson_2016}), but this is a
strictly weaker notion of equivalence than bi\-/interpretability. By finessing
the two ``main facts'' of~\cite{MakkaiReyes1977} (Theorems 2.4.5 and 3.2.8
of~\cite{MakkaiReyes1977}), it is reasonable to believe that \(T\) and
\(\T\C(T)\) are more akin than categorical equivalence. Our first two main
theorems make this precise.

In order to provide a sufficiently extensive account of bi\-/interpretability,
we develop the theory of translations from first principles. This yields our 
bicategories of theories \(\CTh_0\) and \(\CThEq\), where
\(\CTh_0\) is the coherent cousin of the 2\-/category \(\Th_0\) found
in~\cite{halvorson-logic,washington2018}. To the best of our
knowledge, this paper is the first to explicitly define bicategories of
theories while also showing that the relevant coherence diagrams are satisfied.
Aratake's account in~\cite{AraH:2018} of bi\-/interpretability and syntactic
categories is the closest example of related research. However, they
do not prove a biequivalence and only consider proper theories. This
assumption is too strong in that it prevents characterizing
bi\-/interpretability via exact completions of syntactic categories. We explore
this point at the end of Section 6. We do not assume our theories are generally
proper. In a recent paper~\cite{Kamsma2023typespace}, Kamsma defines (though
does not prove the validity of) a 2-category of coherent theories with
interpretations and t-maps rather like ours. However, they only consider purely
relational theories, making this a \textit{strict} 2-category. This is still
useful; our work makes obvious that Kamsma's \(\frak{CohTheory}\) is 2-equivalent to
\(\Coh\), and he proves that it's dual to the 2-category of \textit{type
space functors}, i.e., contravariant functors from the category of finite sets
to the category of spectral spaces. This is a sort of Stone Duality for
theories. 

Lastly, many ideas we detail in this paper are, to a certain degree, already
implicit in Makkai and Reyes' monograph~\cite{MakkaiReyes1977}. Indeed, Makkai
and Reyes suggest that their construction of the syntactic category of a theory
could be considered the object part of a functor~\cite[Chapter
8]{MakkaiReyes1977}. However, that monograph does not provide precise notions
of translation and t\-/map supplied by Visser et al.\ and refined by us. In
this sense, our results can be seen as finally validating an intuition of
Makkai and Reyes regarding the relation between predicate theories and
categorical logic.

\section{Bicategory Theory}
\label{section:bicategory-theory}

This section consists of definitions. Our work uses the formalism in Gray's
book~\cite{Gray_1974}, adjusting some notation and terminology to conform with
modern usage. Namely, what we call \emph{vertical composition},
\emph{horizontal composition}, \emph{pseudofunctor}, \emph{pseudonatural
transformation}, and \emph{strict 2\-/category}, Gray calls \emph{weak
composition}, \emph{strong composition}, \emph{homomorphism},
\emph{quasinatural transformation}, and \emph{2\-/category} respectively. We
shall use \(\cdot\) for vertical composition of 2\-/cells, \(\circ\) for
horizontal composition of 2\-/cells, and juxtaposition or \(\circ\) for
composition of 1\-/cells. We will often use topology\-/motivated language when
reasoning in a bicategory. For example, an invertible 2\-/cell will be called a
\emph{homotopy}, and our notion of \emph{weak equivalence} is called
\emph{homotopy equivalence}.\footnote{%
	Calling weak equivalences ``homotopy equivalences'' comes from the
	influence of homotopy theory in the development of higher categories,
	paralleling, e.g.,~\cite{reihlhomotopy,shulmanhomotopy}. Some
	sources, like~\cite{halvorson-logic}, call bi\-/interpretability
	``weak intertranslatability'' or ``homotopy equivalence'', coming from
	the implicitly 2\-/categorical tools needed in this notion of equivalence of
	theories.
}%

\subsection{Bicategories}
\label{subsection:bicategories}

\begin{definition}[\cite{Gray_1974}, I,3.1]
	A \textbf{bicategory} \(\mathsf{C}\) is a collection of the following.%

	(BC1) A collection \(\ob(\mathsf{C})\) of objects in \(\mathsf{C}\) (instead
	of writing \(X \in \ob(\mathsf{C})\), we adopt the usual abuse of notation and
	write \(X \in \mathsf{C}\)).
	
	(BC2) A hom-category \(\mathsf{C}(X,Y)\) for every pair \(X,Y \in
	\mathsf{C}\). The objects are \textbf{1-cells} \(f: X \to Y\) and the
	morphisms are \textbf{2\-/cells} \(\chi : f \Rightarrow g\). Given a 1-cell 
	\(f: X \to Y\), the identity 2-cell \(f \Rightarrow f\) is written \(\One^f\).
	Hom\-/category composition is called \textbf{vertical composition}, and an
	invertible 2\-/cell is called a \textbf{homotopy}.

	(BC3) A functor \(\circ_{XYZ} : \mathsf{C}(Y,Z) \times \mathsf{C}(X,Y) \to
	\mathsf{C}(X,Z)\) for every triple of objects \(X,Y,Z \in \mathsf{C}\). This 
	functor defines composition of 1-cells on objects and \textbf{horizontal 
	composition} of 2-cells on morphisms. 

	(BC4) A functor \(1_X: 1 \to \mathsf{C}(X,X)\) for every object \(X \in
	\mathsf{C}\), identifying a \textbf{weak identity} 1-cell.

	(BC5) Natural isomorphisms whose components are \(a_{hgf}: (hg)f \Rightarrow
	h(gf)\), which are commonly referred to as \textbf{associators}.

	(BC6) Natural isomorphisms whose components are \(l_f: 1_Y f \Rightarrow f\)
	and \(r_f: f 1_X \Rightarrow f\). They are called \textbf{left} and
	\textbf{right} \textbf{unitors}, respectively.

	(BC7) The associators and unitors satisfy
    coherence laws called the
	\textbf{pentagon} and \textbf{triangle identities} found in,
	e.g.,~\cite{Gray_1974}. Equivalently, \(a\) implies that \(\circ\) is
	associative up to commutativity of the pentagonal diagram, and \(l,r\) imply
	that \(1_X\) is an identity up to the triangular diagram for each \(X
	\in \mathsf{C}\); see Figure \ref{fig:pentagon-triangle-identities}.
\end{definition}

\begin{figure}[h]
	\centering
	\[\begin{array}{cc}
		\begin{tikzpicture}
			\node[
				regular polygon,
				regular polygon sides=5,
				minimum width=32mm,
			] (PG) {}
                (PG.corner 1) node (PG1) {\(((kh)g)f\)}
			(PG.corner 2) node (PG2) {\((k(hg))f\)}
			(PG.corner 3) node (PG3) {\(k((hg)f)\)}
			(PG.corner 4) node (PG4) {\(k(h(gf))\)}
			(PG.corner 5) node (PG5) {\((kh)(gf)\)}
			;
                \draw[->] (PG1)--(PG2) node[midway, above left] {\scriptsize\(a \circ \One\)};
			\draw[->] (PG1)--(PG5) node[midway, above right] {\scriptsize\(a\)};
			\draw[->] (PG2)--(PG3) node[midway, mid left] {\scriptsize\(a\)};
                \draw[->] (PG3)--(PG4) node[midway, below] {\scriptsize\(\One \circ a\)};
                \draw[->] (PG5)--(PG4) node[midway, mid right] {\scriptsize\(a\)};
		\end{tikzpicture} &%
		\begin{tikzpicture}
			\node[
				regular polygon,
				regular polygon sides=3,
				minimum width=32mm,
			] (TG) {}
			(TG.corner 1) node (TG1) {\((g 1_Y) f\)}
			(TG.corner 2) node (TG2) {\(g (1_Y f)\)}
			(TG.corner 3) node (TG3) {\(gf\)}
			;
			\draw[->] (TG1)--(TG2) node[midway, above left] {\scriptsize\(a\)};
			\draw[->] (TG1)--(TG3) node[midway, above right] {\scriptsize\(r \circ \One\)};
			\draw[->] (TG2)--(TG3) node[midway, below] {\scriptsize\(\One \circ l\)};
		\end{tikzpicture}
	\end{array}\]
	\caption{%
		\label{fig:pentagon-triangle-identities}%
		The Pentagon and Triangle Identities%
	}
\end{figure}

\begin{definition}
	We say that a bicategory \(\mathsf{C}\) is \textbf{small} if
	\(\ob(\mathsf{C})\) is a set and each hom-category \(\mathsf{C}(X,Y)\)
	is a small category for any \(X,Y \in \mathsf{C}\).
	Moreover, \(\mathsf{C}\) is a \textbf{strict 2-category} if
	the associators and unitors are identity 2-cells, i.e., \(\mathsf{C}\) is a
	\(\mathsf{Cat}\)\-/enriched category.
\end{definition}

An important example for our work of a strict 2\-/category is \(\Coh\), the
2\-/category of small coherent categories, coherent functors, and
natural transformations. Recall that a category \(C\) is \textbf{coherent} if it has
finite limits as well as pullback\-/stable images and joins, and a functor 
\(\Eff : C \to D\) between coherent categories is \textbf{coherent} if it preserves said properties.
Some use the term \emph{logical}; see, e.g.,~\cite{harnik2011}.

Additionally, if a category admits quotients, it is called a
\textbf{(Barr\-/)exact} category, and if an exact category admits coproducts,
it is called a \textbf{pretopos}. For an introduction to the theory of
coherent categories, exact categories, and pretopoi, see,
e.g.,~\cite{MakkaiReyes1977}.

In this paper, the associators will be trivial, but the unitors will not. This
situation arises naturally from the logical perspective, and is why we work
with bicategories instead of strict 2\-/categories.

\begin{definition}
	Let \(\mathsf{C}\) be a bicategory, and \(X,Y\) objects of \(\mathsf{C}\). A
	\textbf{weak equivalence} (or \textbf{homotopy equivalence}) is a pair of
	1-cells \(f: X \to Y\) and \(g: Y \to X\) such that \(1_X \simeq gf\) and
	\(1_Y \simeq fg\), where \(\simeq\) denotes isomorphism in a hom-category
	(``homotopy''). In this situation, \(f\) and \(g\) are said to be
	\textbf{homotopy equivalences}, and \(g\) is called a \textbf{homotopy
	inverse} of \(f\). Moreover, \(X\) and \(Y\) are said to be \textbf{homotopy
	equivalent}, written \(X \approx Y\).

	We reserve the symbol \(\cong\) for isomorphism in a category or
	bicategory.
\end{definition}

\subsection{Pseudofunctors and Pseudonatural Transformations}
\label{subsection:pseudofunctors}

\begin{definition}[\cite{Gray_1974}, I,3.2]\label{def:pseudofunctor}
	Let \(\mathsf{C},\mathsf{D}\) be bicategories. A \textbf{pseudofunctor}
	\(\sfF: \mathsf{C} \to \mathsf{D}\) is a collection of the following data.%

	(PF1) A function \(\sfF: \ob(\mathsf{C}) \to \ob(\mathsf{D})\) sending objects
	of \(\mathsf{C}\) to objects of \(\mathsf{D}\).

	(PF2) A functor \(\sfF_{XY}: \mathsf{C}(X,Y) \to \mathsf{D}(\sfF(X),\sfF(Y))\)
	for any pair \(X,Y \in \mathsf{C}\).

	(PF3) For every triple of objects \(X,Y,Z \in \mathsf{C}\), a natural
	isomorphism%
	\[%
		\begin{tikzcd}
			\mathsf{C}(Y,Z) \times \mathsf{C}(X,Y) \arrow[r,"\circ"]
			\arrow[d,"\sfF_{YZ} \times \sfF_{XY}",swap] & \mathsf{C}(X,Z)
			\arrow[d,"\sfF_{XZ}"] \arrow[d,"\sfF_{XZ}"]\\
			\mathsf{D}(\sfF(Y),\sfF(Z)) \times \mathsf{D}(\sfF(X),\sfF(Y))
			\arrow[r,"\circ",swap] \arrow[ur,Rightarrow,"\sfF_{XYZ}",swap] &
			\mathsf{D}(\sfF(X),\sfF(Z))
		\end{tikzcd}
	\]%
	defined by homotopies \(\sfF_{gf}: \sfF(g) \sfF(f) \Rightarrow \sfF(gf)\). 
        We call it the \textbf{compositor}.

	(PF4) For every object \(X \in \mathsf{C}\), a homotopy 
	\(\sfF_{1_X}: 1_{\sfF(X)} \Rightarrow \sfF(1_X)\), 
	referred to hereafter as the \textbf{identitor}.

	(PF5) The compositors and identitors make the diagrams in
	Figure~\ref{fig:pseudofunctor-coherence} commute.
\end{definition}

\begin{figure}[h]
    \begin{center}
        \setlength{\tabcolsep}{-3pt} 
    \[\begin{tabular}{ccc}
            \begin{tabular}{@{}l@{}} 
            \begin{tikzpicture}
                \node[
                    regular polygon,
                    regular polygon sides=6,
                    minimum width=32mm,
                    rotate=90,
                ] (HG) {}
                (HG.corner 1) node (HG1) {\((\sfF(h) \sfF(g))\sfF(f)\)}
                (HG.corner 2) node (HG2) {\(\sfF(h) (\sfF(g) \sfF(f))\)}
                (HG.corner 3) node (HG3) {\(\sfF(h) \sfF(gf)\)}
                (HG.corner 4) node (HG4) {\(\sfF(h(gf))\)}
                (HG.corner 5) node (HG5) {\(\sfF((hg)f)\)}
                (HG.corner 6) node (HG6) {\(\sfF(hg)\sfF(f)\)}
                ;
                \draw[->]%
                    (HG1)--(HG6) node[midway, above left, near start] {\scriptsize\(\sfF_{hg}\circ \One\)}; 
                \draw[->]%
                    (HG6)--(HG5) node[midway, above right, near end] {\scriptsize\(\sfF_{(hg)f}\)}; 
                \draw[->]%
                    (HG5)--(HG4) node[midway, left] {\scriptsize\(\sfF(a)\)}; 
                \draw[->]%
                    (HG1)--(HG2) node[midway, right] {\scriptsize\(a\)};
                \draw[->]%
                    (HG2)--(HG3) node[midway, below left, near start] {\scriptsize\(\One\circ\sfF_{gf}\)}; %
                \draw[->]%
                    (HG3)--(HG4) node[midway, below right, near end] {\scriptsize\(\sfF_{h(gf)}\)}; %
            \end{tikzpicture}%
            \end{tabular}%
        & %
            \begin{tabular}{@{}l@{}}
            \begin{tikzpicture}
                \node[
                    regular polygon,
                    regular polygon sides=4,
                    minimum width=28.5mm,
                ] (SLG) {}
                (SLG.corner 1) node (SLG1) {\(\sfF(f 1_X)\)}
                (SLG.corner 2) node (SLG2) {\(\sfF(f) \sfF(1_X)\)}
                (SLG.corner 3) node (SLG3) {\(\sfF(f) 1_{\sfF(X)}\)}
                (SLG.corner 4) node (SLG4) {\(\sfF(f)\)}
                ;
                \draw[->]%
                    (SLG2)--(SLG1) node[midway, above] {\scriptsize\(\sfF_{f 1_X}\)};
                \draw[->]%
                    (SLG1)--(SLG4) node[midway, left] {\scriptsize\(\sfF(r)\)}; 
                \draw[->]%
                    (SLG2)--(SLG3) node[midway, right] {\scriptsize\(\One \circ \sfF_{1_X}\)};%
                \draw[->]%
                    (SLG3)--(SLG4) node[midway, below] {\scriptsize\(r\)};
            \end{tikzpicture}%
            \end{tabular}%
        &%
            \begin{tabular}{@{}l@{}}
            \begin{tikzpicture}
                \node[
                    regular polygon,
                    regular polygon sides=4,
                    minimum width=28.5mm,
                ] (SRG) {}
                (SRG.corner 1) node (SRG1) {\(\sfF(1_Y f)\)}
                (SRG.corner 2) node (SRG2) {\(\sfF(1_Y) \sfF(f)\)}
                (SRG.corner 3) node (SRG3) {\(1_{\sfF(Y)} \sfF(f)\)}
                (SRG.corner 4) node (SRG4) {\(\sfF(f)\)}
                ;
                \draw[->]%
                    (SRG2)--(SRG1) node[midway, above] {\scriptsize\(\sfF_{1_Y f}\)};
                \draw[->]%
                    (SRG1)--(SRG4) node[midway, left] {\scriptsize\(\sfF(l)\)};
                \draw[->]%
                    (SRG2)--(SRG3) node[midway, right] {\scriptsize\(\sfF_{1_Y} \circ \One\)};%
                \draw[->]%
                    (SRG3)--(SRG4) node[midway, below] {\scriptsize\(l\)};
            \end{tikzpicture}%
            \end{tabular}%
    \end{tabular}\]
	\caption{%
		\label{fig:pseudofunctor-coherence}%
		The Pseudofunctor Coherence Laws%
	}
        \end{center}
\end{figure}

\begin{definition}[\cite{Gray_1974}, I,3.2]\label{def:pseudofunctor-compose}
	Given a pair of pseudofunctors \(\mathsf{A} \xrightarrow{\sfF} \mathsf{B}
	\xrightarrow{\sfG} \mathsf{C}\), we can compose them to yield a pseudofunctor
	\(\sfG\sfF: \mathsf{A} \to \mathsf{C}\) as follows:%
	\begin{itemize}
		\item The function \(\sfG\sfF : \ob(\mathsf{A}) \to \ob(\mathsf{C})\) is \(\ob(\mathsf{A})
			\xrightarrow{\sfF} \ob(\mathsf{B}) \xrightarrow{\sfG} \ob(\mathsf{C})\),
		\item For any two objects \(X,Y \in \mathsf{A}\), \((\sfG\sfF)_{XY}\) is
			the composition \(\sfG_{\sfF(X) \sfF(Y)}\sfF_{XY}\),
		\item For the compositor, \((\sfG\sfF)_{XYZ} \eqdef \sfG_{\sfF(X) \sfF(Y) \sfF(Z)} 
        \boxminus \sfF_{XYZ}\), i.e., \((\sfG\sfF)_{gf} = \sfG(\sfF_{gf}) \cdot
			\sfG_{\sfF(g) \sfF(f)}\) for any composable pair of 1-cells \(f,g\) in
			\(\mathsf{A}\),
		\item For the identitor, \((\sfG\sfF)_{1_X} \eqdef \sfG(\sfF_{1_X}) \circ
			\sfG_{1_{\sfF(X)}}\).
	\end{itemize}
\end{definition}

\begin{definition}[\cite{Gray_1974}, I,3.3]\label{def:PNTi}
	Let \(\sfF,\sfG: \mathsf{C} \to \mathsf{D}\) be a pair of pseudofunctors. A
	\textbf{pseudonatural transformation} \(\eta: \sfF \Rightarrow \sfG\) is a collection of the following.%
	
	(PNT1) A 1-cell \(\eta_X: \sfF(X) \to \sfG(X)\) for every object \(X \in \mathsf{C}\).

	(PNT2) A family of natural transformations%
	\[%
		\begin{tikzcd}
			\mathsf{C}(X,Y) \arrow[r,"\sfF_{XY}"] \arrow[d,"\sfG_{XY}",swap] &
			\mathsf{D}(\sfF(X),\sfF(Y)) \arrow[d,"\pwrap{\eta_Y}_*"]\\
			\mathsf{D}(\sfG(X),\sfG(Y)) \arrow[r,"\pwrap{\eta_X}^*"]
			\arrow[ur,Rightarrow,"\eta_{XY}"] & \mathsf{D}(\sfF(X), \sfG(Y)),
		\end{tikzcd}
	\]%
	the components of which are 2-cells \(\eta_f: \sfG(f) \eta_X \Rightarrow
	\eta_Y \sfF(f)\), where each \(\eta_f\) is natural so that for
	any \(\chi: f \Rightarrow g\) we have a commutative diagram; see
	Figure~\ref{fig:pseudonatural-coherence}.

	(PNT3) For any object \(X \in \mathsf{C}\), \(\eta_{1_X}\) is compatible with
	the unitors of \(\mathsf{D}\), making a commutative diagram; see
	Figure~\ref{fig:pseudonatural-coherence}.

	(PNT4) For any composable pair of 1-cells \(f: X \to Y\) and \(g: Y \to Z\) in 
        \(\mathsf{C}\), \(\eta_{gf}\) is compatible with composition, making a commutative diagram; see Figure~\ref{fig:pseudonatural-coherence}.

	In the case that the 2-cells \(\eta_f: \sfG(f) \eta_X \Rightarrow \eta_Y 
        \sfF(f)\) are homotopies for any 1-cell \(f\), we say that \(\eta: \sfF \Rightarrow
	\sfG\) is a \textbf{pseudonatural homotopy}.
\end{definition}

\begin{figure}[h]
	\centering%
	\[\begin{tabular}{cc}%
			\begin{tikzcd}
				\sfG(f) \eta_X \arrow[r,Rightarrow,"\sfG(\chi) \circ \One"]
				\arrow[d,Rightarrow,swap,"\eta_f"] & \sfG(g) \eta_X \arrow[d,Rightarrow,
				"\eta_g"]\\
				\eta_Y \sfF(f) \arrow[r,Rightarrow,"\One \circ \sfF(\chi)",yshift=.1em] & \eta_Y \sfF(g)
			\end{tikzcd} 
                &%
			\begin{tikzcd}
				1_{\sfG(X)} \eta_X \arrow[r,"r^{-1} \cdot l"] \arrow[d,"\sfG_{1_X} \circ
				\One", swap] & \eta_X 1_{\sfF(X)} \arrow[d,"\One \circ \sfF_{1_X}"]\\
				\sfG(1_X) \eta_X \arrow[r,"\eta_{1_X}"] & \eta_X \sfF(1_X)
			\end{tikzcd}%
	\end{tabular}\]%
	\[\begin{tikzcd}[column sep=small]
                \pwrap{\sfG(g)\sfG(f)}\eta_X \arrow[r,"\pwrap{\One \circ \eta_f} \cdot a",yshift=.3em] \arrow[d,"\sfG_{gf}\circ\One"] &%
                \sfG(g) \pwrap{\eta_Y \sfF(f)} \arrow[r,"a^{-1}",yshift=.3em] &%
                \pwrap{\sfG(g) \eta_Y} \sfF(f) \arrow[r,"a \cdot \pwrap{\eta_g\circ \One}",yshift=.3em] &%
                \eta_Z \pwrap{\sfF(g) \sfF(f)} \arrow[d,"\One\circ\sfF_{gf}"]\\%
                \sfG\pwrap{gf}\eta_X \arrow[rrr,"\eta_{gf}"] & & &%
                \eta_Z \sfF\pwrap{gf}
        \end{tikzcd}\]%
	\caption{%
		\label{fig:pseudonatural-coherence}%
		The Pseudonatural Transformation Coherence Laws
	}%
\end{figure}

The canonical definition for equivalence of bicategories (\emph{biequivalence})
does not appear in~\cite{Gray_1974}, but it appears in, e.g., Section 3 of Lack's
article~\cite{zbMATH01879624}.

\begin{definition}\label{def:biequivalence}
	Let \(\mathsf{C},\mathsf{D}\) be bicategories. A \textbf{biequivalence}
	\(\mathsf{C} \rightleftarrows \mathsf{D}\) is a pair of pseudofunctors
	\(\sfF: \mathsf{C} \to \mathsf{D}\) and \(\sfG: \mathsf{D} \to \mathsf{C}\)
	along with a pair of pseudonatural homotopies \(\epsilon: \id_\mathsf{C}
	\Rightarrow \sfG\sfF\) and \(\delta: \id_\mathsf{D} \Rightarrow \sfF\sfG\)
	such that \(\epsilon_X: X \to \sfG\sfF\pwrap{X}\) and
	\(\delta_Y: Y \to \sfF\sfG\pwrap{Y}\) are homotopy equivalences for any
	objects \(X \in \mathsf{C}\) and \(Y \in \mathsf{D}\).
\end{definition}

Biequivalences preserve objects up to homotopy equivalence. For our purposes a
biequivalence using \emph{strict} 2\-/functors (see~\cite{Gray_1974}) is
impossible: see Remark~\ref{remark:WhyWeCantGetStrict2Eq}.

\subsection{The Homotopy Category of a Bicategory}

Given a bicategory \(\mathsf{C}\), we can associate to it a category
\(\h\mathsf{C}\) by, roughly, identifying 1-cells up to homotopy. 
Following~\cite{elgueta2007}, we call \(\h\mathsf{C}\) the \emph{homotopy
category} of \(\mathsf{C}\), similar to the \emph{na\"ive homotopy category}
found in algebraic topology. See~\cite{10.1007/BFb0074299} for details; they use
the term \emph{classifying category}, but we wish to avoid confusion with the
term \emph{classifying topos}. 

\begin{definition}[\cite{10.1007/BFb0074299}, \S7.2]\label{def:homotopy-category}
	Let \(\mathsf{C}\) be a bicategory. The \textbf{homotopy category} 
        \(\h\mathsf{C}\) of \(\mathsf{C}\) is defined by the following:
	
	(HC1) The objects of \(\h\mathsf{C}\) are the objects of \(\mathsf{C}\).
	
	(HC2) The set of morphisms \(\h\mathsf{C}(X,Y)\) between objects \(X,Y \in
	\h\mathsf{C}\) is the set of isomorphism classes of objects in the
	hom-category \(\mathsf{C}(X,Y)\). Given a 1\-/cell \(f: X \to Y\) in
	\(\mathsf{C}\), let \(\ic{f}\) denote its corresponding isomorphism class in
	\(\h\mathsf{C}(X,Y)\).
	
	(HC3) Composition is defined by the rule \(\ic{g}\ic{f} \eqdef \ic{gf}\).
	This rule is strictly associative (due to the pentagon identity), and
	\(\ic{1_X}\) is a strict identity (due to the triangle identity) under this
	composition law.
\end{definition}

\begin{definition}\label{def:homotopy-functor}
	Given a pseudofunctor \(\sfF: \mathsf{C} \to \mathsf{D}\), we obtain an
        induced functor \(\h\sfF: \h\mathsf{C} \to \h\mathsf{D}\) by setting 
        \(\h\sfF \ic{f} \eqdef \ic{ \sfF(f)}\). The compositor and identitor of 
        \(\sfF\) ensure that \(\h\sfF\) preserves composition and identities; hence
        \(\h\sfF\) is a functor.
\end{definition}

The following proposition motivates why we will use homotopy categories for
bi\-/interpretability. Its proof is an immediate application of the homotopies
present in a homotopy equivalence.

\begin{proposition}
	Two objects \(X\) and \(Y\) in a bicategory \(\mathsf{C}\) are isomorphic in
	\(\h\mathsf{C}\) if and only if they are homotopy equivalent in
	\(\mathsf{C}\).
\end{proposition}

\section{Bicategories of Theories}\label{section:bicategories-of-theories}

Unless stated otherwise, we work in the framework of \textbf{coherent logic}
(with equality), outlined in D1 of~\cite{johnstone2002sketches}. Consequently,
given a signature \(\Sigma\), the formulae we consider are the members of
\(\Sigma_{\omega\omega}^g\), i.e., those which are obtained from \(\Sigma\) by finite
combinations of \(\exists, \wedge, \vee, \top,\) and \(\bot\) in the usual
inductive fashion. There are the standard deductive laws for coherent logic
which we take for granted; these are found in, e.g., D1.3
of~\cite{johnstone2002sketches}. In particular we take for granted the existence of
equality relations \(=_\sigma\) for every sort \(\sigma\) in the signature, as
well as their usual introduction and elimination rules. We make a few minor
differences in notation. We call the sort or list of sorts corresponding to a
symbol \(A\) its \textbf{domain} (written \(\Dom A\)), where in the case of a
function symbol \(f : \vec{\sigma} \to \tau\) it additionally has a
\textbf{range} \(\tau\). Given a signature \(\Sigma\), let \(\SSort{\Sigma}^*\)
denote the free monoid generated by \(\SSort{\Sigma}\). The domain of a symbol
belonging to \(\Sigma\) is an element of \(\SSort{\Sigma}^*\). The superscript
of a quantified variable, if present, will indicate the domain of that variable.
We will also assume a bound variable has maximal scope when parentheses are
omitted.

A formula \(\phi\) with domain \(\vec{\sigma}\) will be written as \(\phi
\hookrightarrow \vec{\sigma}\). If \(\vec{x}\) is the canonical context of
\(\phi\), we will often write \(\phi(\vec{x})\) instead of \(\phi\). From this,
the formula \(\phi(\vec{y})\) is defined as \(\phi(\vec{x})[\vec{y}/\vec{x}]\).
It will be helpful to define the equivalence class of a formula \(\phi\) up to
consistent relabeling of its canonical context. We call this the
\textbf{substitution class} of \(\phi\), and denote it by \(\fb{\phi}\) (or
\(\fb{\psi}\) for any \(\psi \in \fb{\phi}\)). For example, \(\phi(\vec{y})\)
and \(\phi(\vec{x})\fb{\vec{y}/\vec{x}}\) refer to the same representative of
\(\fb\phi\).

Let \(\SForm{\Sigma}\) (resp.\ \(\SSub{\Sigma}\)) denote the set of
\(\Sigma\)-formulae (resp.\ \(\Sigma\)-substitution classes).  In order to
disambiguate the equality relations \(=_\sigma\) and identity of symbols, we
will use \(\equiv\) to denote mathematical identity. However, we will also
identify formulae and substitution classes up to \(\alpha\)\-/equivalence, i.e.,
relabeling of bound variables.  We also adopt the abbreviations \(\vec{x}
=_{\vec{\sigma}} \vec{y}\) for \(\bigwedge_{i=1}^n x_i =_{\sigma_i} y_i\) and
\(\exists \vec{x}^{\vec{\sigma}}\) for \(\exists x_1^{\sigma_1} \cdots \exists
x_n^{\sigma_n}\).  For a domain \(\sigma\), the substitution class \(\fb{x
=_\sigma x}\) is special; we refer to it as \(\fb\sigma\).

A (coherent) theory \(T\) is a pair \(\pwrap{\Sigma,\Delta}\), where \(\Sigma\)
is a signature and \(\Delta\) is a set of axioms, i.e., (coherent)
\(\Sigma\)-sequents. We say that a \(\Sigma\)-sequent is \textbf{provable in
\(T\)} if the sequent can be constructed using \(\Delta\) and the rules of
deduction. The provability relation \(\vdash\) between formulae descends to a relation
between substitution classes.

\begin{definition}
	Let \(\phi,\psi \hookrightarrow \vec{\sigma}\) be a pair of
	\(T\)-formulae. We write \(\fb\phi \vdash \fb\psi\) if
	for some context \(\vec{x}\) with domain \(\vec{\sigma}\), the sequent 
        \(\fb\phi\!(\vec{x}) \vdash \fb\psi\!(\vec{x})\) is provable in \(T\).
\end{definition}

The expression \(\phi \dashv\vdash \psi\) will stand for the sequents \(\phi \vdash
\psi\) and \(\psi \vdash \phi\) (and similarly for substitution classes), and we say
that \(\phi\) and \(\psi\) are \textbf{logically equivalent} (relative to \(T\)). This
allows us to define the \textbf{logical equivalence class} \(\ol{\fb{\phi}}\) of
a substitution class \(\fb{\phi}\) in the obvious way. We say that \(\fb{\phi}\) or 
\(\phi(\vec{x})\) (for an appropriate context) \textbf{presents} \(\ol{\fb{\phi}}\).

\subsection{Translations}

Translations, also known as interpretations, are rewritings of one theory's
statements in the language of another. By the soundness and completeness of
first-order logic, and up to the treatment of variables made implicit by Hodges
\cite{hodges_1993} and most explicit by Halvorson \cite{halvorson-logic}, our
notion of translation is the many-sorted analogue of Hodges' left-total
interpretation between theories. Hodges identifies translated symbols up to
variable substitution (see~\cite[Section 5.3, Remark 2]{hodges_1993}).
Halvorson stipulates a map on variables compatible with substitution
(see~\cite[Definition 5.4.2]{halvorson-logic}). We follow Hodges by instead
sending symbols to substitution classes. To recover a reconstrual in the sense
of Halvorson, one needs to provide only a map on variables.

\begin{definition}\label{def:reconstrual}
	Let \(\Sigma_1\) and \(\Sigma_2\) be signatures. A \textbf{reconstrual} \(F:
	\Sigma_1 \to \Sigma_2\) is a collection of the following data:

	(R1) A monoid homomorphism \(\SSort{\Sigma_1}^* \to \SSort{\Sigma_2}^*\) 
        corresponding to a function \(\SSort{\Sigma_1} \to \SSort{\Sigma_2}^*\) via 
        the universal property of free monoids.

	(R2) A function \(\SRel{\Sigma_1} \to \SSub{\Sigma_2}\) such that for any
	\(\Sigma_1\)-relation \(R \hookrightarrow \vec{\sigma}\), we have the
        compatibility condition \(FR \hookrightarrow F\vec{\sigma}\).

	(R3) A function \(\SFunc{\Sigma_1} \to \SSub{\Sigma_2}\) such that for
	any \(f: \vec{\sigma} \to \tau\) in \(\SFunc{\Sigma_1}\), we have the 
        compatibility condition \(Ff \hookrightarrow F(\vec{\sigma},\tau)\).
\end{definition}

Given a reconstrual \(F: \Sigma_1 \to \Sigma_2\), the above data allows us to
define a map \(F^+ : \SForm{\Sigma_1} \to \SSub{\Sigma_2}\) by
declaring that \(F^+\) preserves logical connectives.

\begin{mathrule}[Relations]%
	\label{rule:relations}%
	Given a \(\Sigma_1\)-formula of the form \(R(x_1,\dots,x_n)\), where \(R
	\hookrightarrow \sigma_1,\dots,\sigma_n\) is a relation and
	\(x_1,\dots,x_n\) are distinct variables, define%
	\[%
		F^+\pwrap{R(x_1,\dots,x_n)} \equivdef FR.
	\]%
\end{mathrule}

\begin{mathrule}[Preservation of \(\top\) and \(\bot\)]%
	\label{rule:top-bot-preservation}%
	Define \(F^+ \top \equivdef [\top]\) and \(F^+\bot \equivdef [\bot]\).
\end{mathrule}

\begin{mathrule}[Conjunction preservation]%
	\label{rule:conjunction-preservation}%
	Consider \(\phi \hookrightarrow \vec\sigma\) and \(\psi \hookrightarrow
	\vec\tau\) a pair of \(\Sigma_1\)-formulae, assuming \(F^+\phi\) and
	\(F^+\psi\) are defined. Let \(\vec{x}\) be a \(\vec\sigma\)-context and \(\vec{y}\) a
	\(\vec\tau\)-context such that \(\vec{x}\) and \(\vec{y}\) are disjoint. Moreover, let
	\(\vec{s}\) be an \(F\vec\sigma\)-context and \(\vec{t}\) an \(F\vec\tau\)-context such
	that \(\vec{s}\) and \(\vec{t}\) are disjoint. Define:%
	\[%
		F^+\pwrap{\phi(\vec{x}) \wedge \psi(\vec{y})} \equiv \fb{F^+\phi(\vec{s}) \wedge F^+\psi(\vec{t})}.
	\]%
\end{mathrule}

\begin{mathrule}[Disjunction preservation]%
	\label{rule:disjunction-preservation}%
	Same as Rule~\ref{rule:conjunction-preservation}, replacing \(\wedge\) with
	\(\vee\).
\end{mathrule}

\begin{mathrule}[Quantifier preservation]%
	\label{rule:quantifier-preservation}%
	Let \(\phi \hookrightarrow \vec{\sigma}\) be a \(\Sigma_1\)-formula. We can
	split the domain \(\vec\sigma\) into a list of sorts
	\(\sigma_1,\dots,\sigma_n\). Let \(x_1,\dots,x_n\) be a
	\(\vec\sigma\)-context, and let \(s_1,\dots,s_n\) be an \(F\vec\sigma\)-context
	(so \(s_i\) itself may be an \(F\sigma_i\)-context). Assume \(F^+\phi\) is
	already defined. Then define%
	\[%
		F^+\pwrap{ \exists x_i^{\sigma_i} \phi(x_1,\dots,x_n)} \equivdef
		\fb{\exists s_i^{F\sigma_i} F^+\phi(s_1,\dots,s_n)}.
	\]%
\end{mathrule}

\begin{mathrule}[Context duplication]%
	\label{rule:context-duplication}%
	Let \(\phi \hookrightarrow \vec\sigma\) be
	a \(\Sigma_1\)-formula such that \(\vec{\sigma} \equiv \vec{\tau_1},
	\vec{\sigma_1},\vec{\tau_2},\vec{\sigma_1},\vec{\tau_3}\).
	Assume \(F^+\phi\) is already defined. Then define%
	\[%
		F^+\pwrap{\phi(\vec{y_1},\vec{x_1},\vec{y_2},\vec{x_1},\vec{y_3})} \equivdef
		\fb{F^+\phi(\vec{t_1},\vec{s_1},\vec{t_2},\vec{s_1},\vec{t_3})},
	\]%
	where \(\vec{x_1}\) is a \(\vec{\sigma_1}\)-context, \(\vec{s_1}\) is an
	\(F\vec{\sigma_1}\)-context, \(\vec{y_i}\) is a \(\vec{\tau_i}\)-context, and
	\(\vec{t_i}\) is an \(F\vec{\tau_i}\)-context. We permit any, even all, of the
	\(\vec{y_i}\) to be empty.
\end{mathrule}

\begin{mathrule}[Term reduction]%
	\label{rule:term-reduction}%
	Let \(\phi \hookrightarrow \vec\tau\) be a \(\Sigma_1\)-formula, and assume
	\(F^+\phi\) is already defined. Suppose \(\vec\tau\) splits into a list
	\(\vec{\tau_1},\tau',\vec{\tau_2}\) such that \(\tau'\) is a single sort.  Let
	\(f: \vec{\sigma} \to \tau'\) be a function symbol. If \(f\) is the right-most
	function symbol appearing in the expansion of \(\phi(\vec{x_1},f(\vec{y}),\vec{x_2})\) into
	atomic formulae, set%
	\[%
		F^+\pwrap{\phi(\vec{x_1},f(\vec{y}),\vec{x_2})} \equivdef \fb{\exists {t'}^{F\tau'} \left(F f(\vec{s},t')
		\wedge F^+\phi(\vec{t_1},t',\vec{t_2})\right)},
	\]%
	where \(\vec{x_i}\) are \(\vec{\tau_i}\)-contexts, \(\vec{t_i}\) are
	\(F\vec{\tau_i}\)-contexts (with \(t'\) an \(F\tau'\)-context), \(\vec{y}\) is a
	\(\vec\sigma\)-context, and \(\vec{s}\) is an \(F\vec\sigma\)-context. As in
	Rule~\ref{rule:context-duplication}, we permit empty contexts.
\end{mathrule}

These rules are invariant under substitution of variables in a formula, and for 
any \(\Sigma_1\)\-/formula \(\phi\), \(\Dom F^+ \phi \equiv F\Dom\phi\).
Therefore \(F^+: \SForm{\Sigma_1} \to \SSub{\Sigma_2}\) descends to a map
\(\SSub{\Sigma_1} \to \SSub{\Sigma_2}\) which we call by the same name. Now that
we have defined how this map on substitution classes arises, we often remove the
superscript, writing \(F: \SSub{\Sigma_1} \to \SSub{\Sigma_2}\) for ease of
reference.

\begin{remark}\label{remark:reconstrual-functions-weirdness}
	Rule~\ref{rule:term-reduction} takes the form presented because we do not
	generally assume that reconstruals send the equality symbol to an equality
	symbol. Therefore, when we discuss translations, we will find that the image
	of a function symbol under a translation does not generally satisfy the
	axioms needed to define a function. See Chapter 4 of~\cite{halvorson-logic}
	for an extended discussion.
\end{remark}

Reconstruals can be composed (c.f.~\cite[Definition 5.4.9]{halvorson-logic}).

\begin{definition}\label{def:reconstrual-composition}
	Let \(F: \Sigma_1 \to \Sigma_2\) and \(G: \Sigma_2 \to \Sigma_3\) be a pair
	of reconstruals. The composition \(GF: \Sigma_1 \to \Sigma_3\) is the
	reconstrual defined by the following data: 
	\begin{itemize}
	    \item The monoid homomorphism \(\SSort{\Sigma_1}^* \xrightarrow{GF} 
        \SSort{\Sigma_3}^*\) corresponds to the composition \(\SSort{\Sigma_1} 
        \xrightarrow{F} \SSort{\Sigma_2}^* \xrightarrow{G} \SSort{\Sigma_3}^*\).
	    \item The function \(\SRel{\Sigma_1} \xrightarrow{GF} \SSub{\Sigma_3}\) is 
        \(\SRel{\Sigma_1} \xrightarrow{F} \SSub{\Sigma_2} \xrightarrow{G} 
        \SSub{\Sigma_3}\).
		\item The function \(\SFunc{\Sigma_1} \xrightarrow{GF} \SSub{\Sigma_3}\) is 
        \(\SFunc{\Sigma_1} \xrightarrow{F} \SSub{\Sigma_2} \xrightarrow{G}
		\SSub{\Sigma_3}\).
	\end{itemize}
\end{definition}

\begin{definition}
	Given a reconstrual \(F: \Sigma_1 \to \Sigma_2\) and \(\sigma\)-context \(x\) in
	\(\Sigma_1\), the \textbf{domain class} \(D^F_x\) is the substitution class 
        \(F\fb{x = x}\).
	We use \(D^F_\sigma\) to refer to any \(D^F_x\) for which \(x\) has domain 
        \(\sigma\). Similarly, we let \(E^F_\sigma\) denote the substitution class 
        \(F\fb{x = y}\) for any \(x,y\) with domain \(\sigma\). A formula of the form 
        \(D^F_x(t)\) (for some \(F\sigma\)-context \(t\)) has been called a
	\textbf{domain formula} in, e.g.,~\cite{halvorson-logic,hodges_1993,washington2018}.
\end{definition}

The convenience of \(D^F_\sigma\) comes from Rule \ref{rule:relations}: \(D^F_x \equiv
D^F_y\) for contexts with the same domain. Moreover, \(E^F_\sigma(x,x) \equiv 
D^F_\sigma(x)\) follows by Rules \ref{rule:relations} and \ref{rule:context-duplication}, 
and for \(\vec{\sigma} \equivdef \sigma_1,\dots,\sigma_n\), the identity
\(D^F_{\vec{\sigma}}(\vec{s}) \equiv \bigwedge_{i=1}^n D^F_{\sigma_i}(s_i)\) is due to
Rule \ref{rule:conjunction-preservation}.

A reconstrual need not make any comparison of theories regarding the provability of
sequents. A translation is a reconstrual which preserves sequents.

\begin{definition}\label{def:translation}
	Let \(T_1\) and \(T_2\) be theories. A \textbf{translation} \(F: T_1 \to T_2\)
	is a reconstrual \(F: \Sigma_1 \to \Sigma_2\) such that \(\fb\phi \vdash
	\fb\psi\) implies \(F\fb\phi \vdash F\fb\psi\). A translation is said to be
        \textbf{equality-preserving} (abbreviated to \textbf{e.p.}) if 
        \(E^F_\sigma(s,t) \vdash s =_{F\sigma} t\).
\end{definition}

The previous definition handles how sequents whose formulae have non\-/matching
contexts are translated: for a \(\Sigma_1\)-sequent \(\phi(\vec{x}) \vdash
\psi(\vec{y})\) with relative complements \(\vec{x} - \vec{y}\) and \(\vec{y} -
\vec{x}\), and for \(\vec{s}\) an \(F\Dom \vec{x}\)-context and \(\vec{t}\) an 
\(F\Dom\vec{y}\)-context such that \(s_i = t_i\) only when \(x_i = y_i\), we
have \(F\phi(\vec{s}) \wedge D^F_{\vec{y}-\vec{x}}(\vec{t}-\vec{s}) \vdash
F\psi(\vec{t}) \wedge D^F_{\vec{x}-\vec{y}}(\vec{s}-\vec{t})\), omitting any
domain formula whose subscript is empty. To show that a reconstrual is a
translation, it suffices to prove that in \(T_1\) the axioms with matching
contexts satisfy Definition \ref{def:translation} and the axioms with
non-matching contexts satisfy the previous sequent. This holds by induction and
careful application of the reconstrual rules.

There is a reconstrual \(\Sigma \to \Sigma\) which is almost an identity under composition.

\begin{definition}
	Given a signature \(\Sigma\), the \textbf{identity reconstrual} \(1_\Sigma\)
	sends a sort \(\sigma\) to itself, a relation \(R\) to \(\fb{R}\), and a 
        function symbol \(f\) to \(\fb{f(x) = y}\). For a theory \(T = 
        (\Sigma,\Delta)\), the \textbf{identity translation} \(1_T\) is the identity
	reconstrual \(1_\Sigma\).
\end{definition}

We list below some properties of reconstruals and translations that are invoked 
throughout the following section. Their proofs follow a similar approach to
those in Chapter 5 of~\cite{halvorson-logic}.

\begin{proposition}\label{prop:reconstrual-translation-properties}
	Let \(F: \Sigma_1 \to \Sigma_2\), \(G: \Sigma_2 \to \Sigma_3\) be a pair
	of reconstruals.%
	\begin{itemize}%
		\item[(i)] Composition of reconstruals is associative.
		\item[(ii)] \(1_T: T \to T\) is an e.p.\ translation.
		\item[(iii)] If \(F,G\) are translations, then \(GF\) is a translation.
		\item[(iv)] If \(F,G\) are e.p., then \(GF\) is e.p.
	\end{itemize}%
\end{proposition}

There are three additional properties of reconstruals and translations which have
not appeared in the literature; see \ref{subsection:reconstrual-properties} for their proofs.

\begin{proposition}\label{prop:e.p.-domain-formula}
If \(F\) is e.p., then \(E^F_{\sigma}(s,t) \dashv\vdash s = t \wedge D^F_\sigma(t)\).
\end{proposition}
\begin{proposition}\label{prop:reconstrual-strictly-composes}
	Let \(F: \Sigma_1 \to \Sigma_2\), \(G: \Sigma_2 \to \Sigma_3\) be a pair
	of reconstruals. For any substitution class \(\fb\phi \in \SSub{\Sigma_1}\), we
	have \(\pwrap{GF}^+ \fb\phi \equiv G^+\pwrap{F^+\fb\phi}\). Thus as functions between sets of
	substitution classes, \(\pwrap{GF}^+ = G^+ F^+\).
\end{proposition}

\begin{proposition}\label{prop:reconstrual-log-eq-to-identity}
	Consider a theory \(T = (\Sigma,\Delta)\). For any substitution class
	\(\fb\phi \in \SSub{\Sigma}\), we have that \(1_\Sigma^+ \fb\phi\) and 
        \(\fb\phi\) are logically equivalent.
\end{proposition}

Proposition~\ref{prop:reconstrual-log-eq-to-identity} will help us show
that the identity translation is a weak identity in the sense of Section 
2. This is proven in the next section.

\subsection{The Bicategorical Structure of Theories}
A t\-/map between translations is the syntactic analogue of a natural
transformation. We provide a faithful generalization to coherent logic of the
many-sorted definition provided by Halvorson~\cite[Definition
5.4.11]{halvorson-logic} for classical first\-/order logic. Prior to this,
Friedman and Visser~\cite{Friedman2014WhenBI} defined a single-sorted version
called \textit{i-maps}.

\begin{definition}\label{def:tmap}
	Let \(F,G: T_1 \to T_2\) be a pair of translations. A \textbf{t\-/map}
	\(\chi: F \Rightarrow G\) is a collection of logical equivalence classes of
	\(T_2\)-substitution classes, presented by formulae \(\chi_\sigma
	\hookrightarrow F\sigma,G\sigma\) for every \(\sigma \in \SSort{\Sigma_1}\)
	that satisfy:%
	\begin{align}
		\chi_\sigma(s,t) &\vdash D^F_\sigma(s) \wedge D^G_\sigma(t),\tag{TM1}\\
		E^F_\sigma(s,w) \wedge E^G_\sigma(t,z) \wedge \chi_\sigma(s,t) &\vdash
		\chi_\sigma(w,z),\tag{TM2}\\
		D^F_\sigma(s) &\vdash \exists t^{G\sigma} \chi_\sigma(s,t),\tag{TM3}\\
		\chi_\sigma(s,t) \wedge \chi_\sigma(s,w) &\vdash E^G_\sigma(t,w).\tag{TM4}
	\end{align}
	For a domain \(\vec\sigma \equiv \sigma_1,\ldots,\sigma_n \in
	\SSort{\Sigma_1}^*\), define \(\chi_{\vec\sigma}(\vec{s},\vec{t})\) to be the
	conjunction \(\bigwedge_{i=1}^n \chi_{\sigma_i}(s_i,t_i)\). For any
	\(T_1\)\-/formula \(\phi \hookrightarrow \vec\sigma\), we require the
	sequent%
	\begin{equation}
		F\phi(\vec{s}) \wedge \chi_{\vec\sigma}(\vec{s},\vec{t}) \vdash
		G\phi(\vec{t}).\tag{TM5}
	\end{equation}
\end{definition}

\begin{remark}
	Two t\-/maps \(\eta,\chi: F \Rightarrow G\) are equal if \(\eta_\sigma 
        \dashv\vdash \chi_\sigma\) for all \(\sigma\).
\end{remark}

For any translation \(F: T_1 \to T_2\), the collection \(\One^F_\sigma(s,t)
\equivdef E^F_\sigma(s,t)\) presents a t\-/map \(\One^F: F \Rightarrow F\).
Like reconstruals, t\-/maps can be composed. The proof of
Proposition~\ref{prop:vertical-composition-category} is elementary and formally
similar to Lemma D1.4.1 of~\cite{johnstone2002sketches}.

\begin{definition}\label{def:tmap-vertical-composition}
	Let \(\chi: F \Rightarrow G\) and \(\eta: G \Rightarrow H\) be a pair of
	t\-/maps. Define the \textbf{vertical composition} \(\eta\cdot\chi: F 
        \Rightarrow H\) to be the t\-/map presented by:%
	\[%
		\pwrap{\eta\cdot\chi}_\sigma \pwrap{s,t} \equivdef \exists w^{G\sigma}
		\pwrap{\eta_\sigma(w,t) \wedge \chi_\sigma(s,w)}.
	\]%
\end{definition}
\begin{proposition}\label{prop:vertical-composition-category}
        \(\Hom(T_1,T_2)\) is a category with t\-/maps as 1-cells, t-maps \(\One^F\)
        for \(F: T_1 \to T_2\) as identity 1-cells, and vertical composition as 
        composition.
\end{proposition}
 
As will be shown, t\-/maps correspond to the 2-cells of a bicategory. Thus, if a 
t\-/map \(\chi: F \Rightarrow G\) is invertible with respect to vertical composition,
it is a homotopy in the sense of Section~\ref{section:bicategory-theory}. This agrees
with the single-sorted definition of homotopy appearing in Ahlbrandt and 
Ziegler~\cite{AHLBRANDT198663} when restricted to first-order axiomatizable classes of
structures. In A.4 of~\cite{Friedman2014WhenBI}, Friedman and Visser axiomatize this
notion. Halvorson~\cite[Definition 5.4.12]{halvorson-logic} states a many-sorted
version, in which homotopies are t\-/maps that satisfy three additional conditions:%
\begin{align*}
	D^G_\sigma(t) &\vdash \exists s^{F\sigma} \chi_\sigma(s,t),\tag{TM6}\\
	\chi_\sigma(w,t) \wedge \chi_\sigma(z,t) &\vdash E^F_\sigma(w,z),\tag{TM7}\\
	G\phi(\vec{t}) \wedge \chi_{\vec{\sigma}}(\vec{s},\vec{t}) &\vdash F\phi(\vec{s}).\tag{TM8}
\end{align*}

In order to fit this definition for vertical composition into a bicategory of
theories, we need to define a functor \(\Hom(T_2,T_3) \times \Hom(T_1,T_2) \to
\Hom(T_1,T_3)\) for every triplet of theories \(T_1,T_2,T_3\). The object part
of this functor is composition of translations, and the morphism part will be
horizontal composition of t\-/maps. We dedicate the next subsection to
formulating horizontal composition.

\subsection{Horizontal Composition}

Here we present a novel generalization and proof of an exchange law presumed by 
Visser~\cite[Section 3.1]{Visser_categoriesof}. Whereas they define horizontal 
composition for only invertible t\-/maps over single\-/sorted theories, we extend to
arbitrary t\-/maps over many\-/sorted theories.

\begin{definition}\label{def:tmap-horizontal-composition}
	Let \(F_1,G_1: T_1 \to T_2\) and \(F_2,G_2: T_2 \to T_3\) be a quadruplet of
	translations, and let \(\chi: F_1 \Rightarrow G_1\) and \(\eta: F_2
	\Rightarrow G_2\) be a pair of t\-/maps. Define the \textbf{horizontal   
        composition} \(\eta\circ\chi\) to be the t\-/map presented by:%
	\[%
		\pwrap{\eta\circ\chi}_\sigma\pwrap{s,t} \equivdef
		\exists {t'}^{G_2 G_1\sigma}%
        \pwrap{%
        \exists w^{F_2 G_1\sigma}%
        \pwrap{%
        F_2\chi_\sigma(s,w) \wedge
        \eta_{G_1\sigma}(w,t')
		} \wedge E^{G_2 G_1}_\sigma(t',t)}.
	\]%
\end{definition}
\begin{remark}\label{remark:tmap-horizontal-composition-ep}
		If \(F_1,F_2,G_1,G_2\) are e.p., we have in \(T_3\) a logical equivalence:
		\[%
			\pwrap{\eta\circ\chi}_\sigma\pwrap{s,t} \dashv\vdash \exists w^{F_2
			G_1\sigma} \pwrap{ F_2\chi_\sigma(s,w) \wedge \eta_{G_1\sigma}(w,t) }
			\wedge D_\sigma^{G_2 G_1}(t).
		\]%
\end{remark}

\begin{proposition}\label{prop:horizontal-composition-yields-tmap}
	\(\eta\circ\chi\) is a t\-/map \(F_2 F_1 \Rightarrow G_2 G_1\).
\end{proposition}
\begin{theorem}[Exchange Law]\label{theorem:godements-law}
	Consider any collection of t\-/maps in the following shape:%
    \[%
        \begin{tikzcd}[sep=huge]
            T_1 \arrow[r,"F_1"{name=ted1},bend left=65]
            \arrow[r,"G_1"{name=karl1}] \arrow[r,"H_1"{name=mike1},bend
            right=65,swap]& T_2 \arrow[r,"F_2"{name=ted2},bend left=65]
            \arrow[r,"G_2"{name=karl2}] \arrow[r,"H_2"{name=mike2},bend
            right=65,swap] & T_3 
            \arrow[Rightarrow,from=ted1,to=karl1,"\chi_1"]
            \arrow[Rightarrow,from=karl1,to=mike1,shorten=4pt,"\eta_1"]
            \arrow[Rightarrow,from=ted2,to=karl2,"\chi_2"]
            \arrow[Rightarrow,from=karl2,to=mike2,shorten=4pt,"\eta_2"].
        \end{tikzcd}
    \]%
	Then \(\pwrap{\eta_2\cdot\chi_2}\circ\pwrap{\eta_1\cdot\chi_1} =
	\pwrap{\eta_2\circ\eta_1}\cdot\pwrap{\chi_2\circ\chi_1}\).
\end{theorem}
\begin{proposition}\label{prop:horizontal-composition-unital}
	Let \(F: T_1 \to T_2\) and \(G: T_2 \to T_3\) be a pair of translations. Then
	\(\One^G \circ \One^F = \One^{GF}\).
\end{proposition}

We prove these assertions in
~\ref{subsection:bicategory-proofs}, listing the 
relevant lemmata below.

\begin{lemma}%
	\label{lemma:t-map-lemma-1}%
	\label{lemma:t-map-domain-lemma}%
	For any pair of translations \(T_1 \xrightarrow{F_1} T_2 \xrightarrow{F_2}
	T_3\), the sequent \(D_\sigma^{F_2 F_1} \vdash D^{F_2}_{F_1\sigma}\) is
	provable in \(T_3\) for any domain \(\sigma \in \SSort{\Sigma_1}^*\).
\end{lemma}
\begin{lemma}%
	\label{lemma:t-map-lemma-2}%
	\label{lemma:t-map-equality-lemma}%
	Assume the same conditions as in Lemma~\ref{lemma:t-map-domain-lemma}. The 
        following sequent is provable in \(T_3\) for any domain \(\sigma \in 
        \SSort{\Sigma_1}^*\):%
	\[%
		E^{F_2}_{F_1\sigma}(x,y) \wedge \pwrap{ D^{F_2 F_1}_\sigma(x) \vee D^{F_2
		F_1}_\sigma(y) } \vdash E^{F_2 F_1}_\sigma(x,y).
	\]%
\end{lemma}
\begin{lemma}%
	\label{lemma:t-map-lemma-4}%
	\label{lemma:t-map-Z-lemma}%
	Let \(\chi: F_1 \Rightarrow G_1\) and \(\eta: F_2 \Rightarrow G_2\) be a pair
	of t\-/maps, where \(F_1, G_1: T_1 \to T_2\) and \(F_2, G_2: T_2 \to T_3\). Define:%
	\[%
		Z_\sigma(s,t) \equivdef
                \exists {s'}^{F_2 F_1\sigma} \pwrap{%
                E^{F_2 F_1}_\sigma(s,s') \wedge \exists w^{G_2 F_1\sigma} \pwrap{%
                \eta_{F_1\sigma}(s',w) \wedge
			G_2\chi_\sigma(w,t)
            }%
		}.%
	\]%
	Then \(Z_\sigma\) is logically equivalent to \(\pwrap{\eta\circ\chi}_\sigma\)
	for any domain \(\sigma \in \SSort{\Sigma_1}^*\).
\end{lemma}

Propositions~\ref{prop:horizontal-composition-yields-tmap}
and~\ref{prop:horizontal-composition-unital}, along with
Theorem~\ref{theorem:godements-law}, imply that horizontal composition is
functorial. To complete the construction of a bicategory of theories, we must
demonstrate unitors and associators satisfying BC7.

\begin{proposition}\label{prop:cth0-coherence}
	Composition of translations is strictly associative, i.e., \(\pwrap{F_3
	F_2} F_1 = F_3 \pwrap{F_2 F_1}\) for any composed triplet of translations.
	Hence \(a_{F_3 F_2 F_1} \equivdef \One^{F_3 F_2 F_1}\) is an invertible
	t\-/map satisfying the pentagon identity.

	Furthermore, for any translation \(F: T_1 \to T_2\), the family of
	substitution classes \(E^F_\sigma\) present invertible t\-/maps \(l_F:
	1_{T_2} F \Rightarrow F\) and \(r_F: F 1_{T_1} \Rightarrow F\) satisfying the
	triangle identities.
\end{proposition}
\begin{proof}
	Let \(F\) have domain \(\Sigma\). For any sort \(\sigma \in 
        \SSort{\Sigma}^*\), \(\pwrap{F_3 F_2} F_1
	\sigma \equiv F_3 (F_2 (F_1 \sigma) ) \equiv F_3 \pwrap{F_2 F_1} \sigma\).
	For every relation symbol \(R \in \SRel{\Sigma}\), we defined \(\pwrap{F_3 F_2}
        F_1 R \equiv \pwrap{F_3 F_2}^+ \pwrap{F_1 R}\). By 
        Proposition~\ref{prop:reconstrual-strictly-composes}, \(\pwrap{F_3
	F_2}^+ = F_3^+ F_2^+\). Moreover, \(F_1 R \equiv F_1^+
	\pwrap{R(x)}\) by Rule~\ref{rule:relations}, where \(x\) is any context 
        compatible with \(R\). Thus
	\(\pwrap{F_3 F_2} F_1 R \equiv F_3^+ F_2^+ F_1^+ \pwrap{R(x)}\). Similarly, \(
		F_3 \pwrap{F_2 F_1} R \equivdef F_3^+ \pwrap{F_2 F_1 R} \equiv F_3^+ F_2^+
		F_1^+ \pwrap{R(x)},\)
	so \(\pwrap{F_3 F_2} F_1 R \equiv F_3 \pwrap{F_2 F_1} R\). The same argument
	works for function symbols, proving that \(\pwrap{F_3
	F_2} F_1 = F_3 \pwrap{F_2 F_1}\).

	It is straightforward to check that the provided t\-/maps are invertible.
	Since composition is strictly associative, the pentagon law holds trivially.

	Since the associator is trivial, for a pair of translations \(F: T_1 \to
	T_2\) and \(G: T_2 \to T_3\), the triangle identity reduces to the equation
	\(r_G \circ \One^F = \One^G \circ l_F\). Expanding this in terms of sequents
	shows that we need to demonstrate that the following sequent is provable in
	\(T_3\).%
	\begin{align*}%
		\exists w E^{G 1_{T_2} F}_\sigma(s,w) \wedge E^G_{F\sigma}(w,t) \wedge
		&D^{GF}_\sigma(t)\\%
		&\dashv\vdash \exists w' E^{GF}_\sigma(s,w') \wedge E^G_{F\sigma}(w',t)
		\wedge D^{GF}_\sigma(t)
	\end{align*}%
	Recall \(1_{T_2} (F\sigma) \equiv F\sigma\), so \(w\) and \(w'\) have the
	same domain. Therefore it suffices to prove \(E^{G 1_{T_2} F}_\sigma(s,w)
	\dashv\vdash E^{GF}_\sigma(s,w)\). This follows from
	Proposition~\ref{prop:reconstrual-log-eq-to-identity} and the fact that \(G\)
	is a translation.
\end{proof}

Having verified the coherence laws, we obtain a bicategory of theories.

\begin{theorem}\label{theorem:cth0}
	The collection of small coherent theories forms a bicategory \(\CTh_0\), where
	the 1\-/cells are translations, and the 2\-/cells are t\-/maps.	 The
	composition laws are horizontal and vertical composition of t\-/maps. The
	associator is trivial, and the unitors are described in
	Proposition~\ref{prop:cth0-coherence}.
\end{theorem}

\begin{definition}\label{def:ctheq}
	Let \(\CThEq\) be the 2\-/full sub-bicategory of \(\CTh_0\) spanned by
	equality\-/preserving translations.
\end{definition}

\begin{remark}\label{remark:WhyWeCantGetStrict2Eq}
	Unlike the associator, the unitors of \(\CTh_0\) and \(\CThEq\) are
	nontrivial due to the fact that the identity reconstrual is only a weak
	identity: reconstruals are stipulated to send function symbols to
	substitution classes, so \(F : T_1 \to T_2\) might not behave identically to
	\(1_{T_2} F\) or \(F 1_{T_1}\). Nevertheless, they behave in a logically
	equivalent way, which allows us to obtain unitors.
\end{remark}

Bi-interpretability of theories is historically supported~\cite{AHLBRANDT198663,
Button-Walsh2018,enayat_schmerl_visser_2011,halvorson-logic,hodges_1993,Marker2002,
mceldowney_2020,roque_freire_hamkins_2021,Visser_categoriesof,walsh_2014}. 
See also Button and Walsh~\cite{Button-Walsh2018} for motivation.

\begin{definition}\label{def:bi-interpretability}
	Bi\-/interpretability is homotopy equivalence in \(\CTh_0\). That is, two
	theories \(T_1\) and \(T_2\) are \textbf{bi\-/interpretable} if there exist
	translations \(F: T_1 \to T_2\) and \(G: T_2 \to T_1\) such that \(GF \simeq
	1_{T_1}\) and \(FG \simeq 1_{T_2}\). If, further, both \(F\) and \(G\) are e.p., 
        we say that \(T_1\) and \(T_2\) are \textbf{e.p.~bi\-/interpretable}.
\end{definition}

\section{Biequivalence}

The treatment of syntactic categories in~\cite{MakkaiReyes1977} suggest that the
act of sending a coherent theory \(T\) to its syntactic category may be
considered the object part of a functor (see Proposition 8.1.1
of~\cite{MakkaiReyes1977}). We make this precise by extending the syntactic
category and internal logic relations in~\cite{MakkaiReyes1977} and D1.4
of~\cite{johnstone2002sketches} into pseudofunctors \(\C: \CThEq \to \Coh\) and
\(\T: \Coh \to \CThEq\). 

Since this section works exclusively with 
e.p.\ translations, horizontal composition will assume the simpler presentation 
specified by Remark \ref{remark:tmap-horizontal-composition-ep}.

\subsection{Defining the Pseudofunctors}
We begin by reviewing the definition of syntactic category found in D1.4
of~\cite{johnstone2002sketches}, with one modification for convenience.
Whereas the objects of our syntactic category \(\C(T)\) are substitution classes
\(\fb{\phi(x)}\) of formulae, the objects of Johnstone's \(\mathcal{C}_T\) are
substitution classes \(\bwrap{x.\phi}\) of \emph{formulae in context}---the
context of \(\bwrap{x.\phi}\) may be larger than the domain of \(\phi\). This
distinction is insignificant: the extra free variables of \(\bwrap{x.\phi}\) can
be absorbed by replacing \(\phi\) with the logically equivalent \(\phi \wedge x
= x\). In particular, the map \(\bwrap{x.\phi} \mapsto \fb{\phi \wedge x = x}\) 
defines an equivalence of categories \(\mathcal{C}_T \to \C(T)\).

\begin{definition}\label{def:syntactic-category-0-cell}
For a small coherent theory \(T\), the \textbf{syntatic category} \(\C(T)\) 
is the small category whose objects are \(T\)-substitution classes
\(\fb\phi\) and whose morphisms \(\theta: \fb\phi \to \fb\psi\) are what we
call \(T\)-\textbf{definable maps}: logical equivalence classes of 
\(T\)-substitution classes \(\ol{\fb{\theta}}\) such that any choice 
\(\fb{\theta}\) must satisfy:%
\begin{align}
	\theta(x,y) &\vdash \phi(x) \wedge \psi(y),\tag{DM1}\\
	\theta(x,y_1) \wedge \theta(x,y_2) &\vdash y_1 = y_2,\tag{DM2}\\
	\phi(x) &\vdash \exists y \, \theta(x,y).\tag{DM3}
\end{align}%
\end{definition}

Given \(T\)-definable maps \(\ol{\fb{\alpha}} : \fb\phi \to \fb\psi\)
and \(\ol{\fb{\beta}}: \fb\psi \to \fb\eta\), we can compose them to obtain a
\(T\)-definable map \(\ol{\fb{\beta\alpha}} : \fb\phi \to \fb\eta\). It is presented by
the formula%
\[%
        \beta\alpha(x,z) \equivdef \exists y \pwrap{\alpha(x,y) \wedge \beta(y,z)},
\]%
where \(\fb{\alpha}\) and \(\fb{\beta}\) are substitution classes presenting
\(\ol{\fb{\alpha}}\) and \(\ol{\fb{\beta}}\), respectively. This composition law is
associative, with identity \(1_{\fb\phi}: \fb\phi \to \fb\phi\)
given by \(1_{\fb\phi}(x,y) \equivdef x = y \wedge \phi(x)\) (see Lemma D1.4.1
of~\cite{johnstone2002sketches}).
By Lemma D1.4.10 of~\cite{johnstone2002sketches}, \(\C(T)\) is a coherent category 
when \(T\) is a coherent theory. The (co)limits of \(\C(T)\) can be characterized by
sequents of \(T\); see~\ref{subsection:sc-schemata} for future reference. 

We now establish the 1- and 2-cell components of a
pseudofunctor \(\C\) from \(\CThEq\) to \(\Coh\), proving well\-/definedness in
the next subsection. Proofs of the coherence laws, PF3 through PF5, appear
in~\ref{subsection:coherence-proofs}.

\begin{definition}\label{def:syntactic-category-1-cell}
	Let \(F: T_1 \to T_2\) be an e.p.\ translation. Define a map \(\C(F): \C(T_1)
	\to \C(T_2)\) as follows. For every object \(\fb\phi\) of \(\C(T_1)\), let 
        \(\C(F)[\phi]\) be the substitution class \(F^+[\phi]\). For
	every morphism \(\theta : \fb\phi \to \fb\psi\) of \(\C(T_1)\), pick a
	representative substitution class \(\fb{\theta}\) for \(\ol{\fb{\theta}}\) 
        and define \(\C(F)\theta\) to be the \(T_2\)\-/definable map presented by
	\(F^+[\theta]\).
\end{definition}

\begin{definition}\label{def:syntactic-category-2-cell}
	Let \(\chi: F \Rightarrow G\) be a t\-/map between e.p.\ translations \(F,G:
	T_1 \to T_2\). Define \(\C(\chi)\) to be the map \(\C(F)
	\Rightarrow \C(G)\) whose component along an object \(\fb\phi \hookrightarrow 
        \sigma\) of \(\C(T_1)\) is the \(T_2\)\-/definable map 
        \(\ol{\fb{\C(\chi)_{\fb\phi}}} : \C(F)[\phi] \to \C(G)[\phi]\) presented by 
        the substitution class \(\fb{\chi_\sigma(s,t) \wedge F \phi(s)}\) (picking a
	representative substitution class for each \(\chi_\sigma\)).
\end{definition}

\begin{proposition}[PF3\((\C)\)]\label{prop:C-trivial-compositor}
	Let \(F: T_1 \to T_2\) and \(G: T_2 \to T_3\) be a pair of
	e.p.\ translations. The two functors \(\C(GF)\) and
	\(\C(G)\C(F)\) from \(\C(T_1)\) to \(\C(T_3)\) are equal, and \(\C\) has a
	trivial compositor.
\end{proposition}

\begin{proposition}[PF4\((\C)\)]\label{prop:C-identitor}
	Let \(T\) be a coherent theory. Given an object \(\fb\phi\) of \(\C(T)\), the
	substitution class \(\fb{ \phi(x) \wedge x = y}\) presents a morphism
	\(\fb\phi \to \C(1_T)\fb\phi\). This morphism forms the \(\fb\phi\) component
	of a natural ismorphism \(\C_{1_T}: 1_{\C(T)} \Rightarrow \C(1_T)\), making
	the identitor of \(\C\).
\end{proposition}

\begin{proposition}\label{prop:C-pseudofunctor}
	The maps \(T \mapsto \C(T)\), \(F \mapsto \C(F)\), and \(\chi \mapsto \C(\chi)\)
	define a pseudofunctor \(\C: \CThEq \to \Coh\) called the \textbf{syntactic
	category} pseudofunctor. It has a trivial compositor and an identitor defined
	in Proposition~\ref{prop:C-identitor}.
\end{proposition}

The internal logic operation \(\T: \Coh \to \CThEq\) will be the pseudoinverse
of \(\C: \CThEq \to \Coh\). The 0-cell component comes
from~\cite[Chapter 2]{MakkaiReyes1977}.
\begin{definition}\label{def:internal-logic-0-cell}
	Let \(C\) be a coherent category. Let \(\ul\Sigma_C\) be the signature
	constructed by adding a sort \(\ul A\) and a binary relation \(=_{\ul A}
	\hookrightarrow \ul A, \ul A\) for every object \(A\) of \(C\), as well as a
	function symbol \(\ul f: \ul A \to \ul B\) for every morphism \(f: A \to B\)
	in \(C\). If \(C\) is a small category, \(\ul \Sigma_C\) is a small set. Let
	\(\ul \Delta_C\) be the set of sequents in the signature \(\ul \Sigma_C\)
	defined by the IL axiom schemata 1-10 in~\ref{subsection:IL-axiom-schemata}.
	The \textbf{internal theory} \(\T(C)\) of \(C\) is the theory \(\pwrap{\ul
	\Sigma_C, \ul \Delta_C}\).
\end{definition}

Our formalization of translations allows us to provide a sensible extension of
\(C \mapsto \T(C)\) into a pseudofunctor.

\begin{definition}\label{def:T-1-cells}
	Let \(\Eff: C_1 \to C_2\) be a coherent functor. Define a reconstrual
	\(\T(\Eff): \ul\Sigma_{C_1} \to \ul\Sigma_{C_2}\) in the following manner. For a
	sort \(\ul A\) of \(\ul\Sigma_{C_1}\), set \(\T(\Eff)\ul A\) to be the sort
	\(\ul{\Eff A}\) in \(\ul\Sigma_{C_2}\). For a function symbol \(\ul f: \ul A
	\to \ul B\), set \(\T(\Eff)\ul f\) to be the substitution class
	\(\fb{ \ul{\Eff f}(x) = y }\). For an equality relation \(=_{\ul A}\), define
	\(\T(\Eff)=_{\ul A}\) to be the substitution class \(\fb{ x =_{\ul{\Eff A}}
	y}\). Then the reconstrual \(\T(\Eff)\) is an e.p.\ translation \(\T(C_1) \to \T(
	C_2)\).
\end{definition}

\begin{definition}\label{def:T-2-cells}
	Let \(\chi: \Eff \Rightarrow \Gee\) be a natural transformation between
	coherent functors \(\Eff,\Gee: C_1 \to C_2\). We define a t\-/map \(\T(\chi):
	\T(\Eff) \Rightarrow \T(\Gee)\) as follows. For an object \(A\) of \(C_1\), the
	component of \(\chi\) along \(A\) is a morphism \(\chi_A: \Eff A \to \Gee A\).
	This morphism picks out a function symbol \(\ul{\chi_A}: \ul{\Eff A} \to
	\ul{\Gee A}\) in \(\T(C_2)\). Define the component of \(\T(\chi)\) along the sort
	\(\ul{A}\) to be the logical equivalence class of \(\fb{ \ul{\chi_A}(x) = y }\).
\end{definition}

\begin{proposition}[PF3\((\T)\)]\label{prop:T-compositor}
	Let \(\Eff: C_1 \to C_2\) and \(\Gee: C_2 \to C_3\) be a pair of coherent
	functors. The substitution classes \(\fb{ x =_{\ul{\Gee\Eff A}} y }\) for each
	sort \(\ul A\) of \(\T(C_1)\) present a (t\-/map) homotopy \(\kappa_{\Gee\Eff}:
	\T(\Gee)\T(\Eff) \Rightarrow \T(\Gee\Eff)\). These homotopies form the
	components of the compositor of \(\T\).
\end{proposition}

\begin{proposition}[PF4\((\T)\)]\label{prop:T-identitor}
	Let \(C\) be a coherent category. The translations \(\T(1_C)\)
	and \(1_{\T(C)}\) are identical. Therefore the identitor of \(\T\) is trivial.
\end{proposition}

\begin{proposition}\label{prop:T-pseudofunctor}
	The maps \(C \mapsto \T(C)\), \(\Eff \mapsto \T(\Eff)\), and \(\chi \mapsto
	\T(\chi)\) define a pseudofunctor \(\T: \Coh \to \CThEq\) called the
	\textbf{internal logic} pseudofunctor. It has a compositor defined in
	Proposition~\ref{prop:T-compositor} and a trivial identitor.
\end{proposition}

\subsection{Well-Definedness and Coherence}

We perform the necessary checks to ensure \(\C\) and \(\T\) are pseudofunctors.
First, we show that \(\C\) is well\-/defined.

\begin{proposition}[Well-Definedness of \(\C\)]\label{prop:C-well-defined}
	Let \(F: T_1 \to T_2\) be an e.p.\ translation. Then \(\C(F)\) is a coherent
	functor from \(\C(T_1)\) to \(\C(T_2)\). Let \(\chi: F \Rightarrow G\) be a
	t\-/map between e.p.\ translations. Then \(\C(\chi)\) is a natural
	transformation from \(\C(F)\) to \(\C(G)\). Furthermore, given another t\-/map
	\(\eta\) between e.p.\ translations, \(\C(\eta\cdot\chi) = \C(\eta)\cdot\C(\chi)\)
	whenever \(\eta\cdot\chi\) is defined. Lastly, \(\C(\One^F) = \One^{\C(F)}\).
\end{proposition}

\begin{proof}
	We first show that \(\C(F)\) is a functor. Let \(\alpha: \fb\phi \to \fb\psi\) and
	\(\beta: \fb\psi \to \fb\eta\) be a pair of morphisms in \(\C(T_1)\).
	Applying Rules~\ref{rule:conjunction-preservation} 
        and~\ref{rule:quantifier-preservation}, we have
	\(\C(F)\pwrap{\beta\circ\alpha}(x,y) \equiv \exists w \pwrap{F\alpha(x,w)
	\wedge F\beta(w,y)} \equiv \pwrap{\C(F)\beta
	\circ \C(F)\alpha}(x,y)\). Thus \(\C(F)\) preserves composition of
	morphisms in the syntactic category. Since \(F\) is e.p., \(\C(F)
	1_{\fb{\phi}} = 1_{\C(F)\fb{\phi}}\), so \(\C(F)\) preserves identities;
	hence \(\C(F)\) is a functor.

	To see that \(\C(F)\) is coherent, we first show that
	\(\C(F)\) preserves finite limits. It suffices to prove that pullbacks and
	terminal objects are preserved. The terminal object in \(\C(T)\) is
	isomorphic to \(\fb\top\). By Rule~\ref{rule:top-bot-preservation},
	\(\C(F)\fb\top \equiv \fb\top\), so since \(\C(F)\) is a functor, it must
	preserve terminal objects. Given a cospan \(\fb{\phi_1} \xrightarrow{\theta_1}
        \fb{\psi} \xleftarrow{\theta_2} \fb{\phi_2}\) in \(\C(T)\),
	any associated pullback square is isomorphic to the following square, where
	\(\left(\theta_1\wedge\theta_2\right)(x_1,x_2)\) is defined as \(\exists y
	\left( \theta_1(x_1,y) \wedge \theta_2(x_2,y) \right)\), and
	\(p_i(x_1,x_2,x)\) is defined as \(\pwrap{\theta_1\wedge\theta_2}(x_1,x_2)
	\wedge x_i = x\).%
	\[%
		\begin{tikzcd}
			\pwrap{\theta_1 \wedge \theta_2} \arrow[d,"p_1", swap] \arrow[r,"p_2"] &
			\fb{\phi_2} \arrow[d,"\theta_2"] \\
			\fb{\phi_1} \arrow[r,"\theta_1", swap] & \fb\psi
		\end{tikzcd}
	\]%
	Since \(F\) is e.p., we can invoke Rules~\ref{rule:conjunction-preservation}
	and~\ref{rule:quantifier-preservation} to see that applying \(\C(F)\) to
	every vertex and edge of this square yields a pullback square in \(\C(T_2)\).
	Thus \(\C(F)\) preserves pullbacks and terminal objects, so it preserves
	finite limits.

	To show that \(\C(F)\) preserves finite joins and images, we use the same
	idea as the previous paragraph, invoking
	Rules~\ref{rule:relations}-\ref{rule:term-reduction} and applying \(F\) to
	key diagrams. For a finite family of monics \(\theta_\alpha: \fb{\phi_\alpha}
	\hookrightarrow \fb{\psi}\), their join is given by the monic
	\(\fb{\bigvee_\alpha \exists x_\alpha \theta_\alpha(x_\alpha,y)}
	\hookrightarrow \fb\psi\). We use Rules~\ref{rule:disjunction-preservation}
	and~\ref{rule:quantifier-preservation} to infer that \(\C(F)\) preserves
	finite joins. Lastly, for a monic \(\theta: \fb\phi \hookrightarrow
	\fb\psi\), its image under a morphism \(f: \fb{\psi} \to \fb{\psi'}\) is
	presented by the monic \(\fb{\exists x \exists w \pwrap{\theta(w,x) \wedge
	f(x,x')}} \hookrightarrow \fb{\psi'}\). Therefore
	Rules~\ref{rule:conjunction-preservation} and
	\ref{rule:quantifier-preservation} imply that \(\C(F)\) preserves images. This
	completes the proof that \(\C(F)\) is a coherent functor when \(F\) is an
	e.p.\ translation.

	Now we show that \(\C(\chi)\) is a natural transformation \(\C(F) \Rightarrow
	\C(G)\). Let \(\theta: \fb\phi \to \fb\psi\) be an arbitrary morphism in
	\(\C(T_1)\), where \(\Dom \phi \equivdef \sigma\) and
	\(\Dom \psi \equivdef \tau\). Unpacking the definition of \(\C(F)\) on objects and
	morphisms, we need to show that the following diagram commutes.%
	\[%
		\begin{tikzcd}
			F\fb\phi \arrow[r,"F\theta"] \arrow[d,"\C(\chi)_{\fb\phi}", swap] &
			F\fb\psi \arrow[d,"\C(\chi)_{\fb\psi}"] \\
			G\fb\phi \arrow[r,"G\theta"] & G\fb\psi
		\end{tikzcd}
	\]%
	Since morphisms in a syntactic category with matching domain and codomain are
	equal whenever they are presented by logically\-/equivalent substitution
	classes, it suffices to prove the following sequents in \(T_2\).%
	\[%
		\exists y^{F\tau} F\theta(x,y) \wedge \chi_\tau(y,z) \wedge F\psi(y)
		\dashv\vdash \exists {y'}^{G\sigma} \chi_\sigma(x,y') \wedge F\phi(x)
		\wedge G\theta(y',z)
	\]%
	We begin with the forward sequent. Note \(F\theta(x,y) \vdash D^F_\sigma(x)\)
	since \(F\) is a translation. Therefore TM3(\(\chi\)) implies (by the cut
	rule) \(F\theta(x,y) \vdash \exists {y'}^{G\sigma} \chi_\sigma(x,y')\). By
	TM5(\(\chi\)) we have \(F\theta(x,y) \wedge \chi_\sigma(x,y') \wedge
	\chi_\tau(y,z) \vdash G\theta(y',z)\). From this we infer the forward
	sequent. We proceed to the converse sequent. Since \(F\) is a translation,
	DM3(\(\theta\)) implies \(F\phi(x) \vdash \exists y^{F\tau} F\theta(x,y)\).
	Applying DM1(\(\theta\)) yields \(F\theta(x,y) \vdash F\psi(y)\), and
	\(F\psi(y) \vdash D^F_\tau(y)\). Now we apply TM3(\(\chi\)) to get \(F\phi(x)
	\vdash \exists {z'}^{G\tau} \exists y^{F\tau} F\theta(x,y) \wedge
	\chi_\sigma(x,y') \wedge \chi_\tau(y,z')\). TM5(\(\chi\)) gives us \(\chi_\sigma(x,y')
	\wedge \chi_\tau(y,z') \wedge F\theta(x,y) \vdash G\theta(y',z')\). Thus we
	can use DM2(\(G\theta\)) to replace \(z'\) with \(z\): \(\chi_\sigma(x,y')
	\wedge F\phi(x) \wedge G\theta(y',z) \vdash \exists y^{F\tau} F\theta(x,y)
	\wedge \chi_\tau(y,z)\). DM1(\(F\theta\)) gives us \(F\theta(x,y) \vdash
	F\psi(y)\), completing the proof of the converse sequent.

	All that is left is proving \(\C(\eta\cdot\chi) = \C(\eta)\cdot\C(\chi)\) and
	\(\C(\One^F) = \One^{\C(F)}\). Since \(\C(\eta\cdot\chi)\) and
	\(\C(\eta)\cdot\C(\chi)\) have the same source and target functors, the first
	equation reduces to verifying that the components of \(\C(\eta\cdot\chi)\)
	and \(\C(\eta)\cdot\C(\chi)\) have logically equivalent presentations. This
	is a consequence of TM5 and the observation that vertical composition of
	t\-/maps and composition of morphisms in the syntactic category have the same
	form, syntactically speaking. Similarly, the second equation requires showing
	that \(\C(\One^F)_{\fb\phi} \dashv\vdash \One^{\C(F)}_{\fb\phi}\). Recall
	that \(\One^F : F \Rightarrow F\) is the t\-/map defined by
	\(\One^F_\sigma(x,y) \equivdef E^F_\sigma(x,y) \dashv\vdash x =_{F\sigma} y
	\wedge D^F_\sigma(x)\). Let \(\Dom \phi \equiv \sigma\). Then the component
	of \(\C(\One^F)\) along \(\fb\sigma\) is the definable map presented by \(x
	=_{F\sigma} y \wedge D^F_\sigma(x) \wedge F\phi(x)\). Since \(F\phi(x) \vdash
	D^F_\sigma(x)\), this is logically equivalent (in \(T_2\)) to \(x =_{F\sigma}
	y \wedge F\phi(x)\). This is the definition of the identity \(1_{F\fb\phi}\),
	which is the component of \(\One^{\C(F)}\) along \(\fb\phi\). Hence
	\(\C(\One^F)_{\fb\phi} \dashv\vdash \One^{\C(F)}_{\fb\phi}\).
\end{proof}

The preceding proposition shows that the maps defining \(\C\) are
well\-/defined. The proofs that these maps satisfy the appropriate coherence laws
are found in~\ref{subsection:coherence-proofs}. We now show \(\T\) is a
pseudofunctor, beginning with well-definedness.

\begin{proposition}[Well-Definedness of \(\T\)]\label{prop:T-well-defined}
	Let \(\Eff: C_1 \to C_2\) be a coherent functor. Then \(\T(\Eff): \T(C_1) \to
	\T(C_2)\) is an e.p.\ translation. Let \(\chi: \Eff \Rightarrow \Gee\) be a
	natural transformation between coherent functors. Then \(\T(\chi)\) is a
	t\-/map from \(\T(\Eff)\) to \(\T(\Gee)\). Furthermore, given another natural
	transformation \(\eta\) between coherent functors, \(\T(\eta\cdot\chi) =
	\T(\eta) \cdot \T(\chi)\), whenever \(\eta\cdot\chi\) is defined, as well as
	\(\T(\One^\Eff) = \One^{\T(\Eff)}\). Lastly, \(\T(1_C) = 1_{\T(C)}\) for any
	coherent category \(C\).
\end{proposition}

\begin{proof}
	Abbreviate \(\T(\Eff)\) to \(F\) and \(\T(\Gee)\) to \(G\). We first show
	    that \(F\) is an e.p.\ translation. Recall we defined \(E^F_{\ul
	A}(x,y)\) to be \(x =_{\ul{\Eff A}} y\) in Definition~\ref{def:T-1-cells}, so
	if \(F\) is a translation, then it is e.p.\ as well. To show that \(F\) is a
	translation, we need to show that the images of the IL axiom schemata for
	\(\T(C_1)\) under \(F\) are provable in \(\T(C_2)\). This is true because
	\(\Eff\) is a coherent functor, so it preserves the (co)limits
	mentioned in IL1 through IL10.\footnote{%
		For example, consider the IL3 axiom for a monic \(f: A \hookrightarrow B\)
		in \(C_1\). IL3(\(f\)) is the sequent \(\ul{f}(x_1') =_{\ul B} \ul{f}(x_2')
		\vdash x_1' =_{\ul A} x_2'\). Under \(F\) this sequent translates to%
		\[%
			\exists y_1^{\ul{\Eff B}} \exists y_2^{\ul{\Eff B}} \pwrap{%
				y_1 =_{\ul{\Eff B}} y_2 \wedge F\ul{f}(x_1,y_1) \wedge F\ul{f}(x_2,y_2)
			} \vdash x_1 =_{\ul{\Eff A}} x_2,
		\]%
		where we applied Rule~\ref{rule:term-reduction} twice on the left. From the
		definition of \(F\), we know that the left side is logically equivalent to
		\(\ul{\Eff f}(x_1) =_{\ul{\Eff B}} \ul{\Eff f}(x_2)\), so the translated
		sequent is logically equivalent to the IL3 axiom for \(\ul{\Eff f}\). The
		rest of the IL axioms of \(C_1\) follow a similar argument.
	} %
	Thus \(F\) is an e.p.\ translation.

	We now show that \(\T(\chi)\) is a t\-/map \(F \Rightarrow G\). Since
	\(\ul{\chi_A}\) is a function symbol, TM1(\(\T(\chi)\))	through TM4(\(\T(\chi)\))
	are provable in \(\T(C_2)\). All that remains is TM5(\(\T(\chi)\)): for any 
        \(\T(C_1)\)-formula \(\phi\) with domain \(\vec{\ul{A}} \equivdef
	\ul{A_1},\ldots,\ul{A_n}\), we need%
	\[%
		\T(\chi)_{\ul{\vec{A}}}(\vec{x},\vec{y}) \wedge F\phi(\vec{x}) \vdash
		G\phi(\vec{y}).
	\]%
	It suffices to consider the case where \(\phi\) is an atomic
	formula.\footnote{%
		This is because \(F\) and \(G\) preserve logical connectives, so
		TM5(\(\T(\chi)\)) can be proven for an arbitrary formula by breaking it down
		into a family of sequents of atomic formulae using the introduction and
		elimination rules for \(\wedge, \vee,\) and \(\exists\).
	}
	Atomic formulae in \(\T(C_1)\) take the
	form \(\phi(x_1', x_2') \equiv t_1(x_1') =_{\ul B} t_2(x_2')\) for some pair
	of terms \(t_1: \ul{A_1} \to \ul{B}\) and \(t_2: \ul{A_2} \to \ul{B}\). We
	may assume without loss of generality that \(t_1\) and \(t_2\) are function
	symbols \(\ul{f}\) and \(\ul{g}\), respectively. This is because we may apply IL2
	axioms to reduce a composition of function symbols into a single function
	symbol. Using Rule~\ref{rule:term-reduction} we see that \(F\phi(x_1,x_2)\)
	and \(G\phi(y_1,y_2)\) are logically equivalent to \(\ul{\Eff f}(x_1)
	=_{\ul{\Eff B}} \ul{\Eff g}(x_2)\) and \(\ul{\Gee f}(y_1) =_{\ul{\Gee B}}
	\ul{\Gee g}(y_2)\), respectively. In this case TM5 becomes \(
	\T(\chi)_{\ul{A_1}}(x_1,y_1) \wedge \T(\chi)_{\ul{A_2}}(x_2,y_2) \wedge \ul{\Eff
	f}(x_1) = \ul{\Eff g}(x_2) \vdash \ul{\Gee f}(y_1) = \ul{\Gee g}(y_2). \) We
	can replace \(\T(\chi)_{\ul{A_i}}(x_i,y_i)\) with its definition and apply
	\(=\)-elimination to reduce TM5 to the sequent \(\ul{\Eff f}(x_1) = \ul{\Eff
		g}(x_2) \vdash \ul{\Gee f}(\ul{\chi_{A_1}}(x_1)) = \ul{\Gee
		g}(\ul{\chi_{A_2}}(x_2))\). Deduction yields the sequent \(\ul{\Eff
	f}(x_1) = \ul{\Eff g}(x_2) \vdash \ul{\chi_B}(\ul{\Eff f}(x_1)) =
	\ul{\chi_B}(\ul{\Eff g}(x_2))\). Since \(\chi\) is a natural transformation,
	we can use IL2 axioms to derive the sequents \(\vdash \ul{\chi_B}(\ul{\Eff
	f}(x_1)) = \ul{\Gee f}(\ul{\chi_{A_1}}(x_1))\) and \(\vdash
	\ul{\chi_B}(\ul{\Eff g}(x_2)) = \ul{\Gee g}(\ul{\chi_{A_2}}(x_2))\) in \(\T(
	C_2)\), whence we apply the previous sequent to derive TM5.

	Suppose \(\chi\) is an identity natural transformation \(\One^\Eff: \Eff
	\Rightarrow \Eff\), so its component along an object \(A\) is the identity
	morphism \(1_{\Eff A}\). We can use IL1(\(1_{\Eff A}\)) to see that
	\(\T(\chi)_{\ul{A}}(x,y)\) is logically equivalent to \(x =_{\ul{\Eff A}}
	y\), which presents the identity t\-/map \(\One^F\), since \(F\) is e.p.
	Therefore \(\T(\One^\Eff) = \One^F\).

	Lastly, we show that \(\T(\eta\cdot\chi) = \T(\eta)\cdot\T(\chi)\) and 
        \(\T(1_C) = 1_{\T(C)}\). For the first equation, the component of
	\(\eta\cdot\chi\) along an object \(A\) is \(\eta_A \chi_A\), so the collection
	of IL2 axioms associated to \(\eta_A, \chi_A,\) and \(\eta_A\chi_A\) for
	every \(A\) imply the result. For the second, when \(\Eff\) is the identity functor 
        \(1_C: C \to C\), Definition~\ref{def:T-1-cells} shows that the underlying reconstrual of 
	\(\T(1_C)\) is the identity reconstrual. We conclude.
\end{proof}

As for \(\C\), the proofs of the coherence laws for \(\T\) are found
in~\ref{subsection:coherence-proofs}.

\subsection{Pseudonautral Homotopies and Biequivalence}

We establish the second main theorem of this paper.%
\begin{theorem:CThEqCoh-Correspondence}
	The pseudofunctors \(\C: \CThEq \to \Coh\) and \(\T: \Coh \to \CThEq\) form a
	biequivalence.
\end{theorem:CThEqCoh-Correspondence}%
Before the proof we review some related theory from~\cite{MakkaiReyes1977} 
and~\cite{johnstone2002sketches}. For an object \(\fb\phi\) of \(\C(T)\), let
\(\dom_{\fb\phi}\) denote the morphism \(\fb\phi \to \fb{\Dom\phi}\) presented
by \(\dom_{\fb\phi}(x,y) \equivdef \phi(x) \wedge x = y\). Lemma 1.4.4(iii)
of~\cite{johnstone2002sketches} shows that \(\dom_{\fb\phi}\) is monic. 
In this sense subobjects of \(\C(T)\) generalize the notion of
domain defined earlier for theories: \(\fb\phi \in \Sub\fb{\Dom\phi}\). 

\begin{lemma}[\cite{johnstone2002sketches}, Lemma 1.4.4(iv)]\label{lemma:subobject-lemma}
	Let \(\phi(x),\psi(x) \hookrightarrow \sigma\) be a pair of formulae in \(T\).
	Then \(\fb\phi \vdash \fb\psi\) in \(T\) if and only if \(\dom_{\fb\phi} \leq
	\dom_{\fb\psi}\) as subobjects of \(\fb\sigma\) in \(\C(T)\).
\end{lemma}%
Now we define the pseudonatural homotopies \(\ep: \id_\CThEq \Rightarrow \T\C\)
and \(\delta: \id_\Coh \Rightarrow \C\T\), beginning with \(\ep\).

\subsubsection{The Pseudonatural Homotopy \(\ep\)}\label{subsub:epsilon}

The product \(\fb{\sigma_1} \times \ldots \times
\fb{\sigma_n}\) is presented by the conjunction \(\fb{\bigwedge_{i=1}^n
x_i^{\sigma_i} = x^{\sigma_i}_i} \equiv \fb{\vec x = \vec x} \equiv
\fb{\vec\sigma}\). When we need an explicit presentation of a product, we will
use this conjunction. For example, in \(\T\C(T)\) for some theory \(T\),
\(\ul{\fb{\sigma_1} \times \ldots \times \fb{\sigma_n}} \equivdef
\ul{\fb{\vec\sigma}}\). Given a domain \(\sigma_1,\ldots,\sigma_n\) in \(T\), we
have projection morphisms \(\pi^{\vec\sigma}_{\sigma_i}: \fb{\vec\sigma} \to
\fb{\sigma_i}\); these are presented by \(\pi^{\vec\sigma}_{\sigma_i}(\vec x, y)
\equivdef \vec x = \vec x \wedge x_i = y\). Lastly, for a function symbol
\(f: \vec\sigma \to \tau\) in \(T\), let \(\theta_f : \fb{\vec\sigma} \to
\fb\tau\) denote the morphism presented by \(\theta_f(x,y) \equivdef f(x) =
y\).

\begin{definition}\label{def:ep_T}
	Let \(T = \pwrap{\Sigma,\Delta}\) be a coherent theory. Define a reconstrual
	\(\ep_T: T \to \T\C(T)\) in the following manner.  For a sort \(\sigma \in
	\SSort{\Sigma}\), we have the object \(\fb\sigma\) in \(\C(T)\). In \(\T\C(T)\),
	this corresponds to a sort \(\ul{\fb\sigma}\). Set \(\ep_T \sigma \equivdef
	\ul{\fb\sigma}\).  For a relation \(R \hookrightarrow
	\sigma_1,\ldots,\sigma_n\) in \(\Sigma\), we have the monic \(\dom_{\fb R}:
	\fb R \to \fb{\vec\sigma}\) in \(\C(T)\). Set%
	\[%
		\pwrap{\ep_T R}(\vec x) \equivdef%
			\exists y^{\ul{\fb R}} \bigwedge_{i=1}^n
			\ul{\pi^{\vec\sigma}_{\sigma_i}}\pwrap{ \ul{\dom_{\fb R}}(y) } = x_i.
	\]%
	Note for \(n=1\), the projection is (by default) the identity morphism
	\(1_{\fb\sigma}\). Lastly for a function symbol \(f :
	\sigma_1,\ldots,\sigma_n \to \tau\) of \(\Sigma\) we have the morphism
	\(\theta_f\) in \(\C(T)\). Define%
	\[%
		\pwrap{\ep_T f}(\vec x,y) \equivdef \exists z^{\ul{\fb{\vec\sigma}}}
		\pwrap{%
			\bigwedge_{i=1}^n \ul{\pi^{\vec\sigma}_{\sigma_i}}(z) = x_i \wedge
			\ul{\theta_f}(z) = y
		}.%
	\]%
\end{definition}

\begin{lemma}\label{lemma:ep-presentation}
	For any formula \(\phi \hookrightarrow \sigma_1,\ldots,\sigma_n\) of \(T\),
	\(\pwrap{\ep_T\phi}\pwrap{x_1,\ldots,x_n}\) is logically equivalent in
	\(\T\C(T)\) to%
	\[%
		\exists y^{\ul{\fb\phi}} \bigwedge_{i=1}^n
		\ul{\pi^{\vec\sigma}_{\sigma_i}}\pwrap{\ul{\dom_{\fb\phi}}\pwrap{y}} = x_i.
	\]%
\end{lemma}
See~\ref{subsection:biequivalence-proofs} for the proof. It relates
Rules~\ref{rule:relations}\-/\ref{rule:term-reduction} to
\axiom{IL} axiom schemata and inducts on formula complexity.

\begin{proposition}\label{prop:epsilon-translation}
	\(\ep_T: T \to \T\C(T)\) is an e.p.\ translation. In particular,
	\(E^{\ep_T}_\sigma(s,t) \dashv\vdash s = t\) for all sorts \(\sigma \in \SSort{\Sigma}^*\).
\end{proposition}%

\begin{proof}[Proof of Proposition~\ref{prop:epsilon-translation}]
	We first show that \(\ep_T\) is a translation. Suppose that \(\phi(x) \vdash
	\psi(x)\) in \(T\), where \(\phi\) and \(\psi\) have domain \(\sigma\). By
	Lemma~\ref{lemma:subobject-lemma}, \(\dom_{\fb\phi} \leq
	\dom_{\fb\psi}\) in \(\Sub{\fb\sigma}\) in \(\C(T)\). So there exists a
	morphism \(f: \fb\phi \to \fb\psi\) in \(\C(T)\) such that \(\dom_{\fb\phi} =
	\dom_{\fb\psi} f\). \axiom{IL2} for this triplet of morphisms implies
	\[%
		\ul{\dom_{\fb\phi}}(y_1) = x \vdash \exists y_2^{\ul{\fb\psi}} \pwrap{%
			y_2 = \ul{f}(y_1) \wedge \ul{\dom_{\fb\psi}}(y_2) = x
		}.%
	\]%
	If we quantify over \(y_1\), we see that the following sequent is provable in
	\(\T\C(T)\).%
	\[%
		\exists y_1^{\ul{\fb\phi}} \ul{\dom_{\fb\phi}}(y_1) = x \vdash%
		\exists y_2^{\ul{\fb\psi}} \ul{\dom_{\fb\psi}}(y_2) = x.
	\]%
	By Lemma~\ref{lemma:ep-presentation}, if \(\phi\) and \(\psi\) are unary, we
	are done. Otherwise, we can apply \axiom{IL6} axioms to the domain \(\sigma\)
	to decompose it into its factors and conclude again using
	Lemma~\ref{lemma:ep-presentation}.

	Lastly, we show that \(E^{\ep_T}_\sigma(s,t) \dashv\vdash s = t\) for any sort
	\(\sigma\). Note that \(\dom_{\fb{x_1 =_\sigma x_2}} : \fb{x_1 =_\sigma x_2}
	\to \fb\sigma \times \fb\sigma\) is the equalizer for the diagram \(\fb\sigma
	\times \fb\sigma \rightrightarrows \fb\sigma\). The definition of \(\ep_T\)
	applied to a relation implies%
	\[%
		E^{\ep_T}_\sigma(s,t) \equiv \exists y^{\ul{\fb{x_1 = x_2}}}
		\ul{\pi^{\sigma,\sigma}_{\sigma,x_1}}\pwrap{\ul{\dom_{\fb{x_1=x_2}}}(y)} = s
		\wedge
		\ul{\pi^{\sigma,\sigma}_{\sigma,x_2}}\pwrap{\ul{\dom_{\fb{x_1=x_2}}}(y)} = t.
	\]%
	From this we see that \axiom{IL5} applied to the equalizer implies the desired
	sequent.
\end{proof}
We proceed to the 1\-/cell component of \(\ep\). Let \(F: T_1 \to T_2\) be an
e.p.\ translation. First we define a homotopy \(\ep_F: \T\C(F) \ep_{T_1}
\Rightarrow \ep_{T_2} F\). Given a \(T_1\)-sort \(\sigma\), where \(F\sigma
\equiv \tau_1,\ldots,\tau_n\), define the \(\sigma\) component of the t\-/map
to be%
\[%
	\pwrap{\ep_F}_\sigma(x,y_1,\ldots,y_n) \equivdef \bigwedge_{i=1}^n
	\ul{\pi^{\vec\tau}_{\tau_i}}\pwrap{%
		\ul{\dom_{D^F_\sigma}}(x)
	} = y_i.%
\]%
This is well\-/defined, since \(\T\C(F)\ep_{T_1}\sigma \equiv
\ul{D^F_\sigma}\), and \(\ep_{T_2} F\sigma \equiv
\ep_{T_2}\pwrap{\tau_1,\ldots,\tau_n} \equiv
\ul{\fb{\tau_1}},\ldots,\ul{\fb{\tau_n}}\).

\begin{proposition}\label{prop:ep_F-definable-iso}
	\(\pwrap{\ep_F}_\sigma\) presents a \(\T\C(T_2)\)\-/definable isomorphism
	\(\pwrap{\ep_F}_\sigma: D_\sigma^{\ep_{T_2}F} \to
	D_\sigma^{\T\C(F)\ep_{T_1}}\).
\end{proposition}
\begin{proposition}\label{prop:ep_F-homotopy}
	\(\ep_F\) is a homotopy.
\end{proposition}%

\begin{proposition}\label{prop:ep-pseudonatural}
	\(\ep\) is a pseudonatural homotopy \(\id_\CThEq \Rightarrow \T\C\).
\end{proposition}%

For proofs of these propositions, see~\ref{subsection:biequivalence-proofs}. We
have reached the final step for \(\ep\).

\begin{proposition}\label{prop:final-step-ep}
	\(\ep_T: T \to \T\C(T)\) is a homotopy equivalence (e.p.\
	bi\-/interpretation) for any coherent theory \(T\).
\end{proposition}
\begin{proof}
	Define a homotopy inverse \(\gamma_T: \T\C(T) \to T\) via the following
	reconstrual. For a sort \(\ul{\fb\phi}\), set \(\gamma_T \ul{\fb\phi}
	\equivdef \Dom\fb\phi\). For a function symbol \(\ul{f}: \ul{\fb\phi} \to
	\ul{\fb\psi}\), there is a corresponding morphism \(f: \fb\phi \to \fb\psi\)
	in \(\C(T)\); hence there is a \(T\)\-/substitution class \(\fb f
	\hookrightarrow \Dom\fb\phi,\Dom\fb\psi\) presenting \(f\). Set \(\gamma_T
	\ul{f} \equivdef \fb f\). Finally set \(E^{\gamma_T}_{\ul{\fb\phi}}(x,y)
	\equivdef \phi(x) \wedge x = y\). Thus, \(\gamma_T\), assuming it is a
	translation, is e.p. To show that \(\gamma_T\) is a translation, it suffices
	to prove that \(\gamma_T\) translates all instances of the axiom schemata
	\axiom{IL1} through \axiom{IL10} into provable sequents of \(T\). This is
	done by matching each \axiom{ILi} axiom with the corresponding sequents in
	\axiom{SCi}. This amounts to proving
	Proposition~\ref{prop:sc-axiom-schemata}, which is elementary.

	We need to find homotopy t\-/maps \(\chi: 1_T \Rightarrow \gamma_T\ep_T\) and
	\(\eta: 1_{\T\C(T)} \Rightarrow \ep_T \gamma_T\). Given a \(T\)\-/sort
	\(\sigma\), set \(\chi_\sigma(s,t) \equivdef s = t\). Given a
	\(\T\C(T)\)\-/sort \(\ul{\fb\phi}\), set%
	\[%
		\eta_{\ul{\fb\phi}}(s,t_1,\ldots,t_n) \equivdef \bigwedge_{i=1}^n
		\ul{\pi^{\vec\sigma}_{\sigma_i}}\pwrap{\ul{\dom_{\fb\phi}}(s)} = t_i,
	\]%
	where \(\Dom\fb\phi \equiv \vec\sigma \equiv \sigma_1,\ldots,\sigma_n\).
	Axioms \axiom[\chi]{TM1} through \axiom[\chi]{TM4}, \axiom[\chi]{TM6}, and
	\axiom[\chi]{TM7} are satisfied due to
	Proposition~\ref{prop:epsilon-translation},
	Rules~\ref{rule:conjunction-preservation} and~\ref{rule:context-duplication},
	and Proposition~\ref{prop:reconstrual-log-eq-to-identity}. For
	\axiom[\chi]{TM5} and \axiom[\chi]{TM8}, it suffices to show that \(\gamma_T
	\ep_T \phi(\vec{t}) \dashv\vdash \phi(\vec{t})\) is provable. By
	Lemma~\ref{lemma:ep-presentation} and
	Rules~\ref{rule:conjunction-preservation},
	\ref{rule:quantifier-preservation}, and~\ref{rule:term-reduction},%
	\[%
		\gamma_T \ep_T \phi(\vec{t}) \dashv\vdash%
		\exists y^{\vec\sigma}\pwrap{%
			\bigwedge_{i=1}^n \exists {t'}^{\vec\sigma} \pwrap{%
				\dom_{\fb\phi}(y,t') \wedge \pi^{\vec\sigma}_{\sigma_i}(t',t_i)
			}%
		} \dashv\vdash \phi(\vec{t}).%
	\]%
	The case for \(\eta\) is less simple. While \axiom[\eta]{TM2} and
	\axiom[\eta]{TM3} are straightforward, \axiom[\eta]{TM1} and
	\axiom[\eta]{TM6} follow from Lemma~\ref{lemma:ep-presentation}.
	\axiom[\eta]{TM4} follows from \axiom[\eta]{TM1} and
	Proposition~\ref{prop:epsilon-translation}. \axiom[\eta]{TM7} is a
	consequence of \axiom{IL6} axioms for the product \(\fb{\vec\sigma}
	\rightrightarrows \fb{\sigma_i}\) and \axiom[\dom_{\fb\phi}]{IL3}. This
	leaves \axiom[\eta]{TM5} and \axiom[\eta]{TM8}.

	For \axiom{TM5}, we need to show that given a \(\T\C(T)\)-substitution class
	\(\fb{A} \hookrightarrow \ul{\fb{\phi_1}},\ldots,\ul{\fb{\phi_n}}\), where
	\(\fb{\phi_i} \hookrightarrow \sigma_{i 1}, \ldots, \sigma_{i m_i}\),%
	\[%
		1_{\T\C(T)}A(\vec{s}) \wedge \bigwedge_{i=1}^n\bigwedge_{j=1}^{m_i}%
		\ul{\pi^{\vec{\sigma_i}}_{\sigma_{i j}}}\pwrap{%
			\ul{\dom_{\fb{\phi_i}}}(s_i)
		} = t_{i j} \vdash%
		\ep_T \gamma_T A(\vec{t}).
	\]%
	We prove this claim by induction on the complexity of formulae. The base case
	is the one where \(A(\vec{s})\) is an atomic formula. The only relation
	symbols in the signature of \(\T\C(T)\) are equality relations. Therefore,
	using \axiom{IL1} and \axiom{IL2} axioms if necessary, 
        if \(A(\vec{s})\) is atomic, then it
	is logically equivalent to a formula of the form \(\ul{f}(s_1) =
	\ul{g}(s_2)\), where \(\ul{f}: \ul{\fb{\phi_1}} \to \ul{\fb\psi}\) and
	\(\ul{g}: \ul{\fb{\phi_2}} \to \ul{\fb\psi}\) are function symbols.
	Thus assume without loss of generality that \(A(\vec{s}) \equiv
	A(s_1,s_2) \equiv \ul{f}(s_1) = \ul{g}(s_2)\).	

	The \(\T\C(T)\)\-/sort \(\ul{\fb\psi}\) comes from a substitution class
	\(\fb\psi\) in \(T\). Let \(\tau_1,\ldots,\tau_m\) denote the domain of
	\(\fb\psi\). Using Rules~\ref{rule:conjunction-preservation},
	\ref{rule:context-duplication}, and~\ref{rule:term-reduction}, we can
	rearrange the right side of \axiom[\eta]{TM5} to show%
	\[%
		\ep_T \gamma_T A(\vec{t}) \dashv\vdash%
		\exists z_1^{\ep_T \gamma_T \ul{\fb\psi}}%
		\exists z_2^{\ep_T \gamma_T \ul{\fb\psi}}\pwrap{%
			\ep_T f(\vec{t_1},z_1) \wedge \ep_T g(\vec{t_2},z_2) \wedge%
			E^{\ep_T \gamma_T}_{\ul{\fb\psi}}(z_1,z_2)
		}.%
	\]%
	Since \(E^{\ep_T\gamma_T}_{\ul{\fb\psi}}(z_1,z_2) \equiv \ep_T\psi(z_1)
	\wedge z_1 = z_2\), the above is logically equivalent to%
	\[%
		\exists z_1^{\ep_T\gamma_T\ul{\fb\psi}}\pwrap{%
			\ep_T f(\vec{t_1},z_1) \wedge \ep_T g(\vec{t_2},z_1) \wedge \ep_T
			\psi(z_1)
		}.%
	\]%
	Recall the universal property of a product ensures
	\(\pi^{\vec\sigma,\vec\tau}_{\tau_i} = \pi^{\vec\tau}_{\tau_i}
	\pi^{\vec\sigma,\vec\tau}_{\vec\tau}\). With this in mind, we can expand each
	conjunct using Definition~\ref{def:ep_T} and
	Lemma~\ref{lemma:ep-presentation} and simplify. If we also consider
	\axiom{IL6}, we can further simplify, showing%
	\begin{align}\label{eq:TM5-RHS}%
		\ep_T \gamma_T A(\vec t) \dashv\vdash%
		\exists y_1^{\ul{\fb f}} \exists y_2^{\ul{\fb g}}\Bigg(%
			&\ul{\pi^{\vec{\sigma_1},\vec\tau}_{\vec\tau}}\pwrap{%
				\ul{\dom_{\fb f}}(y_1)
			} = \ul{\pi^{\vec{\sigma_2},\vec\tau}_{\vec\tau}}\pwrap{%
				\ul{\dom_{\fb g}}(y_2)	
			}\notag\\%
			&\wedge \bigwedge_{j=1}^{m_1}%
			\ul{\pi^{\vec{\sigma_1}}_{\sigma_{1 j}}}\pwrap{%
				\ul{\pi^{\vec{\sigma_1},\vec\tau}_{\vec{\sigma_1}}}\pwrap{%
					\ul{\dom_{\fb f}}(y_1)
				}%
			} = t_{1 j}\notag\\%
			&\wedge \bigwedge_{j=1}^{m_2}%
			\ul{\pi^{\vec{\sigma_2}}_{\sigma_{2 j}}}\pwrap{%
				\ul{\pi^{\vec{\sigma_2},\vec\tau}_{\vec{\sigma_2}}}\pwrap{%
					\ul{\dom_{\fb g}}(y_2)
				}%
			} = t_{2 j}\notag\\%
			&\wedge\exists w^{\ul{\fb\psi}} \pwrap{%
				\ul{\dom_{\fb\psi}}(w) =%
				\pi^{\vec{\sigma_1},\vec\tau}_{\vec\tau}\pwrap{%
					\ul{\dom_{\fb f}}(y_1)
				}%
			}%
		\Bigg).%
	\end{align}%
	The last conjunct on the right side is tautological, due to \axiom[f]{DM1},
	so we can omit it. Now we expand the left side of \axiom{TM5}.
	Proposition~\ref{prop:reconstrual-log-eq-to-identity} allows us to replace
	\(1_{\T\C(T)} A\) with \(A\). Thus the left side of \axiom{TM5} is
	logically equivalent to%
	\begin{equation}\label{eq:TM5-LHS}%
		\ul{f}(s_1) = \ul{g}(s_2) \wedge \bigwedge_{j=1}^{m_1}%
		\ul{\pi^{\vec{\sigma_1}}_{\sigma_{1 j}}}\pwrap{%
			\ul{\dom_{\fb{\phi_1}}}(s_1)
		} = t_{1 j} \wedge \bigwedge_{j=1}^{m_2}%
		\ul{\pi^{\vec{\sigma_2}}_{\sigma_{2 j}}}\pwrap{%
			\ul{\dom_{\fb{\phi_2}}}(s_2)
		} = t_{2 j}.%
	\end{equation}%
	The conjunctions in Formula~\ref{eq:TM5-LHS} are similar in form to those in
	Formula~\ref{eq:TM5-RHS}. The key to proving \axiom{TM5} is making that
	similarity precise. Consider the morphism \(p_1: \fb{f} \to \fb{\phi_1}\)
	presented by \(p_1(x_2,x_2,y) \equivdef f(x_1,x_2) \wedge x_1 = y\). Since
	\axiom[p_1]{SC9} is provable, \(p_1\) is a regular epimorphism (see
	Proposition~\ref{prop:sc-axiom-schemata}), so \axiom[p_1]{IL9} is an axiom of
	\(\T\C(T)\). Furthermore, \(\dom_{\fb{\phi_1}} p_1 =
	\pi^{\vec{\sigma_1}\vec\tau}_{\vec{\sigma_1}} \dom_{\fb{f}}\). Combining
	\axiom{IL2} of this equation with \axiom[p_1]{IL9} shows that in \(\T\C(T)\)
	we have%
	\[%
		\vdash \exists y_1^{\ul{\fb{f}}}%
		\ul{\pi^{\vec{\sigma_1}\vec\tau}_{\vec{\sigma_1}}}\pwrap{%
			\ul{\dom_{\fb{f}}}(y_1)
		} = \ul{\dom_{\fb{\phi_1}}}(s_1).
	\]%
	Therefore the second conjunct of Formula~\ref{eq:TM5-LHS} entails the second
	conjunct of Formula~\ref{eq:TM5-RHS}. If we let \(p_2: \fb{g} \to
	\fb{\phi_2}\) be the morphism presented by \(p_2(x_1,x_2,y) \equivdef
	g(x_1,x_2) \wedge x_1 = y\), then a similar argument shows that the third
	conjunct of Formula~\ref{eq:TM5-LHS} entails the third conjunct of
	Formula~\ref{eq:TM5-RHS}. 

	This leaves the first conjunct. The other conjuncts were proven using
	\axiom[p_1]{IL9} and \axiom[p_2]{IL9}, therefore it suffices to prove the
	following sequent.
	\begin{align*}%
		\ul{f}(s_1) = \ul{g}(s_2) \wedge \ul{p_1}(y_1) =\,\,& s_1 \wedge \ul{p_2}(y_2) =
		s_2\\%
		&\vdash \ul{\pi^{\vec{\sigma_1}\vec\tau}_{\vec\tau}}\pwrap{%
			\ul{\dom_{\fb{f}}}(y_1)
		} = \ul{\pi^{\vec{\sigma_2}\vec\tau}_{\vec\tau}}\pwrap{%
			\ul{\dom_{\fb{g}}}(y_2)
		}%
	\end{align*}%
	We can eliminate the variables \(s_1\) and \(s_2\), so it suffices to prove
	the sequent%
	\[%
		\ul{f}\pwrap{\ul{p_1}(y_1)} = \ul{g}\pwrap{\ul{p_2}(y_2)} \vdash%
		\ul{\pi^{\vec{\sigma_1}\vec\tau}_{\vec\tau}}\pwrap{%
			\ul{\dom_{\fb{f}}}(y_1)
		} = \ul{\pi^{\vec{\sigma_2}\vec\tau}_{\vec\tau}}\pwrap{%
			\ul{\dom_{\fb{g}}}(y_2)
		}.%
	\]%
	Note \(\dom_{\fb\psi} f p_1 = \pi^{\vec{\sigma_1}\vec\tau}_{\vec\tau}
	\dom_{\fb{f}}\) and \(\dom_{\fb\psi} g \ p_2 =
	\pi^{\vec{\sigma_2}\vec\tau}_{\vec\tau} \dom_{\fb{g}}\). Therefore if we
	apply \(\ul{\dom_{\fb\psi}}\) to both terms on the left side of the above
	sequent, we obtain the right side using \axiom{IL2}. Thus
	Formula~\ref{eq:TM5-LHS} entails Formula~\ref{eq:TM5-RHS}, establishing the
	base case for \axiom[\eta]{TM5}.

	The inductive step for \axiom[\eta]{TM5} follows from 
        Rules~\ref{rule:relations}\-/\ref{rule:term-reduction}, similar to earlier
	proofs. The proof for \axiom{TM8} is similar enough to \axiom{TM5} to omit,
	so we conclude.
\end{proof}

\subsubsection{The Pseudonatural Homotopy \(\delta\)}
\begin{definition}\label{def:delta}
	Let \(C\) be a small coherent category. Define a coherent functor \(\delta_C:
	C \to \C\T(C)\) in the following manner. For an object \(A\) of \(C\), set
	\(\delta_C A\) to be the object \(\fb{\ul{A}}\). For a morphism \(f: A \to
	B\) of \(C\), set \(\delta_C f\) to be the morphism \(\theta_{\ul{f}}:
	\fb{\ul A} \to \fb{\ul B}\), which we recall is presented by \(\theta_{\ul
	f}(x,y) \equiv \ul{f}(x) = y\).
\end{definition}%

We will show that \(\delta_C\) is a coherent functor when we prove
\axiom[\delta]{PNT1} later. All that remains are the pseudonatural homotopies
\(\delta_\Eff: \C\T(\Eff) \delta_{C_1} \Rightarrow \delta_{C_2} \Eff\). By the
following proposition, we set \(\delta_\Eff \eqdef \One^\Eff\).

\begin{proposition}\label{prop:delta-F-one}
	For a coherent functor \(\Eff : C_1 \to C_2\), \(\C\T(\Eff)\delta_{C_1} =
	\delta_{C_2} \Eff\).
\end{proposition}
\begin{proof}
	For an object \(A \in C_1\), note that \(\C\T(\Eff)\delta_{C_1} A \equiv
	\C\T(\Eff)\fb{\ul{A}} \equiv \fb{\T(\Eff)\ul{A}} \equiv \fb{\ul {\Eff A}}
	\equiv \delta_{C_2} \Eff A\) by the definitions of \(\C\) and \(\T\) on
	0-cells. Similarly, for any morphism \(f : A \to B\) in \(C_1\), we have
	\(\C\T(\Eff)\delta_{C_1} f \equiv \C\T(\Eff)\theta_ {\ul{f}} \equiv
	\theta_{\T(\Eff)\ul{f}} \equiv \theta_{\ul{\Eff f}} \equiv \delta_{C_2} \Eff
	f\). Thus \(\delta_{C_2} \Eff = \C\T(\Eff) \delta_{C_1}\).
\end{proof}
\begin{proposition}\label{prop:delta-pseudonatural}
	\(\delta\) is a pseudonatural homotopy \(\delta: \id_\Coh \Rightarrow
	\C\T\).
\end{proposition}

Proposition~\ref{prop:delta-pseudonatural} is proven
in~\ref{subsection:biequivalence-proofs}. We have reached the
final step for \(\delta\).

\begin{proposition}\label{prop:final-step-delta}
	\(\delta_C: C \to \C\T(C)\) is a homotopy equivalence (equivalence of
	categories), for any small coherent category \(C\).
\end{proposition}
\begin{proof}
	There are two key elements to this proof.
	First,~\cite{MakkaiReyes1977} generalizes the notion of \(\Set\)\-/valued
	models of a coherent theory \(T\) to models valued in an arbitrary coherent
	category. This is done by sending the logical symbols of \(T\) to subobjects
	and morphisms in the coherent category. Moreover, these models may be
	pushed forward along a coherent functor via composition (see Chapter 8
	of~\cite{MakkaiReyes1977}). Second, there exist models \(M_0: \T(C) \to
	\C\T(C)\) and \(M: \T(C) \to C\) which are universal with respect to this
	pushforward operation (see Propositions 8.1.2 and 8.2.3
	of~\cite{MakkaiReyes1977}).

	\(M_0: \T(C) \to \C\T(C)\) is the model from~\cite[p.\
	243]{MakkaiReyes1977} which sends the sort \(\ul{A}\) to the object
	\(\fb{\ul{A}}\) and the function symbol \(\ul{f}: \ul{A} \to \ul{B}\) to the
	definable map presented by \(\theta_{\ul{f}}\). \(M: \T(C) \to C\) is the
	\emph{canonical interpretation} of~\cite[p.\ 82]{MakkaiReyes1977}, sending
	\(\ul{A}\) to \(A\) and \(\ul{f}\) to \(f\). We first note that \(M_0 =
	\delta_C \circ M\). On sorts we have \(\delta_C M \ul{A} \equiv \delta_C
	A \equiv \fb{\ul{A}} \equiv M_0 \ul{A}\). On function symbols we have
	\(\delta_C M \ul{f} \equiv \delta_C f \equiv \theta_{\ul{f}} \equiv M_0
	\ul{f}\). On the equality symbol \(M_0(=_{\ul{A}})\) is the subobject of
	\(M_0 \ul{A} \times M_0 \ul{A}\) presented by \(\fb{x =_{\ul{A}} y}\). This
	is the diagonal subobject. On the other hand, \(M(=_{\ul{A}})\) is defined to
	be the diagonal subobject of \(A \times A\). Since \(\delta_C\) preserves
	pullbacks, \(\delta_C\) preserves diagonal subobjects, so \(\delta_C
	M(=_{\ul{A}}) \equiv M_0(=_{\ul{A}})\). Secondly, we note that Proposition
	8.2.3 of~\cite{MakkaiReyes1977} implies the existence of a coherent functor
	\(\iota_C: \C\T(C) \to C\) such that \(M = \iota_C \circ M_0\). 
        With the last equation, this implies \(M = \iota_C \delta_C \circ M\)
	and \(M_0 = \delta_C \iota_C \circ M_0\). Propositions 8.1.2 and 8.2.3
	of~\cite{MakkaiReyes1977} imply that \(\iota_C \delta_C\) and \(\delta_C
	\iota_C\) are uniquely determined up to natural isomorphism by these
	equations. On the other hand, the identities \(1_C\) and \(1_{\C\T(C)}\) also
	satisfy these equations. Therefore \(\iota_C \delta_C \simeq 1_C\) and
	\(\delta_C \iota_C \simeq 1_{\C\T(C)}\), so \(\delta_C\) is an equivalence
	with homotopy inverse \(\iota_C\).
\end{proof}
By Propositions~\ref{prop:ep-pseudonatural}, \ref{prop:final-step-ep},
\ref{prop:delta-pseudonatural}, and \ref{prop:final-step-delta}, we have
established the \(\CThEq-\Coh\) correspondence.
\begin{theorem}\label{theorem:CThEqCoh-Correspondence}%
	The pseudofunctors \(\C: \CThEq \to \Coh\) and \(\T: \Coh \to \CThEq\) form a
	biequivalence.
\end{theorem}
\begin{corollary}
	Two coherent theories \(T_1\) and \(T_2\) are e.p.\ bi\-/interpretable if and
	only if \(\C(T_1)\) and \(\C(T_2)\) are equivalent categories.
\end{corollary}
\begin{remark}
	Theorem \ref{theorem:CThEqCoh-Correspondence} is compatible with Makkai and
	Reyes' \cite{MakkaiReyes1977} in the sense that the inverse of the component
	\(\delta_C\) of the pseudonatural homotopy \(\delta\) is the functor
	\(\iota_C: \C\T(C) \to C\) induced by the canonical interpretation via
	Proposition \ref{prop:final-step-delta}. As a result, the ``unsatisfactory''
	definition of interpretation alluded to in Proposition 8.1.1
	of~\cite{MakkaiReyes1977} is subsumed by our notion of equality-preserving
	translation, in a manner that is unique up to (t-map) homotopy.
\end{remark}
\begin{remark}\label{remark:WhyWeCantGetStrict2Eq}
	With our canonical notion of translation, we should not expect a strict
	version of Theorem~\ref{theorem:CThEqCoh-Correspondence}. This is because a
	translation sends function symbols to substitution classes, whereas a functor
	sends morphisms to morphisms. Since the function symbols of \(\T(C)\)
	correspond to morphisms of \(C\), the compositor for \(\T\) cannot be
	trivial.
\end{remark}

\section{Bi\-/Interpretability}

Using the results developed so far, we prove the third main theorem of this
paper.%

\begin{theorem}\label{theorem:BiInterpretability-result}%
	Let \(T_1\) and \(T_2\) be small coherent theories. \(T_1\) and \(T_2\) are
	bi\-/interpretable if and only if \(\C(T_1)^\ex \approx \C(T_1^\eq)\) and
	\(\C(T_2)^\ex \approx \C(T_2^\eq)\) are equivalent categories.
\end{theorem}%

We will actually prove a stronger result: that the homotopy categories
\(\h\CTh_0\) and \(\h\Exact\) are equivalent. First we introduce
pertinent notation.

\begin{definition}\label{def:quotients}
	A \textbf{congruence} over an object \(A\)
	in a coherent category \(C\) is a
	monomorphism \(R \hookrightarrow A \times A\) such that, via the Yoneda
	embedding, \(\Hom\pwrap{X,R} \hookrightarrow \Hom\pwrap{X,A\times A} \cong
	\Hom\pwrap{X,A} \times \Hom\pwrap{X,A}\) defines an equivalence relation for
	any object \(X\) of \(C\). This can be axiomatized using the internal logic
	\(\T(C)\) (Definition 3.3.6 of~\cite{MakkaiReyes1977}). We say that a
	congruence \(R\) admits a \textbf{quotient} \(A / R\) if there exists a
	pullback diagram%
	\[%
		\begin{tikzcd}
			R \arrow[r,"p_2", shift left=1] \arrow[r,"p_1", shift right=1, swap] & A%
			\arrow[r,"q"] & B,
		\end{tikzcd}
	\]%
	where \(q\) is the coequalizer of \(p_1\) and \(p_2\). More generally, a
	morphism \(q: A \to B\) is an \textbf{effective epimorphism} if it is the
	coequalizer of its kernel pair \(A \times_q A \rightrightarrows A\). In the
	language of~\cite{borceux_1994}, a quotient diagram is also called an
	\emph{exact sequence}.
\end{definition}

Congruences and quotients are compatible with pullbacks in the following sense.
Given a monomorphism \(\alpha: X \hookrightarrow A\), a congruence \(p_1
\times p_2 : R \hookrightarrow A \times A\), and a pullback square%
\[%
		\begin{tikzcd}
			\pwrap{\alpha\times\alpha}^* R
			\arrow[d,"\pwrap{\alpha\times\alpha}^*\pwrap{p_1 \times p_2}", swap]
			\arrow[r] & R \arrow[d,"p_1\times p_2"]\\%
			X \times X \arrow[r,"\alpha \times \alpha"] & A \times A,
		\end{tikzcd}
\]%
\(\pwrap{\alpha\times\alpha}^* R\) is a congruence over \(X\). For that reason,
we abbreviate the pullback \(\pwrap{\alpha\times\alpha}^* R\) to
\(\restr{R}{X}\). If we write \(p_1^X,p_2^X : X \times X \to X\) for the
projection morphisms and define \(\restr{p_i}{X} \eqdef p_i^X \circ
\pwrap{\alpha\times\alpha}^*\pwrap{p_1\times p_2}\), then for any quotient \(q:
A \twoheadrightarrow B\) of \(R\), we obtain a pullback diagram (by a simple
diagram chase)%
\[%
	\begin{tikzcd}
		\restr{R}{X} \arrow[r,"\restr{p_2}{X}", shift left=1]
		\arrow[r,"\restr{p_1}{X}", shift right=1,swap] & X \arrow[r,"q\alpha"] &
		B.
	\end{tikzcd}
\]%
Furthermore, the coimage \(\restr{q}{X} : X \twoheadrightarrow \exists_q X\)
from the image factorization \(q\alpha = \pwrap{\exists_q
\alpha}\pwrap{\restr{q}{X}}\) is a coequalizer for \(\restr{R}{X}\), so the
quotient \(X / \restr{R}{X}\) is presented by the regular epimorphism
\(\restr{q}{X}\). See Propositions 2.5.7 and 2.5.8 of~\cite{borceux_vol1} and
Chapter 2 of~\cite{borceux_1994} for details.

We will also need a more careful treatment of the \emph{canonical
interpretation} \(M: \T(C) \to C\) of~\cite[p.~82]{MakkaiReyes1977}. In
particular, subobjects of \(\C(T)\) can be described in terms of substitution
classes in both \(T\) and \(\T\C(T)\).

\begin{definition}[\cite{MakkaiReyes1977}, Chapter 2, Section 4]
	Let \(C\) be a coherent category, and let \(\alpha: A \hookrightarrow X_1
	\times \ldots \times X_n\) be a monomorphism in \(C\), presenting a subobject
	called, by the usual abuse of notation, \(A \in \Sub\pwrap{X_1\times \ldots
	\times X_n}\). Define the \(\T(C)\)\-/formula
	\(
		\utilde{A}(x_1,\ldots,x_n) \equivdef \exists a^{\ul{A}} \bigwedge_{i=1}^n
		\ul{\pi^{X_1 \times \ldots \times
		X_n}_{X_i}\mkern-55mu}\mkern55mu\pwrap{\ul{\alpha}(a)} = x_i.
	\)
	If \(C = \C(T)\), a \(T\)\-/substitution class \(\fb\phi \hookrightarrow
	\sigma_1,\ldots,\sigma_n\) presents a subobject given by \(\dom_{\fb\phi}: \fb\phi
	\hookrightarrow \fb{\sigma_1} \times \ldots \times \fb{\sigma_n}\). The
	construction above identifies a \(\T\C(T)\)\-/formula which, for the sake of
	simple notation, we denote by \(\utilde\phi\):
	\(
		\utilde{\phi}(x_1,\ldots,x_n) \equivdef \exists t^{\ul{\fb\phi}}
		\bigwedge_{i=1}^n \ul{\pi^{\vec\sigma}_{\sigma_i}}\pwrap{\ul{\dom_{\fb\phi}}(t)} =
		x_i.
	\)
\end{definition}

By making their construction more explicit, it becomes clear by
Lemma~\ref{lemma:ep-presentation} that \(\utilde{\phi}(\vec x)\) is
logically equivalent to \(\ep_T \phi(\vec x)\). We make this identification
moving forward. Makkai and Reyes construct the canonical interpretation \(M\) such that
\(M\pwrap{\ep_T A}\) is the subobject \(A\) (in fact, it is this property
that justifies the definition of the \emph{extended canonical language}). In
particular, \(M(\ep_T\fb\phi)\) is the subobject presented by
\(\dom_{\fb\phi}\) in a syntactic category \(\C(T)\) for a \(T\)\-/substitution
class \(\fb\phi\). This motivates the following lemma.

\begin{lemma}\label{lemma:super-subobject-lemma}
	Let \(T\) be a coherent theory. Let \(\phi,\psi \hookrightarrow
	\sigma_1,\ldots,\sigma_n\) be a pair of \(T\)\-/formulae. The following are
	equivalent: (1) \(\fb\phi \vdash \fb\psi\) in \(T\); (2) \(\dom_{\fb\phi}
	\leq \dom_{\fb\psi}\) as subobjects of \(\fb{\sigma_1} \times \ldots \times
	\fb{\sigma_n}\) in \(\C(T)\); (3) \(\ep_T \fb\phi \vdash \ep_T\fb\psi\) in
	\(\T\C(T)\).
\end{lemma}

\begin{proof}
	\((1) \implies (2)\) is half of Lemma~\ref{lemma:subobject-lemma}. \((2)
	\implies (3)\) follows from the proof of
	Proposition~\ref{prop:epsilon-translation}. For \((3) \implies (1)\), we
	recall the proof of Proposition~\ref{prop:final-step-ep}, which showed that
	\(\ep_T: T \to \T\C(T)\) has a homotopy inverse \(\gamma_T : \T\C(T) \to T\)
	such that \(\gamma_T \ep_T \fb\phi \dashv\vdash \fb\phi\). Thus \(\ep_T
	\fb\phi \vdash \ep_T \fb\psi\) implies \(\gamma_T \ep_T \fb\phi \vdash
	\gamma_T \ep_T \fb\psi\), which holds if and only if \(\fb\phi \vdash
	\fb\psi\).
\end{proof}

\subsection{The functor \(\X\)}

\(\X\) is a categorification of \emph{exact completion}
(from~\cite{MakkaiReyes1977}), except we consider translations which are not
necessarily e.p. Recall the construction in~\cite{MakkaiReyes1977}.

\begin{proposition}[\cite{MakkaiReyes1977}, Theorem 8.4.3]%
	\label{prop:exact-completion}
	Given a coherent category \(C\), there exists a Barr\-/exact coherent
	category \(C^\ex\) and a coherent functor \(I: C \to C^\ex\) satisfying the
	following properties.%

	(\axiom{EC1}) \(I\) is conservative, fully faithful, and full on subobjects.
	
	(\axiom{EC2}) Any object of \(C^\ex\) is isomorphic to a quotient \(I(A) /
	I(R)\), where \(A\) is an object of \(C\), and \(R \hookrightarrow A \times
	A\) is a congruence in \(C\).
	
	(\axiom{EC3}) For any coherent functor \(\Eff: C \to D\), where \(D\) is a
	Barr\-/exact coherent category, there exists an extension \(\Eff^\ex: C^\ex \to
	D\), unique up to natural isomorphism, such that \(\Eff=\Eff^\ex\circ I\).%
	\[%
		\begin{tikzcd}%
			C^\ex \arrow[r,dashed,"\Eff^\ex"] & D \\
			C \arrow[u,"I"] \arrow[ur,"\Eff",swap] &
		\end{tikzcd}%
	\]%

	(\axiom{EC4}) There exists a coherent theory \(T = \T(C)^\ex\) such that
	\(C^\ex = \C(T)\), where \(T\) is a conservative extension of \(\T(C)\)
	obtained by adjoining for every pair \(\pwrap{A,R}\), with \(r: R
	\hookrightarrow A \times A\) a congruence over \(A\) in \(C\), a sort symbol
	\(\ul{A} / \ul{R}\), an equality relation for this new sort, and a function
	symbol \(q_R: \ul{A} \to \ul{A}/\ul{R}\) along with three axioms (per pair):%
	\begin{align}%
		&\vdash \exists a \; q_R(a) = x \tag{Q1},\\
		q_R(a) = q_R(a') \dashv\!\!\!\:&\vdash \exists t^{\ul{R}} \pwrap{%
			\ul{\pi_1}\pwrap{\ul{r}(t)} = a \wedge \ul{\pi_2}\pwrap{\ul{r}(t)} = a'
		} \tag{Q2}.
	\end{align}%
\end{proposition}

\begin{remark}\label{remark:definition-of-I}
	Since \(\T(C)^\ex\) is an extension of \(\T(C)\),~\cite{MakkaiReyes1977}
	defined the coherent functor \(I: C \to C^\ex\) by observing an ``inclusion''
	interpretation from \(\T(C)\) to \(\T(C)^\ex\). An interpretation
	in~\cite{MakkaiReyes1977} is similar to our notion of an e.p.\ translation,
	except interpretations send function symbols to function symbols. Unpacking
	how this interpretation defines \(I\), we can see that this identifies \(C\)
	with the full subcategory of \(C^\ex\) spanned by those objects which are not
	formal quotients \(\fb{\ul A / \ul R}\). Under this identification, the
	inclusion functor \(I: C \to \C^\ex\) is identified with \(\delta_C\) from
	Proposition~\ref{prop:final-step-delta}, sending an object \(A\) of \(C\) to
	the object \(\fb{\ul A}\) of \(C^\ex = \C(T)\) and a morphism \(f: A \to B\)
	to the morphism in \(C^\ex\) presented by the substitution class \(\fb{
	\ul{f}(x) = y }\). This characterization suffices for our work.
\end{remark}

Because of EC1, we will identify \(C\) with a subcategory of \(C^\ex\). We need
to extend this completion to lift 2\-/cells.

\begin{proposition}\label{prop:exact-completion-2-cells}
	Let \(\chi: \Eff \Rightarrow \Gee\) be a natural transformation between coherent
	functors \(\Eff,\Gee: C \to D\), where \(D\) is Barr\-/exact. There exists a
	natural transformation \(\chi^\ex: \Eff^\ex \Rightarrow \Gee^\ex\). Furthermore
	\(\chi \mapsto \chi^\ex\) defines a functor \(\Hom\pwrap{C,D} \to
	\Hom\pwrap{C^\ex,D}\).
\end{proposition}

The proof of this proposition is based on property EC2 for \(C^\ex\): any
object of \(C^\ex\) is a coequalizer of a diagram \(p_1,p_2: R
\rightrightarrows A\) in \(C\). Therefore via the universal property of
colimits, we can define the components of \(\chi^\ex\) by using the components
of \(\chi\) along \(R\) and \(A\).
See~\ref{subsection:exact-completion-properties} for details.

\begin{corollary}\label{cor:exact-preserves-homotopy}
	Let \(\Eff,\Gee : C^\ex \to D\) be coherent functors where \(D\) is Barr\-/exact.
	Let \(I: C \to C^\ex\) be the inclusion functor from
	Proposition~\ref{prop:exact-completion}. If \(\Eff \circ I\) is naturally
	isomorphic to \(\Gee \circ I\), then \(\Eff\) is naturally isomorphic to \(\Gee\).
\end{corollary}

\begin{proposition}\label{prop:exact-reflects-homotopy}
	The functor \(\chi \mapsto \chi^\ex\) reflects isomorphism. Therefore
	\(\Eff^\ex\) is naturally isomorphic to \(\Gee^\ex\) if and only if \(\Eff\) is
	naturally isomorphic to \(\Gee\).
\end{proposition}

\begin{proof}
	The inclusion \(I : C \to C^\ex\) is faithful. Therefore, given a natural
	isomorphism \(\chi^\ex: \Eff^\ex \Rightarrow \Gee^\ex\), the component along
	an object \(I X\), where \(X\) is an object of \(C\), is an isomorphism
	\(\pwrap{\Eff^\ex \circ I} X \to \pwrap{\Gee^\ex \circ I} X\). Since
	\(\Eff^\ex \circ I = \Eff\) and \(\Gee^\ex \circ I = \Gee\), this is an
	isomorphism \(\Eff X \to \Gee X\) in \(D\). Call this isomorphism \(\chi_X\).
	This defines a natural isomorphism \(\chi: \Eff \Rightarrow \Gee\) since, for
	any morphism \(f: X \to Y\) in \(C\), \(\Gee f \circ \chi_X = \Gee^\ex I f
	\circ \chi^\ex_{I X} = \chi^\ex_{I Y} \circ \Eff^\ex I f = \chi_Y \circ \Eff
	f \).
\end{proof}

On 0\-/cells, the functor \(\X\) is defined by \(\X(T) \eqdef \C(T)^\ex\). On
1\-/cells we need to lift a (generally not e.p.) translation \(T_1 \to T_2\) to
a coherent functor \(\C(T_1) \to \C(T_2)^\ex\), then we will invoke property
EC3 of \(\C(T_1)^\ex\). We will invoke
Theorem~\ref{theorem:CThEqCoh-Correspondence}. Let \(F: T_1 \to T_2\) be a
translation. For any \(T_1\)\-/sort \(\sigma\), \(E^F_\sigma\) is a congruence
over \(D^F_\sigma\), where \(r : E^F_\sigma \to D^F_\sigma \times D^F_\sigma\)
is presented by \(r(x_1,x_2,y_1,y_2) \equiv E^F_\sigma(x_1,x_2) \wedge x_1 =
y_1 \wedge x_2 = y_2\), so the theory \(\T(\C(T_2))^\ex\) introduces a quotient
sort symbol \(\ul{D^F_\sigma} / \ul{E^F_\sigma}\). Abbreviate this symbol to
\(\ul{Q^F_\sigma}\). Call the associated quotient function symbol \(q^F_\sigma:
\ul{D^F_\sigma} \to \ul{Q^F_\sigma}\). Moreover, given \(n\) function symbols
\(f_i : \vec{\sigma_i} \to \tau_i\), we shall make the following abbreviation
for the product morphism \(f_1 \times \ldots \times f_n\):
\[%
	\widearrow{f}\pwrap{\vec{x},\vec{y}} \equivdef \pwrap{%
			\bigwedge_{i=1}^n f_i(x_i) = y_i 
	} : \fb{\vec\sigma_1} \times \ldots \times \fb{\vec\sigma_n} \to \fb{\tau_1}
	\times \ldots \times \fb{\tau_n}.
\]%

\begin{definition}\label{def:F-prime}
	Let \(F: T_1 \to T_2\) be a translation. We define a new reconstrual \(F^\eq:
	T_1 \to \T(\C(T_2))^\ex\) in the following manner. For a \(T_1\)\-/sort
	symbol \(\sigma\), a \(T_1\)\-/relation \(R \hookrightarrow
	\sigma_1,\ldots,\sigma_n\) (including the equality relations), or a function
	symbol \(f: \sigma_1,\ldots,\sigma_n \to \tau\), we define \(F^\eq\)
	according to the rules:
	\begin{align*}%
		F^\eq\sigma \equivdef \ul{Q^F_\sigma}, \quad F^\eq R (\vec{x}) &\equivdef
		\exists {x_1'}^{\ul{D^F_{\sigma_1}}} \exists {x_1''}^{\ul{\fb{F\sigma_1}}}
		\ldots \exists {x_n'}^{\ul{D^F_{\sigma_n}}} \exists
		{x_n''}^{\ul{\fb{F\sigma_n}}} \,\Big(\ep_{T_1} FR (\pvec{x}') \\ 
		&\wedge \ul{\widearrow{\dom}_{D^F_{\sigma}}}(\pvec{x}'',\pvec{x}') \wedge
		\widearrow{q}^F_{\sigma} (\pvec{x}'',\vec{x})\Big),
	\end{align*}%
	\begin{align*}%
		F^\eq f(\vec{x},y) &\equivdef \exists {x_1'}^{\ul{D^F_{\sigma_1}}} \exists
		{x_1''}^{\ul{\fb{F\sigma_1}}} \ldots \exists {x_n'}^{\ul{D^F_{\sigma_n}}}
		{x_n''}^{\ul{\fb{F\sigma_n}}} \exists
		{y'}^{\ul{D^F_\tau}}\exists{y''}^{\ul{\fb{F\tau}}}\,\Big( \ep_{T_1} F
		f(\pvec{x}',y) \\
		&\wedge \ul{\widearrow{\dom}_{D^F_{\sigma}}}(\pvec{x}'',\pvec{x}') \wedge
		\widearrow{q}^F_{\sigma}(\pvec{x}',\vec{x}) \wedge \ul{\dom_{D^F_\tau}}(y'')
		= y' \wedge q^F_\tau(y'') = y\Big).
	\end{align*}%
\end{definition}

Similar to Lemma~\ref{lemma:ep-presentation}, \(F^\eq\) has a uniform
presentation for any \(T_1\)\-/formula.

\begin{lemma}\label{lemma:eq-presentation}
	Let \(\phi \hookrightarrow \sigma_1,\ldots,\sigma_n\) be a \(T_1\)\-/formula
	(or substitution class). The following logical equivalence is provable in
	\(\T(\C(T_2))^\ex\):
	\begin{align*}%
		F^\eq \phi(\vec{x}) &\dashv\vdash \exists {x_1'}^{\ul{D^F_{\sigma_1}}}
		\exists {x_1''}^{\ul{\fb{F\sigma_1}}} \ldots \exists
		{x_n'}^{\ul{D^F_{\sigma_n}}} \exists {x_n''}^{\ul{\fb{F\sigma_n}}}\\
		&\Big(\ep_{T_1}F\phi (\vec{x}') \wedge
		\ul{\widearrow{\dom}_{D^F_{\sigma}}}(\vec{x}'',\vec{x}') \wedge
	\widearrow{q}^F_{\sigma}(\vec{x}'',\vec{x})\Big).
	\end{align*}%
\end{lemma}

\begin{proof}
	This proof is based on induction on reconstrual rules, analogous to the proof
	of Lemma~\ref{lemma:ep-presentation}. The inductive steps for
	Rules~\ref{rule:relations} and~\ref{rule:top-bot-preservation} hold
	trivially. Rules~\ref{rule:conjunction-preservation}
	and~\ref{rule:quantifier-preservation} follow from \emph{Frobenius
	reciprocity} (Lemma A1.3.3 of~\cite{johnstone2002sketches}), but they can
	also be proven directly. This leaves
	Rules~\ref{rule:disjunction-preservation},~\ref{rule:context-duplication},
	and~\ref{rule:term-reduction}. We show the inductive steps for
	Rules~\ref{rule:disjunction-preservation}
	and~\ref{rule:context-duplication};~\ref{rule:term-reduction} follows a
	similar argument as the one for~\ref{rule:context-duplication}.
        
	(Rule~\ref{rule:disjunction-preservation}) Suppose \(\phi(\vec{x})
	\hookrightarrow \vec\sigma\) and \(\psi(\vec{y}) \hookrightarrow \vec{\tau}\)
	are \(T_1\)\-/formulae satisfying the lemma statement. The inductive step
	amounts to showing that \(\phi(\vec{x}) \vee \psi(\vec{y})\) also satisfies
	the lemma statement. Applying Rule~\ref{rule:disjunction-preservation} to the
	reconstrual \(F^\eq\), we have \(F^\eq\fb{\phi(\vec{x}) \vee
	\psi(\vec{y})}(\vec{s},\vec{t}) \equiv F^\eq\phi(\vec{s}) \vee
	F^\eq\psi(\vec{t})\). Now apply the induction hypothesis to each conjunct,
	and consider the following proof tree.
	\begin{prooftree}
		\AxiomC{}
		\RightLabel{\axiom[q^F_{\tau_1}]{Q1},...,\axiom[q^F_{\tau_n}]{Q1}}
		\UI$\fCenter \exists \pvec{t}'' \widearrow{q}^F_\tau(\pvec{t}'',\vec{t})$
		\RightLabel{\(\exists\)-elim.}
		\UI$\fCenter \exists \pvec{t}' \exists \pvec{t}''
		\left(\widearrow{q}^F_\tau(\pvec{t}'',\vec{t}) \wedge
		\ul{\widearrow{\dom}_{D^F_\tau}}(\pvec{t}'',\pvec{t}')\right)$
		\UI$ F^\eq\phi(\vec{s}) \fCenter F^\eq \phi(\vec{s}) \wedge \exists
		\pvec{t}' \exists \pvec{t}'' \left(\widearrow{q}^F_\tau(\pvec{t}'',\vec{t})
		\wedge \ul{\widearrow{\dom}_{D^F_\tau}}(\pvec{t}'',\pvec{t}')\right)$
	\end{prooftree}
	By \(\vee\)\-/introduction, \(\exists\)\-/introduction, and
	\(\exists\)\-/elimination, we arrive at the sequent%
	\begin{gather*}%
		F^\eq\phi(\vec{s}) \fCenter \exists \pvec{s}' \exists \pvec{s}'' \exists
		\pvec{t}' \exists \pvec{t}'' \Big(\pwrap{\ep_{T_1} F\phi(\pvec{s}') \vee
		\ep_{T_1} F\psi(\pvec{t}')} \wedge
		\widearrow{q}^F_\sigma(\pvec{s}'',\vec{s})\\ 
		\wedge \ul{\widearrow{\dom}_{D^F_\sigma}}(\pvec{s}'',\pvec{s}') \wedge
		\widearrow{q}^F_\tau(\pvec{t}'',\vec{t}) \wedge
		\ul{\widearrow{\dom}_{D^F_\tau}}(\pvec{t}'',\pvec{t}')\Big).
	\end{gather*}%
	There is an analogous proof if we replace \(F^\eq\phi(\vec{s})\) with
	\(F^\eq\psi(\vec{t})\), where we instead use \axiom[q^F_{\sigma_i}]{Q1}. By
	\(\vee\)\-/introduction, we deduce the forward sequent of the inductive step.
	The converse sequent is a consequence of \(\vee\)\-/introduction and
	\(\vee\)\-/elimination.

	(Rule~\ref{rule:context-duplication}) Suppose that
	\(\phi(\vec{y_1},\vec{x_1},\vec{y_2},\vec{x_2},\vec{y_3}) \hookrightarrow
	\vec{\tau_1},\vec{\sigma_1},\vec{\tau_2},\vec{\sigma_1},\vec{\tau_3}\) is a
	\(T_1\)-formula satisfying the lemma statement. We consider the case where
	\(\vec{x_2} \equiv \vec{x_1}\). Applying Rule~\ref{rule:context-duplication},
	\[F^\eq\fb{\phi(\vec{y_1},\vec{x_1},\vec{y_2},\vec{x_1},\vec{y_3})}
		(\vec{t_1},\vec{s_1},\vec{t_2},\vec{s_1},\vec{t_3}) \equiv
	F^\eq\phi(\vec{t_1},\vec{s_1},\vec{t_2},\vec{s_1},\vec{t_3}),\] which by the
	inductive hypothesis, is logically equivalent to
	\begin{gather*}%
		\exists \mvec{t_1}' \exists \mvec{t_1}'' \exists \mvec{s_1}' \exists
		\mvec{s_1}'' \exists \mvec{t_2}' \exists \mvec{t_2}'' \exists \mvec{s_2}'
		\exists \mvec{s_2}'' \exists \mvec{t_3}' \exists \mvec{t_3}''
		\Bigg(\ep_{T_1} F\phi(\mvec{t_1}',
			\mvec{s_1}',\mvec{t_2}',\mvec{s_2}',\mvec{t_3}')\\
		\wedge \bigwedge_{i=1}^3
		\left(\ul{\widearrow{\dom}_{D^F_{\tau_i}}}(\mvec{t_i}'', \mvec{t_i}')
			\wedge \widearrow{q}^F_{\tau_i}(\mvec{t_i}'',\vec{t_i})\right) \wedge
			\bigwedge_{i=1}^2
			\left(\ul{\widearrow{\dom}_{D^F_{\sigma_1}}}(\mvec{s_i}'', \mvec{s_i}')
		\wedge \widearrow{q}^F_{\sigma_1}(\mvec{s_i}'',\vec{s_1})\right)\Bigg).
	\end{gather*}%

	For ease of reference, set \(\vec{\upsilon} \equivdef \vec{\sigma_1}\). By
	\axiom[q^F_{\upsilon_i}]{Q2}, we have the sequent
	\[%
		\widearrow{q}^F_{\upsilon}(\mvec{s_1}'',\vec{s_1}) \wedge
		\widearrow{q}^F_{\upsilon}(\mvec{s_2}'',\vec{s_1}) \vdash \bigwedge_{i=1}^n
		\exists e_i^{\ul{E^F_{\upsilon_i}}}\left(\ul{\pi_1}(\ul{r_i}(e_i)) =
		s_{1i}'' \wedge \ul{\pi_2}(\ul{r_i}(e_i)) = s_{2i}''\right),
	\]%
	where we note that \(r_i : E^F_{\upsilon_i} \hookrightarrow D^F_{\upsilon_i}
	\times D^F_{\upsilon_i}\) satisfies \(\dom_{E^F_{\upsilon_i}} =
	\dom_{D^F_{\upsilon_i} \times D^F_{\upsilon_i}} \: r_i\). Since
	\(\dom_{D^F_{\upsilon_i} \times D^F_{\upsilon_i}}\) is a product morphism
	\(\dom_{D^F_{\upsilon_i}} \times \dom_{D^F_{\upsilon_i}}\), we further obtain
	the equations \(\dom_{D^F_{\upsilon_i}} \pi_j r_i =
	\pi^{\vec{\omega_i}}_{\vec{\omega_{ij}}} \dom_{E^F_{\upsilon_i}}\), where \(j
	\in \{1,2\}\) and \(\vec{\omega_i} \equivdef \vec{\omega_{i1}},
	\vec{\omega_{i2}} \equivdef F\upsilon_i,F\upsilon_i\). Thus, using IL2,
	\(=\)-elimination, the previous expression, and the cut rule, we find that
	the last conjunction of our expansion implies 
	\[%
		\bigwedge_{i=1}^n \exists
		e_i^{\ul{E^F_{\upsilon_i}}}\left(\ul{\pi^{\vec{\omega_i}}_{\vec{\omega_{i1}}}}
		\left(\ul{\dom_{E^F_{\upsilon_i}}}(e_i)\right) = s_{1i}' \wedge
	\ul{\pi^{\vec{\omega_i}}_{\vec{\omega_{i2}}}}
\left(\ul{\dom_{E^F_{\upsilon_i}}}(e_i)\right) = s_{2i}'\right),
	\]%
	which is logically equivalent to \(\bigwedge_{i=1}^n \ep_{T_1}
	E^F_{\upsilon_i} (s_{1i}',s_{2i}')\). By
	Lemma~\ref{lemma:super-subobject-lemma}, we see that
	\[%
		\ep_{T_1}
		F\phi(\mvec{t_1}',\mvec{s_1}',\mvec{t_2}',\mvec{s_2}',\mvec{t_3}') \wedge
		\bigwedge_{i=1}^n \ep_{T_1} E^F_{\upsilon_i}(s_{1i}',s_{2i}') \vdash
		\ep_{T_1}
		F\phi(\mvec{t_1}',\mvec{s_1}',\mvec{t_2}',\mvec{s_1}',\mvec{t_3}'),
	\]%
	allowing us to deduce the forward sequent of the induction. The converse
	sequent follows by \(\exists\)-introduction and \(\exists\)-elimination.
\end{proof}

\begin{proposition}\label{prop:eq-translation}
	\(F^\eq\) is an e.p.\ translation \(T_1 \to \T(\C(T_2))^\ex\).	
\end{proposition}
\begin{proof}
	We first show that \(F^\eq\) is a translation. Suppose \(\fb\phi \vdash
	\fb\psi\) is a provable sequent in \(T_1\), where \(\phi\) and \(\psi\) have
	domain \(\sigma_1,\ldots,\sigma_n\). By the preceding
	Lemma~\ref{lemma:eq-presentation} and deduction, proving \(F^\eq \phi \vdash
	F^\eq \psi\) is equivalent to proving the sequent \(\ep_{T_1}
	F\phi(\pvec{x}') \vdash \ep_{T_1} F\psi(\pvec{x}')\), which follows from
	Lemma~\ref{lemma:super-subobject-lemma}.

	We now show that \(F^\eq\) is e.p. Let \(\vec{\tau} \equivdef
	\ul{Q^F_\sigma},\ul{Q^F_\sigma}\). Then \(E^{F^\eq}_\sigma(x,y)\) is defined
	as
	\begin{gather*}%
		\exists x_1' \exists x_1'' \exists x_2' \exists x_2'' \left(\ep_{T_1}
		E^F_\sigma(x_1',x_2') \wedge
		\ul{\widearrow{\dom}_{D^F_\tau}}(x_1'',x_2'',x_1',x_2') \wedge
		\widearrow{q}^F_\tau(x_1'',x_2'',x,y)\right)\\
		\dashv\vdash \exists x_1''\exists x_2''
		\left(\ep_{T_1} E^F_\sigma\left(\ul{\dom_{D^F_{\tau_1}}}(x_1''),
		\ul{\dom_{D^F_{\tau_2}}}(x_2'')\right) \wedge 
		\widearrow{q}^F_\tau(x_1'',x_2'',x,y)\right).
	\end{gather*}
	As in the proof of Lemma~\ref{lemma:eq-presentation}, we use a logically
	equivalent presentation for \(\ep_{T_1} E^F_\sigma\), and then invoke IL2 and
	\(=\)\-/elimination. We additionally use IL3:
	\[%
		\exists e_1^{\ul{E^F_{\tau_1}}} \ul{\dom_{D^F_{\tau_1}}}
		(\ul{\pi_i}(\ul{r_1}(e_1))) = \ul{\dom_{D^F_{\tau_1}}}(x_i'') \vdash
		\exists e_1^{\ul{E^F_{\tau_1}}} \ul{\pi_i}(\ul{r_1}(e_1)) = x_i'',
	\]%
	where \(i \in \{1,2\}\). Apply \axiom[q^F_{\tau_1}]{Q2} and conclude via
	\(=\)\-/elimination.
\end{proof}

By the preceding Proposition~\ref{prop:eq-translation}, \(F^\eq: T_1 \to
\T(\C(T_2))^\ex\) induces a coherent functor \(\hw{F} \eqdef \C(F^\eq): \C(T_1)
\to \C\pwrap{\T(\C(T_2))^\ex} = \C(T_2)^\ex\). This will be the basis for the
definition of \(\X\) on 1\-/cells.

\begin{proposition}\label{prop:Fhat-quotient-diagram}
	Let \(F: T_1 \to T_2\) be a translation, and let \(\fb\phi \hookrightarrow
	\vec{\sigma}\) be a \(T_1\)\-/substitution class. The monomorphism
	\(\hw{F}\dom_{\fb\phi}: \hw{F}\fb\phi \hookrightarrow
	\fb{\ul{Q}^F_{\sigma_1},\ldots,\ul{Q}^F_{\sigma_n}}\) presents the subobject
	\(\exists_{\widearrow{q}^F_\sigma} F\fb\phi\) (where \(F\fb\phi\) is realized
	as a subobject of \(D^F_{\vec{\sigma}}\) by factoring the \(\C(T_2)\)
	morphism \(\dom_{F\phi}\)).
\end{proposition}
\begin{proof}
	Using the identification in Remark~\ref{remark:definition-of-I},
	\(D^F_{\vec\sigma}\) is identified with \(\fb{\ul{D^F_{\vec\sigma}}}\), and
	\(F\fb\phi\) is identified with \(\fb{\ul{F\phi}} \cong \fb{\ep_{T_2}
	F\fb\phi }\). Furthermore, \(D^F_{\vec\sigma}\) is isomorphic to the product
	\(D^F_{\sigma_1} \times \ldots \times D^F_{\sigma_n}\), which is identified
	with \(\fb{\ul{D^F_{\sigma_1}}} \times \ldots \times \fb{\ul{D^F_{\sigma_n}}}
	\cong \fb{\ep_{T_2} D^F_{\vec\sigma}}\), where the explicit isomorphism
	\(\varphi: \fb{\ep_{T_2} D^F_{\vec\sigma}} \xrightarrow\sim
	\fb{\ul{D^F_{\sigma_1}}} \times \ldots \times \fb{\ul{D^F_{\sigma_n}}}\) is
	presented by \( \varphi(\vec{x},\vec{y}) \equivdef
	\ul{\widearrow{\dom}_{D^F_\sigma}}(\vec{y},\vec{x})\). Consequently, the
	subobject \(F\fb\phi \hookrightarrow D^F_{\vec\sigma}\) is also presented by
	the monic \(e: \fb{\ep_{T_2} F\fb\phi} \hookrightarrow
	\fb{\ul{D^F_{\sigma_1}}} \times \ldots \times \fb{\ul{D^F_{\sigma_n}}},\)
	where \( e(\vec{x},\vec{y}) \equivdef \ep_{T_2} F\phi(\vec{x}) \wedge
	\ul{\widearrow{\dom}_{D^F_\sigma}}(\vec{y},\vec{x}) \). Thus
	\(\exists_{\widearrow{q}^F_\sigma} F\fb\phi\) is presented by the
	\(\T(\C(T_2))^\ex\)\-/formula%
	\begin{gather*}%
		\exists \vec{x'} \exists \vec{x''} e(\vec{x'},\vec{x''}) \wedge
		\widearrow{q}^F_\sigma(\vec{x''},\vec{x})\\
		\dashv\vdash \exists \vec{x'} \exists \vec{x''} \ep_{T_2}F\phi(\vec{x'})
		\wedge \ul{\widearrow{\dom}_{D^F_\sigma}}(\vec{x''},\vec{x'}) \wedge
		\widearrow{q}^F_\sigma(\vec{x''},\vec{x}).
	\end{gather*}
	\(\hw{F}\dom_\phi\) presents the subobject \(\exists \vec{t}\,
	\hw{F}\phi(\vec{t}) \wedge \vec{t} = \vec{x}\). By
	Lemma~\ref{lemma:eq-presentation}, this is logically equivalent to the
	preceding formula. In particular, the subobjects \(\hw{F}\dom_\phi\) and
	\(\exists_{\widearrow{q}^F_\sigma} F\fb\phi\) are presented by logically
	equivalent formulae in \(\T(\C(T_2))^\ex\), so they are the same subobject.
\end{proof}

\begin{proposition}\label{prop:Feq-preserves-homotopy}
	Let \(F,G: T_1 \to T_2\) be a pair of translations. Then \(F\) and \(G\) are
	homotopic if and only if \(F^\eq\) and \(G^\eq\) are homotopic.
\end{proposition}
\begin{proof}
	Suppose \(\chi: F \Rightarrow G\) is a homotopy t\-/map. We define a new
	family of \(\T(\C(T_2))^\ex\)\-/formulae:
	\begin{gather*}
		\chi^\eq_\sigma(x,y) \equivdef \exists {x'}^{\ul{D^F_\sigma}} \exists
		{x''}^{\ul{D^F_\sigma}} \exists {y'}^{\ul{D^G_\sigma}} \exists
		{y''}^{\ul{D^G_\sigma}} \Big(\ep_{T_1} \chi_\sigma(x',y') \wedge
		\ul{\dom_{D^F_\sigma}}(x'') = x' \\ 
		\wedge \ul{\dom_{D^G_\sigma}}(y'') = y' \wedge q^F_\sigma(x'') = x \wedge
		q^G_\sigma(y'') = y\Big).
	\end{gather*}
	Lemma~\ref{lemma:super-subobject-lemma} applied to the \axiom{TMi} axioms for
	\(\chi\), along with        Proposition~\ref{prop:eq-translation} and
	Lemmas~\ref{lemma:eq-presentation} and \ref{lemma:ep-presentation}, imply
	that \(\chi^\eq\) is a homotopy t\-/map \(F^\eq \Rightarrow G^\eq\).

	Conversely, suppose \(\eta: F^\eq \Rightarrow G^\eq\) is a homotopy t\-/map.
	Let \(\sigma\) be a \(T_1\)\-/sort. Then the component \(\eta_\sigma\) is a
	\(\T(\C(T_2))^\ex\)\-/substitution class with domain
	\(\ul{Q^F_\sigma},\ul{Q^G_\sigma}\). Moreover, since \(F^\eq\) and \(G^\eq\)
	are e.p.~(Proposition~\ref{prop:eq-translation}), \(\eta_\sigma\) is a
	\(\T(\C(T_2))^\ex\)\-/definable isomorphism \(Q^F_\sigma \to Q^G_\sigma\).
	Consider the `preimage' of \(\eta_\sigma\):
	\begin{gather*}%
		\tilde{\eta}_\sigma\pwrap{\pvec{x}',\pvec{y}'} \equivdef \exists
		x^{\ul{Q^F_\sigma}} \exists {x''}^{\ul{D^F_\sigma}} \exists
		y^{\ul{Q^G_\sigma}} \exists {y''}^{\ul{D^G_\sigma}}\\
		\bigwedge_{i=1}^n \pwrap{%
			\ul{\pi^{\vec{\tau}}_{\tau_i}} \pwrap{%
				\ul{\dom_{D^F_\sigma}}(x'')
			} = x_i' 
		} \wedge \bigwedge_{j=1}^m \pwrap{
			\ul{\pi^{\vec\upsilon}_{\upsilon_j}}\pwrap{%
				\ul{\dom_{D^G_\sigma}}(y'')
			} = y_j'
		}\\
		\wedge q^F_\sigma(x'') = x \wedge q^G_\sigma(y'') = y \wedge
		\eta_\sigma(x,y),
	\end{gather*}%
	where \(F\sigma \equivdef \tau_1,\ldots,\tau_n\) and \(G\sigma \equivdef
	\upsilon_1,\ldots,\upsilon_m\). Since \(\tilde{\eta}_\sigma\) is a
	\(\T(\C(T_2))^\ex\)\-/formula with domain
	\(\ul{\vec{\tau}},\ul{\vec{\upsilon}}\), it determines a unique subobject of
	\(\fb{\ul{\vec{\tau}}} \times \fb{\ul{\vec{\upsilon}}}\) in
	\(\T(\C(T_2))^\ex\). Note that the collection \(\{\tilde{\eta}_\sigma\}\)
	satisfies the \axiom{TMi} axioms. Specifically, \(\tilde{\eta}\) defines a
	homotopy \(\ep_{T_2} F \Rightarrow \ep_{T_2} G\), where \(\ep_{T_2} F\) and
	\(\ep_{T_2} G\) are interpreted as translations \(T_1 \to \T(\C(T_2))^\ex\).
	We convert this to a homotopy \(F \Rightarrow G\) using \axiom{EC1} of
	Proposition~\ref{prop:exact-completion}, which indicates that \(\C(T_2)
	\hookrightarrow \C(T_2)^\ex\) is full on subobjects. Since
	\(\fb{\ul{\vec{\tau}}} \times \fb{\ul{\vec{\upsilon}}}\) is an object of
	\(\C(T_2)\), the subobject presented by \(\tilde{\eta}_\sigma\) is also
	presented by a monomorphism \(\alpha: A \hookrightarrow \fb{\ul{\vec{\tau}}}
	\times \fb{\ul{\vec{\upsilon}}} = \fb{\ul{\vec{\tau}},\ul{\vec{\upsilon}}}\)
	in \(\C(T_2)\). Define \(\chi_\sigma\) to be the substitution class presented
	by
	\[%
		\chi_\sigma(x',y') \equivdef \exists a^{\Dom A} \alpha(a,x',y').
	\]%
	This ensures that \(\dom_{\chi_\sigma}\) presents the same subobject as
	\(\alpha\). In particular, \(\ep_{T_2} \chi_\sigma\) is logically equivalent
	to \(\tilde{\eta}_\sigma\) in \(\T(\C(T_2))^\ex\). Using
	Lemma~\ref{lemma:super-subobject-lemma} and \(\{\tilde{\eta}_\sigma\}\), we
	deduce that the collection \(\chi \eqdef \{\overline{\chi_\sigma}\}\)
	satisfies the \axiom{TMi} axioms for a homotopy \(\chi: F \Rightarrow G\), as
	desired.
\end{proof}

We are now set to completely define the functor \(\X: \h\CTh_0 \to \h\Exact\).

\begin{definition}\label{def:X-functor}
	Given a theory \(T\), set \(\X(T) \eqdef \C(T)^\ex\). Given a homotopy class
	of a translation \(\fbb{F}: T_1 \to T_2\), let \(\hw{F}: \C(T_1) \to
	\C(T_2)^\ex\) be the coherent functor \(\C(F^\eq)\). The 1\-/cell \(\X\fbb{F}:
	\X(T_1) \to \X(T_2)\) is defined to be natural isomorphism class
	\(\fbb{\hw{F}^\ex}\), where \(\hw{F}^\ex: \C(T_1)^\ex \to \C(T_2)^\ex\) is the
	lift of \(\hw{F}\) using \axiom{EC3} of
	Proposition~\ref{prop:exact-completion}.
\end{definition}

By Proposition~\ref{prop:Feq-preserves-homotopy}, the map \(F \mapsto F^\eq\)
descends to a well\-/defined map on homotopy classes. By
Theorem~\ref{theorem:CThEqCoh-Correspondence}, the map \(F^\eq \mapsto
\hw{F}\) also descends to a map on homotopy classes. Finally, by
Proposition~\ref{prop:exact-reflects-homotopy}, \(\hw{F} \mapsto \hw{F}^\ex\)
descends to a map on natural isomorphism classes. Therefore the functor \(\X\)
is well\-/defined.

\subsection{Functoriality and Equivalence}

We have shown that the map \(\X: \h\CTh_0 \to \h\Exact\) is well-defined. Now we
need to show that it is a functor.

\begin{proposition}\label{prop:X-identities}
	Let \(T\) be a coherent theory. Then \(\X\fbb{1_T} = \fbb{1_{\X(T)}}\).
\end{proposition}

\begin{proof}
	We need to find a natural isomorphism \(\hw{1_T}^\ex \Rightarrow 1_{\X(T)}\).
	Let \(I: \C(T) \to \X(T)\) be the inclusion functor of
	Proposition~\ref{prop:exact-completion}. By
	Corollary~\ref{cor:exact-preserves-homotopy}, it suffices to find a natural
	isomorphism \(\hw{1_T}^\ex \circ I \Rightarrow I\). Furthermore EC3 stipulates
	\(\hw{1_T}^\ex \circ I = \hw{1_T}\). Given an object \(\fb\phi \hookrightarrow
	\vec{\sigma}\) of \(\C(T)\), \(\hw{1_T}\fb\phi\) fits into a (trivial) quotient
	diagram%
	\[%
		\fb{\vec{x} =_{\vec{\sigma}} \vec{y} \wedge \phi(\vec{x}) \wedge
		\phi(\vec{y})} \rightrightarrows 1_T \fb\phi \twoheadrightarrow
		\hw{1_T}\fb\phi.
	\]%
	Furthermore \(1_T \fb\phi\) is logically equivalent to \(\fb\phi\)
	(Proposition~\ref{prop:reconstrual-log-eq-to-identity}), so
	\(\hw{1_T}\fb\phi\) is the quotient of \(\fb\phi\) along the diagonal
	\(\fb\phi \hookrightarrow \fb\phi \times \fb\phi\). As \(\fb\phi\) satisfies the
	universal property of this quotient, \(\hw{1_T}\fb\phi\) is isomorphic to
	\(\fb\phi\), and this isomorphism is natural in \(\fb\phi\) since colimits are
	unique up to natural isomorphism.
\end{proof}

\begin{proposition}\label{prop:X-composition}
	Given a composable pair of translations \(F: T_1 \to T_2\) and \(G: T_2 \to
	T_3\), \(\X\fbb{G}\X\fbb{F} = \X\fbb{GF}\).
\end{proposition}

\begin{proof}
	It suffices to find a natural isomorphism \(\hw{G}^\ex\hw{F}^\ex \Rightarrow
	\hw{GF}^\ex\).  Corollary~\ref{cor:exact-preserves-homotopy} reduces this task
	to finding a natural isomorphism \(\kappa: \hw{G}^\ex\hw{F}^\ex I \Rightarrow
	\hw{GF}^\ex I\), where \(I : \C(T_1) \to \C(T_1)^\ex\) is the inclusion
	functor of Proposition~\ref{prop:exact-completion}. Since \(I\) is fully
	faithful, we identify \(\C(T_1)\) with a subcategory of \(\C(T_1)^\ex\). Let
	\(A \hookrightarrow \vec{\sigma}\) be an object of \(\C(T_1)\). The component of
	\(\kappa\) along \(A\) is a morphism \(\kappa_A: \hw{G}^\ex\hw{F}A \to
	\hw{GF}A\). Proposition~\ref{prop:Fhat-quotient-diagram} implies that
	\(\hw{GF}A\) is a quotient \(\exists_{\widearrow{q}^{GF}_\sigma} GF A\), fitting into a
	diagram%
	\[%
		\restr{E^{GF}_{\vec{\sigma}}}{GF A} \rightrightarrows GF A
		\xrightarrow{\restr{\widearrow{q}^{GF}_\sigma}{GF A}} \hw{GF} A \cong
		\exists_{\widearrow{q}^{GF}_\sigma} GF A.
	\]%
	Since colimits are unique up to natural isomorphism, we can find
	\(\kappa_A\)---such that \(\kappa_A\) is natural in \(A\)---by showing that
	\(\hw{G}^\ex\hw{F}A\) is also a quotient of \(GF A/\restr{E^{GF}_{\vec{\sigma}}}{GF A}\).

	First we need to find a morphism \(Q: GF A \to \hw{G}^\ex\hw{F}A\) to serve
	as the quotient morphism. To find it, we split \(Q\) into two components: a
	\(\hw{G}^\ex\) piece and a \(\widearrow{q}^G_{\tau}\) piece for \(\vec{\tau}
	\equivdef F\vec{\sigma}\). We begin by noting that \(\hw{F} A\) is a
	quotient:%
	\[%
		\restr{E^F_{\vec{\sigma}}}{F A} \rightrightarrows F A
		\xrightarrow{\restr{\widearrow{q}^F_\sigma}{F A}} \hw{F} A.
	\]%
	Since \(\hw{G}^\ex\) preserves quotients, we have another quotient diagram in
	\(\C(T_3)^\ex\)%
	\[%
		\begin{tikzcd}
			\hw{G}\pwrap{\restr{E^F_{\vec{\sigma}}}{F A}} \arrow[r,"p_1'", shift right,swap]
			\arrow[r, "p_2'", shift left] & \hw{G}F A \arrow[r, "\hw{G}^\ex q"] &
			\hw{G}^\ex\hw{F}A,
		\end{tikzcd}
	\]%
	where \(q\) is an abbreviation for \(\restr{\widearrow{q}^F_\sigma}{F A}\)
	and \(p_1',p_2'\) are the pullback's projection morphisms. \(\hw{G}\) is a
	coherent functor, so it preserves pullbacks; hence
	\(\hw{G}\pwrap{\restr{E^F_{\vec{\sigma}}}{F A}} \cong
	\restr{\hw{G}E^F_{\vec{\sigma}}}{\hw{G}F A}\). Then the last quotient diagram
	indicates that we have an effective epimorphism \(\hw{G}^{\ex}q : \hw{G} F A
	\to \hw{G}^\ex\hw{F} A\). We now find a morphism \(GF A \to \hw{G} F A\).
	Similar to the case for \(\hw{F}A\), \(\hw{G}F A\) fits into a quotient%
	\[%
		\begin{tikzcd}
			\restr{E^G_{\vec{\tau}}}{GF A} \arrow[r,"p_2", shift left] \arrow[r,"p_1",
			shift right, swap] & GF A \arrow[r,"\restr{\widearrow{q}^G_{\tau}}{GF A}"
			{yshift=1pt}] & \hw{G}
			F A.
		\end{tikzcd}
	\]%
	These define the needed components. Set \(Q \eqdef
	\hw{G}^\ex\pwrap{\restr{\widearrow{q}^F_{\sigma}}{F
	A}}\circ\restr{\widearrow{q}^G_{\tau}}{GF A}\).

	We now verify that \(Q: GF A \to \hw{G}^\ex F A\) satisfies the same
	universal property as the quotient \(\restr{\widearrow{q}^{GF}_\sigma}{GF A}:
	GF A \to \hw{GF} A\).

	Suppose we have a morphism \(f: GF A \to Z\) in \(\C(T_3)^\ex = \X(T_3)\)
	which coequalizes \(\restr{E^{GF}_{\vec{\sigma}}}{GF A}\):%
	\begin{equation}\label{eq:EGF-f-Z-diagram}%
		\restr{E^{GF}_{\vec{\sigma}}}{GF A} \doublerightarrow{p_2''}{p_1''} GF A
		\xrightarrow{f} Z.
	\end{equation}%
	We need to prove that there exists a unique morphism \(f'': \hw{G}^\ex\hw{F}
	A \to Z\) which factors \(f\) through \(Q\):%
	\[%
		\begin{tikzcd}
			\restr{E^{GF}_{\vec{\sigma}}}{GF A} \arrow[r,"p_2''",shift left]
			\arrow[r,"p_1''",shift right,swap] &
			GF A \arrow[r, "Q"] \arrow[rr,bend right=30,"f"] & \hw{G}^\ex\hw{F} A
			\arrow[r,dashed,"\exists! f''"] & Z.
		\end{tikzcd}
	\]%
	Recall Lemma~\ref{lemma:t-map-lemma-2}: \(E^G_{\vec{\tau}}(x,y) \wedge
	D^{GF}_{\vec{\sigma}}(x) \vdash E^{GF}_{\vec{\sigma}}(x,y)\). Since \(G\) and
	\(F\) are translations, \(GF A \vdash D^{GF}_{\vec{\sigma}}\), so the
	preceding sequent implies \(E^G_{\vec{\tau}}(x,y) \wedge GF A(x) \vdash
	E^{GF}_{\vec{\sigma}}(x,y)\). Therefore as subobjects of \(\Sub\pwrap{GF A
	\times GF A}\), \(\restr{E^G_{\vec{\tau}}}{GF A}\) factors through
	\(\restr{E^{GF}_{\vec{\sigma}}}{GF A}\), so there exists a monomorphism
	\(\iota\) allowing us to attach \(\restr{E^G_{\vec{\tau}}}{GF A}\) to
	Diagram~\ref{eq:EGF-f-Z-diagram}.
	\[%
		\begin{tikzcd}[sep=small]
			\restr{E^G_{\vec{\tau}}}{GF A}
			\arrow[rr,"p_2"]\arrow[dd,"p_1",swap]\arrow[dr,"\iota",dashed,hook] & & GF
			A \arrow[dd,"f"]\\
			& \restr{E^{GF}_{\vec{\sigma}}}{GF A}\arrow[dl,"p_1''",swap]
			\arrow[ur,"p_2''"] & \\
			GF A \arrow[rr,"f"] & & Z
		\end{tikzcd}
	\]%
	Commutativity of the outer square implies \(f\) also coequalizes
	\(\restr{E^G_{\vec{\tau}}}{GF A}\). Therefore the universal property of the
	quotient \(\hw{G}F A\) implies the existence of a \emph{unique} morphism
	\(f': \hw{G} F A \to Z\) such that \(f = f' \circ
	\restr{\widearrow{q}^G_{\tau}}{GF A}\).
	\[%
		\begin{tikzcd}
			\restr{E^G_{\vec{\tau}}}{GF A} \arrow[r,shift left,"p_2"] \arrow[r, shift
			right,"p_1",swap] & GF A \arrow[r, "\restr{\widearrow{q}^G_{\tau}}{GF A}"
			{yshift=1pt}]
			\arrow[rr,bend right=30,"f"] & \hw{G}F A \arrow[r,dashed,"\exists! f'"] &
			Z
		\end{tikzcd}
	\]%
	In particular \(f p_1 = f p_2\). The next step is to show that \(f'\)
	coequalizes \(p_1'\) and \(p_2'\) from earlier.
	\(\restr{E^F_{\vec{\sigma}}}{F A}\) is an object of \(\C(T_2)\), thus
	\(\hw{G}\restr{E^F_{\vec{\sigma}}}{\hw{G}F A}\) is itself a quotient. Let
	\begin{gather*}
		R \eqdef \restr{\pwrap{%
			\restr{E^G_{\vec\tau}}{GF A} \times \restr{E^G_{\vec\tau}}{GF A}
		}}{%
			\restr{E^{GF}_{\vec\sigma}}{GF A}
		},\\%
		\tilde{q} \eqdef \restr{\widearrow{q}^G_\tau}{GF A} \times
		\restr{\widearrow{q}^G_\tau}{GF A},\quad%
		\tilde{p}_i \eqdef \restr{\pwrap{p_i \times p_i}}{%
			\restr{E^{GF}_{\vec\sigma}}{GF A}
		}.%
	\end{gather*}
	Semantically \(R\) represents pairs of contexts for
	\(\restr{E^{GF}_{\vec\sigma}}{GF A}\) which are equivalent modulo
	\(\restr{E^G_{\vec\tau}}{GF A}\). The quotient diagram for
	\(\hw{G}\restr{E^F_{\vec{\sigma}}}{\hw{G}F A}\) is
	\begin{equation}\label{eq:R-quotient}
		\begin{tikzcd}[ampersand replacement=\&]
			R \arrow[r,shift left,"{\tilde{p}_2}"] \arrow[r, shift right, swap,
			"{\tilde{p}_1}"] \& {\restr{E^{GF}_{\vec\sigma}}{GF A}} \arrow[r, two
			heads, "{\restr{\tilde{q}}{E^{GF}_{\vec\sigma}}}"] \&
			{\hw{G}\restr{E^F_{\vec{\sigma}}}{\hw{G}F A}}.
		\end{tikzcd}
	\end{equation}
	Since \(f' p_1'\) and \(f' p_2'\) are morphisms out of the quotient
	\(\hw{G}\restr{E^F_{\vec{\sigma}}}{\hw{G}F A}\), \(f' p_1'\) and \(f' p_2'\)
	correspond uniquely to a pair of squares, respectively:
	\[\arraycolsep=-0.2em%
		\begin{array}{cc}
			\begin{tikzcd}
				R \arrow[r,"\tilde{p}_2"]\arrow[d,"\tilde{p}_1",swap] &
				\restr{E^{GF}_{\vec{\sigma}}}{GF A} \arrow[d,"f' p_1'
				\restr{\tilde{q}}{E^{GF}_{\vec\sigma}}"]\\
				\restr{E^{GF}_{\vec{\sigma}}}{GF A} \arrow[r,"f' p_1'
				\restr{\tilde{q}}{E^{GF}_{\vec\sigma}}",swap] & Z,
			\end{tikzcd} &
			\begin{tikzcd}
				R \arrow[r,"\tilde{p}_2"]\arrow[d,"\tilde{p}_1",swap] &
				\restr{E^{GF}_{\vec{\sigma}}}{GF A} \arrow[d,"f' p_2'
				\restr{\tilde{q}}{E^{GF}_{\vec\sigma}}"]\\
				\restr{E^{GF}_{\vec{\sigma}}}{GF A} \arrow[r,"f' p_2'
				\restr{\tilde{q}}{E^{GF}_{\vec\sigma}}",swap] & Z.
			\end{tikzcd}
		\end{array}
	\]%

	Since \(\tilde{q}\) is a product of quotient morphisms and \(p_i'\) the
	projections associated to a pullback, we can permute \(p_i'\circ
	\restr{\tilde{q}}{E^{GF}_{\vec\sigma}}\) to a projection followed by a
	quotient. To that end, we note that the following three diagrams commute,
	where the first commutes because \(\restr{\tilde{q}}{E^{GF}_{\vec\sigma}}\)
	is a coimage.
	\[\setlength{\arraycolsep}{-0.2em}%
		\begin{array}{@{}c@{}c@{}c@{}}
			\begin{tikzcd}[ampersand replacement=\&, column sep=small]
				\restr{E^{GF}_{\vec\sigma}}{GF A} \arrow[r,hook]
				\arrow[d,"{\restr{\tilde{q}}{E^{GF}_{\vec\sigma}}}",swap, two heads] \&
				\pwrap{GF A}^2 \arrow[d, two heads, "\tilde{q}"] \\
				\restr{\hw{G}E^F_{\vec\sigma}}{\hw{G} F A} \arrow[r, hook] \&
				\pwrap{\hw{G} F A}^2
			\end{tikzcd} &
			\begin{tikzcd}[ampersand replacement=\&, column sep=small]
				\restr{\hw{G}E^F_{\vec\sigma}}{\hw{G}F A} \arrow[r,hook]
				\arrow[dr,"p_i", swap] \& \pwrap{\hw{G}F A}^2 \arrow[d,two
				heads,"\pi_i'"] \\
				\& \hw{G}F A
			\end{tikzcd} &
			\begin{tikzcd}[ampersand replacement=\&, column sep=small]
				\restr{E^{GF}_{\vec\sigma}}{GF A} \arrow[r,hook] \arrow[dr,"p_i''",
				swap] \& \pwrap{GF A}^2 \arrow[d,"\pi_i"] \\
				\& GF A
			\end{tikzcd}
		\end{array}
	\]%
	Since \(\tilde{q}\) is a product, \(\pi_i' \circ \tilde{q} =
	\restr{q^G_{\vec\tau}}{GF A} \circ \pi_i\). Stitching the above diagrams together,
	this equation implies \(p_i' \circ \restr{\tilde{q}}{E^{GF}_{\vec\sigma}} =
	\restr{q^G_{\vec\tau}}{GF A} \circ p_i''\).

	Swapping the projection and quotient morphisms using the above argument
	yields the equation
	\[%
		f' \circ p_i' \circ \restr{\tilde{q}}{E^{GF}_{\vec\sigma}} = f' \circ
		\restr{q^G_{\vec\tau}}{GF A} \circ p_i'' = f \circ p_i''.
	\]%
	Diagram~\ref{eq:EGF-f-Z-diagram} implies \(f \circ p_1'' = f \circ p_2''\),
	therefore the universal property of Diagram~\ref{eq:R-quotient} implies \(f'
	p_1' = f' p_2'\). Thus \(f'\) coequalizes \(p_1'\) and \(p_2'\), so the
	universal property of \(\hw{G}^\ex\pwrap{\restr{\widearrow{q}^F_\sigma}{F
	A}}\) implies that there exists a \emph{unique} morphism \(f'' :
	\hw{G}^\ex\hw{F} A \to Z\) fitting into the diagram
	\[%
		\begin{tikzcd}[sep=large]
			\hw{G}\restr{E^F_{\vec{\sigma}}}{\hw{G}F A} \arrow[r,shift left]
			\arrow[r, shift right] & \hw{G}F A
			\arrow[r,"\hw{G}^\ex\pwrap{\restr{\widearrow{q}^F_\sigma}{F A}}"
			{yshift=1pt}] \arrow[rr,"f'", bend right=30] & \hw{G}^\ex \hw{F} A
			\arrow[r,"\exists! f''", dashed] & Z.
		\end{tikzcd}
	\]%
	In particular, \(f = f'' \circ
	\hw{G}^\ex\pwrap{\restr{\widearrow{q}^F_\sigma}{F A}} \circ
	\restr{\widearrow{q}^G_{\tau}}{G F A} = f'' \circ Q\). \(f''\) is uniquely
	determined by \(f'\), which is uniquely determined by \(f\); therefore
	\(f''\) is uniquely determined by \(f\). This is the same universal property
	(shown below) as the quotient \(\hw{GF} A\), as desired.%
	\[%
		\begin{tikzcd}
			\restr{E^{GF}_{\vec{\sigma}}}{GF A} \arrow[r, shift left] \arrow[r, shift
			right] & GF A \arrow[r,"Q"] \arrow[rr,"f", bend right=30] & \hw{G}^\ex
			\hw{F} A \arrow[r,dashed,"\exists! f''"] & Z
		\end{tikzcd}
	\]%
	Since \(\hw{GF} A\) and \(\hw{G}^\ex\hw{F} A\) satisfy the same universal
	property, there exists an isomorphism \(\hw{G}^\ex\hw{F} A \to \hw{GF} A\).
	Setting \(\kappa_A\) to be this isomorphism yields the desired natural
	isomorphism.
\end{proof}

\begin{corollary}\label{prop:X-functorial}
	The operation \(\X: \h\CTh_0 \to \h\Exact\) is a functor.
\end{corollary}

We are ready to complete the second main theorem.

\begin{theorem}\label{theorem:bi-interpretability-equivalence}
	The functor \(\X: \h\CTh_0 \to \h\Exact\) is fully faithful and essentially
	surjective, so it is an equivalence of categories.
\end{theorem}

\begin{proof}
	We first show that \(\X\) is faithful. Suppose \(\X\fbb{F} = \X\fbb{G}\), where
	\(F,G: T_1 \to T_2\) are a pair of translations. Then \(\hw{F}^\ex\) and
	\(\hw{G}^\ex\) are naturally isomorphic. By
	Proposition~\ref{prop:exact-reflects-homotopy}, this implies \(\hw{F} \simeq
	\hw{G}\). By Proposition~\ref{prop:Feq-preserves-homotopy} and
	Theorem~\ref{theorem:CThEqCoh-Correspondence}, this implies \(F \simeq G\),
	so \(\fbb{F} = \fbb{G}\); hence \(\X\) is faithful. 
	
	Now we show that \(\X\) is essentially surjective. Let \(C\) be a
	Barr\-/exact coherent category. Then \(T \eqdef \T(C)\) is a coherent theory.
	By Theorem~\ref{theorem:CThEqCoh-Correspondence}, \(\C(T)\) is equivalent to
	\(C\), so \(\C(T)\) is Barr\-/exact. In particular, \(\X(T) = \C(T)^\ex\) and
	\(\C(T)\) are equivalent categories, so \(\X(T)\) and \(C\) are isomorphic in
	\(\h\Exact\).

	Finally we show that \(\X\) is full. Let \(\fbb{\Eff} : \X(T_1) \to \X(T_2)\)
	be a homotopy class presented by a coherent functor \(\Eff: \C(T_1)^\ex \to
	\C(T_2)^\ex\). Let \(I: \C(T_1) \to \C(T_1)^\ex\) be the inclusion functor,
	and let \(T\) be the conservative extension \(\T(\C(T_2))^\ex\) of
	\(\T\C(T_2)\) such that \(\C(T) = \C(T_2)^\ex\).  Given the coherent functor
	\(\Eff \circ I : \C(T_1) \to \C(T_2)^\ex\), by
	Theorem~\ref{theorem:CThEqCoh-Correspondence}, there exists an e.p.\
	translation \(F: T_1 \to T\) such that \(\C(F)\) is naturally isomorphic to
	\(\Eff \circ I\). The rest of this proof will construct a translation \(T \to
	T_2\) which we can compose with \(F\) to get the desired translation \(T_1
	\to T_2\). 
	
	Recall that the proof of Theorem~\ref{theorem:CThEqCoh-Correspondence}
	constructed an e.p.\ translation \(\gamma: \T\C(T_2) \to T_2\). Consequently
	we will find a translation \(\rho: T \to \T\C(T_2)\). Once we find \(\rho\),
	we will show that this choice yields \(\X(\gamma\rho) = \fbb{ 1_{\X(T_2)} }\);
	hence by Proposition~\ref{prop:X-composition} we conclude
	\(\X\pwrap{\gamma\rho F} = \X(\gamma\rho)\X(F) = \fbb{ 1_{\X(T_2)} }
	\fbb{\C(F)^\ex} = \fbb\Eff\), as desired. We need to define a reconstrual for
	\(\rho\). \(T\) is an extension of \(\T\C(T_2)\): given any symbol \(\ast\) in
	the signature of \(T\) which is contained in the signature of \(\T\C(T_2)\),
	set \(\rho(\ast) \equivdef 1_{\T\C(T_2)}(\ast)\). The only remaining symbols
	in the signature of \(T\) come from quotients. Given a congruence \(r: R
	\hookrightarrow A \times A\) in \(\C(T_2)\), define:%
	\[%
		\begin{array}{ccc}
			\rho\pwrap{\ul A/\ul R} \equivdef \ul A,&
			E^\rho_{{\ul A}/{\ul R}} \equivdef \fb{ \exists t^{\ul R}
			\bigwedge_{i=1,2} \ul{\pi_i}\pwrap{\ul{r}(t)} =_{\ul A} x_i},&
			\rho\pwrap{q_{R}} \equivdef E^\rho_{\ul A/ \ul R}.
		\end{array}
	\]%
	Note that this definition of \(\rho\) sends the axioms \axiom[R]{Q1} and
	\axiom[R]{Q2} to provable sequents in \(\T\C(T_2)\). Since \(1_{\T\C(T_2)}\)
	is a translation, and \(\rho\) is identical to \(1_{\T\C(T_2)}\) when
	restricting to the signature of \(\T\C(T_2)\), \(\rho\) sends all the axioms
	of \(T\) to provable sequents in \(\T\C(T_2)\), so \(\rho\) is a translation. 
	
	The next step is showing \(\X\pwrap{\gamma\rho} = \fbb{ 1_{\X(T_2)} }\), or
	equivalently that \(\hw{\gamma\rho}^\ex \simeq 1_{\X(T_2)}\). Since \(\C(T)\)
	is Barr\-/exact, \(\C(T)\) is its own exact completion in the sense of
	Proposition~\ref{prop:exact-completion}; therefore \(\hw{\gamma\rho}^\ex\) is
	naturally isomorphic to \(\hw{\gamma\rho}: \C(T) \to \C(T_2)^\ex\). Recall
	\(\C(T_2)^\ex\) was defined to be \(\C(T)\), so \(\hw{\gamma\rho}\) is
	a functor from \(\C(T)\) to itself; composing \(\hw{\gamma\rho}\) with the
	inclusion \(I_2: \C(T_2) \to \C(T)\) and invoking \axiom{EC2} shows
	\(\pwrap{\hw{\gamma\rho} \circ I_2}^\ex\) is naturally isomorphic to
	\(\hw{\gamma\rho}\). By Proposition~\ref{prop:exact-reflects-homotopy}, this
	implies that \(\hw{\gamma\rho} \simeq 1_{\X(T_2)}\) if and only if
	\(\hw{\gamma\rho} \circ I_2\) is naturally isomorphic to \(I_2\). Let \(A
	\equiv \fb\phi \) be an arbitrary object of \(\C(T_2)\).  Then
	\(\pwrap{\hw{\gamma\rho} \circ I_2} A \equiv \hw{\gamma\rho}\fb{\ul A}\), and the
	latter is the quotient \(\hw{\gamma\rho}\fb{\ul A} =
	\exists_{q} D^{\gamma\rho}_{\ul A}\), where \(q\) stands for 
        \(\widearrow{q}^{\gamma\rho}_{\ul A}\). However
	\(E^{\gamma\rho}_{\ul A} \equiv \gamma E^{\rho}_{\ul A} \equiv \gamma
	E^{1_{\T\C(T_2)}}_{\ul A} \equiv E^\gamma_{\ul A}\). Recalling the definition
	of \(\gamma\), this implies that%
	\[%
			E^{\gamma\rho}_{\ul A}(x,y) \equiv E^\gamma_{\ul A}(x,y) \dashv\vdash
			\phi(x) \wedge x = y.
	\]%
	Thus the congruence \(E^{\gamma\rho}_{\ul A} \hookrightarrow \fb\phi \times
	\fb\phi\) presents the same subobject of \(\fb\phi \times \fb\phi\) as the
	diagonal \(\fb\phi \hookrightarrow \fb\phi \times \fb\phi\). In particular
	\(\exists_q D^{\gamma\rho}_{\ul A}\) is isomorphic to \(I_2 \fb\phi \equiv I_2
	A\).  Furthermore, quotients are unique up to natural isomorphism, so
	\(\exists_q D^{\gamma\rho}_{\ul A}\) is \emph{naturally} isomorphic to \(I_2
	A\). We conclude that \(\hw{\gamma\rho}\circ I_2\) is naturally isomorphic to
	\(I_2\). This completes the proof that \(\X(\gamma\rho) = \fbb{1_{\X(T_2)}}\).
	Consequently, \(\gamma\rho F: T_1 \to T_2\) is a translation such that
	\(\X(\gamma\rho F) = \X(\gamma\rho)\X(F) = \fbb{1_{\X(T_2)}}\fbb{\Eff} =
	\fbb{\Eff}\), so \(\X\) is full. Hence \(\X\) is fully
	faithful and essentially surjective, so it is an equivalence of categories.
\end{proof}

Since homotopy equivalence in \(\CTh_0\) is bi\-/interpretability and homotopy
equivalence in \(\Exact\) is equivalence of categories, we have proven our
last main theorem.

\begin{corollary*}[Theorem~\ref{theorem:BiInterpretability-result}]
	Two small coherent theories \(T_1\) and \(T_2\) are bi\-/interpretable if and
	only if \(\C(T_1)^\ex\) and \(\C(T_2)^\ex\) are equivalent categories.
\end{corollary*}

\section{Conclusions}

The three main theorems of this paper, Theorems~\ref{theorem:cth0},
\ref{theorem:CThEqCoh-Correspondence},
and~\ref{theorem:bi-interpretability-equivalence}, are novel results that
elucidate the bicategorical structure of coherent theories. Having established
that \(\CThEq\) is a bicategory, a natural question to ask is what kinds of
2-limits and 2-colimits exist in \(\CThEq\). Using
Theorem~\ref{theorem:CThEqCoh-Correspondence}, we can answer this question
using results for \(\Coh\). In Theorem~4.9 of~\cite{kanalas2021}, it is shown
that \(\Coh\) admits weak colimits. Since biequivalences preserve weak
colimits, \(\CThEq\) also admits weak colimits.

\begin{theorem}\label{theorem:CThEq-Cocomplete}
	\(\CThEq\) has weak colimits.
\end{theorem}

This is a generalization of work by~\cite{Visser_categoriesof}, where it is
shown that a single\-/sorted analogue of \(\CTh_0\) admits coproducts (called
\emph{sum theories} in~\cite{Visser_categoriesof}).
Theorem~\ref{theorem:CThEq-Cocomplete} shows that the same is true for
many\-/sorted theories. Furthermore, Theorem~\ref{theorem:CThEq-Cocomplete}
implies additional constructions of coherent theories exist, including quotients
and non-disjoint unions, i.e., pushouts, at least in a bicategorical sense.

On the other hand,~\cite{kanalas2021} presents results for the existence of
certain weak limits and homotopy limits, e.g., Theorem 4.25
of~\cite{kanalas2021}. Theorem~\ref{theorem:CThEqCoh-Correspondence} allows
these existence theorems to be ported over to coherent theories, as long as one
is careful to interpret what a \emph{homotopy limit} (in the abstract
homotopy-theoretic sense) of coherent theories is.

There is already research in the nexus of mathematical logic and homotopy
theory. Campion, Cousins, and Ye~\cite{campion_cousins_ye_2021} associate to
any theory a topological space determined by the theory's category of models.
This association exhausts all homotopy types in the sense that any homotopy
type of a space is presented by an abstract elementary class. Our work
demonstrates that theories have an \emph{intrinsic} homotopy\-/theoretic
structure, setting the foundation for explaining how these homotopy-theoretic
ideas may arise.

\subsection{Biequivalences in Related Categories}

The proofs for the three main theorems can be extended to logic fragments
related to coherent logic. For example, we may consider \emph{classical logic},
which introduces the negation connective \(\neg\).\footnote{%
	The universal quantifier \(\forall\) and the conditional \(\Longrightarrow\) are also
	standard connectives in classical logic, but these can be defined in terms of
	coherent connectives and negation.
}%
\[%
	\begin{bprooftree}
		\AxiomC{}
		\UnaryInfC{$\fCenter \phi \vee \neg\phi$}
	\end{bprooftree} \qquad%
	\begin{bprooftree}
		\AxiomC{}
		\UnaryInfC{$\phi \wedge \neg\phi \fCenter \bot$}
	\end{bprooftree} \qquad%
	\begin{bprooftree}
		\AxiomC{}
		\UnaryInfC{$\neg\neg \phi \fCenter \phi$}
        \end{bprooftree} \qquad%
        \begin{bprooftree}
        \AxiomC{}
        \UnaryInfC{$\phi \fCenter \neg\neg\phi$}
	\end{bprooftree}
\]%
The negation connective gives the syntactic category of a classical theory
additional structure: each subobject has a complement, making the syntactic
category a \emph{Boolean (coherent) category}. In order to recover a
biequivalence, we introduce a new reconstrual rule and an axiom schema for the
internal logic.
\begin{mathrule}[Negation preservation]\label{rule:negation-preservation}
	Given a reconstrual \(F: \Sigma_1 \to \Sigma_2\) between classical signatures
	and a \(\Sigma_1\)-formula \(\phi\), define%
	\[%
		F^+\pwrap{\neg\phi(\vec{x})} \equivdef \fb{ \neg F^+\phi(\vec{t}) }.
	\]%
\end{mathrule}

\begin{enumerate}
\item[IL12] \emph{Axiom for complements}: For a pair of monomorphisms \(f: A
\hookrightarrow X\) and \(g: B \hookrightarrow X\) such that \(g\)
presents the complement of \(f\), the sequents%
\[%
	\exists b^{\ul B} \ul g(b) = x \dashv\vdash \neg \pwrap{ \exists a^{\ul
	A} \ul f(a) = x }.
\]%
\end{enumerate}

Axiom schema \axiom{IL12} ensures that the internal logic of the syntactic
category of a classical theory \(T\) is still e.p.\ bi\-/interpretable with \(T\).
Strictly speaking, Rule~\ref{rule:negation-preservation} is not necessary, since
for any coherent translation \(F: T_1 \to T_2\) between classical theories,
\(F(\neg\phi)\) is logically equivalent to \(\neg F\phi\) via the law of
excluded middle and Rule~\ref{rule:disjunction-preservation}---hence any
coherent translation between classical theories is homotopic in a natural way to
a classical translation. Nevertheless, introducing
Rule~\ref{rule:negation-preservation} is harmless, and the proofs for
Theorems~\ref{theorem:cth0}, \ref{theorem:CThEqCoh-Correspondence},
and~\ref{theorem:bi-interpretability-equivalence} prove the analogous theorems for
classical logic.
\begin{theorem}\label{theorem:classical-logic}
	Let \(\Th_0\) be the collection of (small) classical theories. \(\Th_0\) is a
	bicategory, whose 1\-/cells are `classical translations' (coherent
	translations which also satisfy Rule~\ref{rule:negation-preservation}) and
	2\-/cells are t\-/maps. The full sub-bicategory \(\ThEq\) spanned by e.p.\
	classical translations is part of a biequivalence: \(\C^\Bool: \ThEq \to
	\Bool\) and \(\T^\Bool: \Bool \to \ThEq\), where \(\C^\Bool\) is the
	syntactic category pseudofunctor as in
	Proposition~\ref{prop:C-pseudofunctor}, and \(\T^\Bool\) is the internal
	logic pseudofunctor as in Proposition~\ref{prop:T-pseudofunctor}, where we
	also include axiom schema \axiom{IL12}.

	Furthermore, the functor \(\X: \h\CTh_0 \to \h\Exact\) from Section 5
	restricts to an equivalence \(\h\Th_0 \to \h\mathsf{ExactBool}\), where
	\(\mathsf{ExactBool}\) is the bicategory of Boolean Barr-exact categories. In
	particular, two classical theories \(T_1\) and \(T_2\) are bi\-/interpretable
	if and only if \(\C^\Bool(T_1)^\ex\) and \(\C^\Bool(T_2)^\ex\) are equivalent
	categories.
\end{theorem}

This procedure of introducing additional reconstrual rules and internal logic
axioms for more expressive fragments of predicate logic yield similar results
as Theorem~\ref{theorem:classical-logic}. For example, first\-/order
intuitionistic theories with e.p.\ translations are biequivalent to Heyting
categories with Heyting functors. In this case we must introduce reconstrual
rules for preserving the new connectives \(\forall, \Longrightarrow\), and
\(\neg\) (see \axiom{IL11} of~\ref{subsection:IL-axiom-schemata}). On the other
hand, if we consider \(\kappa\)-coherent logic \(L^g_{\kappa\omega}\), then we
obtain a biequivalence with \(\kappa\)-coherent categories (called
\emph{\(\kappa\)-logical} in~\cite{MakkaiReyes1977}) as long as we modify
Rule~\ref{rule:disjunction-preservation} to account for \(\kappa\)-ary
disjunctions. The proofs of these results are essentially the same as in the
case \(L^g_{\omega\omega}\) which we have covered.

\subsection{Morleyization}

Since Boolean categories \(\Bool\) form a (full) sub\-/bicategory of
coherent categories \(\Coh\), the inclusion pseudofunctor \(\iota: \Bool \to
\Coh\) can be used to embed classical logic, \(\ThEq\), into coherent logic,
\(\CThEq\). This is the idea behind \emph{Morleyization}. One account of
Morleyization is Lemma D1.5.13 of~\cite{johnstone2002sketches}. We consider only
the case of a classical theory.
\begin{lemma*}
	Given a classical theory \(T\) with signature \(\Sigma\), there exists a
	signature \(\Sigma'\) containing \(\Sigma\) and a coherent theory \(T'\) with
	signature \(\Sigma'\) such that for any Boolean category \(C\) there
	is an equivalence of categories%
	\[%
		\Mod(T,C) \approx \Mod_{\text{elem}}(T',C),
	\]%
	where \(\Mod_{\text{elem}}(T',C)\) is the subcategory of models of \(T'\) in
	\(C\) with morphisms elementary embeddings.
\end{lemma*}
Using the correspondence between models and coherent functors established
in~\cite{MakkaiReyes1977}, this equivalence of categories can be described using
syntactic categories:%
\[%
	\Bool\pwrap{\C^\Bool(T),C} \approx \Coh\pwrap{\C(T'),\iota(C)}.
\]%
Moreover, the proof of D1.5.13 implies the stronger result that this equivalence
is natural in \(C\). On the other hand, \(\C(T')\) is a Boolean
category, since the signature \(\Sigma'\) adds a
complement to every \(\Sigma\)\-/formula. The 2\-/categorical Yoneda lemma thus
implies that \(\C^\Bool(T)\) and \(\C(T')\) are equivalent categories. Using the
biequivalences \(\C^\Bool\) and \(\T\), this yields the following
characterization of \(T'\).

\begin{theorem}\label{theorem:morleyization-is-unit}
	Let \(T\) be a classical theory. Its Morleyization \(T'\) is e.p.\
	bi\-/interpretable with the coherent theory \(\T\iota\C^\Bool(T)\).
\end{theorem}

This allows us to extend the Morleyization operation into a pseudofunctor, up to
e.p.\ bi\-/interpretability.

\begin{corollary*}
	Morleyization extends to a pseudofunctor \(\M: \ThEq \to \CThEq\) given by
	\(\M \eqdef \T \iota \C^\Bool\).
\end{corollary*}

Since \(\T\) and \(\C^\Bool\) are biequivalences, and \(\iota\) an inclusion,
we obtain a finer description of Morleyization.

\begin{corollary*}
	Two classical theories \(T_1, T_2\) are e.p.\ bi\-/interpretable if and only
	if their Morleyizations \(\M(T_1)\) and \(\M(T_2)\) are e.p.\
	bi\-/interpretable.
\end{corollary*}

The closest related result in the literature is found
in~\cite{tsementzis2015syntactic}, where it is shown that \(T_1\) and \(T_2\)
are Morita equivalent if and only if \(\M(T_1)\) and \(\M(T_2)\) are Morita
equivalent (Morita equivalence of classical theories may be understood as
in~\cite{barrett_halvorson_2016}).

Using the Morleyization pseudofunctor \(\M\), the relationship between \(T\) and
\(\M(T)\) resembles an adjunction. Let \(\CThEq_\sim\) and \(\ThEq_\sim\) denote
the \(\pwrap{2,1}\)-categories obtained by remembering only the invertible
\(2\)-cells in the bicategories \(\CThEq\) and \(\ThEq\) respectively. There is
a pseudofunctor \(\Cl: \CThEq_\sim \to \ThEq_\sim\), where \(\Cl(T)\) is the
theory obtained by considering the coherent theory \(T\) as a classical theory.
Given a translation \(F: T_1 \to T_2\), the same underlying reconstrual defines
a translation \(\Cl(F): \Cl(T_1) \to \Cl(T_2)\), and any invertible t\-/map
\(\chi: F \Rightarrow G\) defines an invertible t\-/map \(\Cl(\chi): \Cl(F)\
\Rightarrow \Cl(G)\).

\begin{proposition}[Morleyization adjunction]
	For any coherent theory \(T_1\) and any classical theory \(T_2\), there is an
	equivalence of categories
	\[%
		\ThEq_\sim\pwrap{\Cl(T_1), T_2} \approx \CThEq_\sim\pwrap{T_1, \M(T_2)}.
	\]%
\end{proposition}
\begin{proof}[Proof Sketch]
	\(\M(T_2)\) is constructed by introducing predicates \(C_\phi\) and \(D_\phi\)
	for every \(T_2\)-formula \(\phi\) so that \(C_\phi\) and \(D_\phi\) behave
	like coherent analogues of \(\phi\) and \(\neg\phi\) respectively. With this
	in mind, given a translation \(F: \Cl(T_1) \to T_2\), we obtain the
	reconstrual \(F^\M: T_1 \to \M(T_2)\) by setting \(F^\M R\) (for a
	\(T_1\)-relation \(R\)) to be the substitution class obtained by replacing any
	copy of \(\neg S(x)\) with \(D_S(x)\) for every atomic \(T_2\)-formula \(S\)
	appearing in \(FR\). On the other hand, given a translation \(G: T_1 \to
	\M(T_2)\), we obtain the associated \(G^{\Cl}: \Cl(T_1) \to T_2\) by replacing
	any copy of \(C_S(x)\) with \(S(x)\) and any copy of \(D_S(x)\) with \(\neg
	S(x)\) appearing in the expansion of \(G R\), for \(R\) a \(T_1\)-relation.
\end{proof}

If the above equivalence of categories satisfies the appropriate coherence
conditions, then this shows that Morleyization fits into a bi-adjunction \(\Cl
\dashv \M\). The existence of such a bi-adjunction would show that Morleyization
is essentially unique---in the sense that it is the best approximation to an
inverse to the operation \(T \mapsto \Cl(T)\) which ``forgets'' that \(T\) is
restricted to a smaller fragment of logic. There is evidence that this is the
case. Proposition 5.6 of~\cite{tsementzis2015syntactic} proves that, for a
classical theory \(T\), \(\Cl\M(T)\) is Morita equivalent to \(T\). In fact,
this Morita equivalence does not involve coproduct Morita extensions. Pairing
this with Proposition 5.12 of~\cite{mceldowney_2020}, this implies that \(T\)
and \(\Cl\M(T)\) are bi\-/interpretable.

On the other hand, the adjunction is a framework for generalizing Morleyization
to other fragments more expressive than coherent logic. By identifying these
more-expressive fragments with sub-bicategories of \(\Coh\), we propose
investigating the following generalized Morleyization adjunction.

\begin{question}
	Let \(\mathsf{D}\) be a sub-bicategory of \(\Coh\) closed under
	equivalence of categories, and let \(d: \mathsf{D}_\sim \to \Coh_\sim\) be the
	restriction of the inclusion pseudofunctor to invertible 2\-/cells. Under what
	conditions on \(\mathsf{D}_\sim\) does \(d\) admit a (weak) left\-/adjoint?
\end{question}

The pseudofunctor \(d\) plays the role of Morleyization from
\(\mathsf{D}\)-theories to coherent theories. The left adjoint plays the
``forgetful'' role of interpreting a coherent theory \(T\) as a
\(\mathsf{D}\)-theory. Examples of this generalized Morleyization adjunction
include \(\mathsf D \eqdef \Pretopos\), where the left-adjoint is pretopos
completion, and \(\mathsf D \eqdef \Exact\) with left-adjoint the exact
completion.

\subsection{Morita Equivalence}

The last decade has seen papers~\cite{AraH:2018,barrett_halvorson_2016,halvorson-logic,
harnik2011,mceldowney_2020,tsementzis2015syntactic, washington2018} 
compare Morita equivalence and bi\-/interpretability to various
operations on syntactic categories. Between bi\-/interpretability and Morita
equivalence, the state\-/of\-/the\-/art is found in~\cite{mceldowney_2020}
and~\cite{halvorson-logic,washington2018}.

\begin{theorem*}[\cite{mceldowney_2020}, Proposition 5.12]
	Suppose two theories \(T_1\) and \(T_2\) are Morita equivalent via sequences
	of Morita extensions which do not define any coproduct sorts. Then \(T_1\) and
	\(T_2\) are bi\-/interpretable.
\end{theorem*}

\begin{theorem*}%
	[\cite{washington2018}, Theorem 7;~\cite{halvorson-logic}, Theorem 7.5.5]
	Suppose \(T_1\) and \(T_2\) are bi\-/interpretable theories. Then \(T_1\) and
	\(T_2\) are Morita equivalent.
\end{theorem*}

These results were proven for classical theories, but the proofs also work for
coherent theories. It is important to permit Morita extensions which define
\(n\)-ary product sorts or coproduct sorts at
once---otherwise~\cite{mceldowney_2020} provides a counterexample following the
statement of Proposition 5.12.

From the category theory perspective, these results were anticipated
in~\cite{tsementzis2015syntactic} and~\cite{harnik2011}.

\begin{theorem*}[\cite{tsementzis2015syntactic}, Theorem 3.9]
	Two coherent theories \(T_1\) and \(T_2\) are Morita equivalent if and only if
	the pretopos completions of \(\C(T_1)\) and \(\C(T_2)\) are equivalent
	categories.
\end{theorem*}

When restricting to the class of \emph{proper} theories (those for which any
model has at least two distinct elements),~\cite{harnik2011} shows that the
pretopos completion of \(\C(T)\) is the syntactic category \(\C(T^\eq)\), where
\(T \to T^\eq\) is Shelah's elimination of imaginaries construction. \(T^\eq\)
defines a new sort for every quotient definable in \(T\). The restriction to
proper theories ensures that \(\C(T^\eq)\) is a pretopos by ensuring that
coproducts can be defined in terms of a quotient (see~\cite[Theorem 5.3]{harnik2011}).  
For a general coherent (resp.\ classical) theory \(T\),
\(\C(T^\eq)\) is just a Barr\-/exact coherent (resp.\ Boolean) category.
However, \(\C(T^\eq)\) is still the \emph{exact} completion \(\C(T)^\ex\). In
the language of~\cite{harnik2011}, the inclusion \(\C(T) \to \C(T^\eq)\) is a
tight extension and \(\C(T^\eq)\) is closed under taking quotients. A theory
\(T\) being proper is the only obstruction between the exact completion and
pretopos completion of \(\C(T)\).

\begin{proposition}\label{prop:properness-characterization}
	A consistent coherent (or classical) theory \(T\) is proper if and only if
	\(\C(T)^\ex\) is a pretopos.
\end{proposition}
\begin{proof}
	The forward direction is Theorem 5.3 of~\cite{harnik2011}. Suppose
	\(\C(T)^\ex\) is a pretopos. Then the coproduct \(\2 = \1 \amalg \1\) is
	definable in \(\C(T)^\ex\). This means that \(\2\) is a quotient of a
	congruence \(R \rightrightarrows \fb\phi\) in \(\C(T)\). In particular, there
	exists a pair of subobjects \(\fb{\psi_1}\) and \(\fb{\psi_2}\) of \(\fb\phi\)
	which project to each copy of \(\1\) in the quotient. Since these copies are
	disjoint, \(\psi_1\) and \(\psi_2\) must satisfy \(\psi_1(x) \wedge \psi_2(x)
	\vdash \bot\). On the other hand, since \(T\) is consistent, \(\1\) is never
	empty, so \(\psi_1\) and \(\psi_2\) are never empty. Thus \(T\) is proper,
	witnessed by \(\psi_1\) and \(\psi_2\).
\end{proof}

By understanding the interplay between proper theories and elimination of
imaginaries, our Theorem~\ref{theorem:BiInterpretability-result} gives us a more
complete picture of Morita equivalence and bi\-/interpretability. Namely,
\(T_1\) and \(T_2\) are bi\-/interpretable if and only if \(\C(T_1)^\ex\) and
\(\C(T_2)^\ex\) are equivalent categories. When \(T_1\) and \(T_2\) are proper,
this is equivalent to the condition that the pretopos completions are
equivalent---which is the same as Morita equivalence.

\begin{corollary}
	Two theories \(T_1\) and \(T_2\) are bi\-/interpretable if and only if
	\(\C(T_1^\eq)\) and \(\C(T_2^\eq)\) are equivalent categories.
\end{corollary}

We conclude with the observation that Harnik essentially conjectured the
preceding corollary in~\cite{harnik2011}. Following Definition 6.2
of~\cite{harnik2011}, Harnik proposes that, for theories with finite signatures,
the most general, reasonable, notion of interpretability \(T_1 \to T_2\) is a
coherent functor \(\C(T_1) \to \C(T_2^\eq)\). Moreover, Harnik uses this
proposal to conclude that \(T_1\) and \(T_2\) are bi\-/interpretable if and only
if \(\C(T_1)^\eq\) and \(\C(T_2)^\eq\) are equivalent. Using
Theorem~\ref{theorem:bi-interpretability-equivalence}, we see that Harnik's
proposed general notion of interpretation is our notion of a translation \(T_1
\to T_2\). From this we recover Harnik's conjecture, and we recover it in a
framework which permits theories with infinite signatures that are not
necessarily proper.
\subsection*{Acknowledgements}

\addcontentsline{toc}{section}{Acknowledgements}

The authors would like to thank Hans Halvorson for his invaluable guidance
and knowledge, the Princeton University Mathematics and Philosophy
Departments for their generous financial support, and Krist\'of Kanalas for
helpful remarks on the 0\-/cells of a syntactic category.

\appendix
\section{IL Axioms and Additional Proofs}

\subsection{IL Axiom Schemata}
\label{subsection:IL-axiom-schemata}

Let \(C\) be a small coherent category, and let \(\ul{\Sigma}_C\) be the internal
signature defined in Definition~\ref{def:internal-logic-0-cell}. Makkai and
Reyes defined eleven axiom schemata (see \S 2.4 of~\cite{MakkaiReyes1977})
encoding the possible axioms for the theory \(\T(C)\) based on the (co)limits of \(C\). We list them
below.
\begin{enumerate}
	\item[IL1] \emph{Axiom for identity}: For every \(1_A: A \to A\), the sequent \(\vdash \ul{1_A}(x) =_{\ul{A}} x\).
	\item[IL2] \emph{Axiom for commutative diagrams}: For every triplet of
		morphisms \(f,g,h\) such that \(fg = h\), the sequent \(\vdash
		\ul{h}(x) = \ul{f}(\ul{g}(x))\).
	\item[IL3] \emph{Axiom for a monomorphism}: For every monic \(f: A \to B\),
		the sequent \(\ul{f}(x_1) =_{\ul{B}} \ul{f}(x_2) \vdash x_1 =_{\ul{A}} x_2\).
	\item[IL4] \emph{Axioms for a terminal object}: For a terminal object \(A\),
		the sequents \(\vdash x =_{\ul A} y\) and \(\vdash \exists x^{\ul A} x =
		x\).
	\item[IL5] \emph{Axioms for an equalizer}: For an equalizer \(E \xrightarrow{e} A \doublerightarrow{f}{g} B\), the sequents%
		\[%
			\begin{array}{rl}%
				\ul{e}(x) =_{\ul{A}} \ul{e}(y) &\vdash x =_{\ul{E}} y,\\
				&\vdash \ul{f}(\ul{e}(x)) =_{\ul{B}} \ul{g}(\ul{e}(x)),\\
				\ul{f}(x) =_{\ul{B}} \ul{g}(x) &\vdash \exists z^{\ul E} \ul{e}(z) =_{\ul{A}} x.
			\end{array}%
		\]%
	\item[IL6] \emph{Axioms for a product}: For a product \(A \xleftarrow{f} X
		\xrightarrow{g} B\), the sequents%
		\[%
			\begin{array}{rl}%
				\ul{f}(x) =_{\ul{A}} \ul{f}(y) \wedge \ul{g}(x) =_{\ul{B}}
				\ul{g}(y) &\vdash x =_{\ul{X}} y,\\
				&\vdash \exists x^{\ul X} \left(\ul{f}(x) =_{\ul{A}} a \wedge
				\ul{g}(x) =_{\ul{B}} b\right).
			\end{array}%
		\]%
	\item[IL7] \emph{Axiom for an initial object}: For an initial object \(\ul A\),
		the sequent \(x =_{\ul A} x \vdash \bot\).
	\item[IL8] \emph{Axioms for unions}: For a family of subobjects \(f_i: A_i
		\hookrightarrow X\) and \(g: B \hookrightarrow X\) such that \(B \cong
		\bigvee_i A_i\), the sequents%
		\[%
			\bigvee_i \exists a^{\ul{A_i}} \ul{f_i}(a) =_{\ul X} x
			\dashv\vdash \exists b^{\ul B} \ul{g}(b) =_{\ul X} x.
		\]%
	\item[IL9] \emph{Axiom for images}: For a regular epimorphism\footnote{%
		Called a \emph{surjective morphism} in~\cite{MakkaiReyes1977}.
		} \(f: A \to B\), the sequent \( \vdash \exists a^{\ul A} \ul{f}(a)
		= b\).
	\item[IL10] \emph{Axioms for intersections}: For a family of subobjects \(
	f_i: A_i \hookrightarrow X\) and \(g: B \hookrightarrow X\)
	such that \(B \cong \bigwedge A_i\), the sequents%
		\[%
			\bigwedge_i \exists a^{\ul{A_i}} \ul{f_i}(a) =_{\ul X} x
			\dashv\vdash \exists b^{\ul B} \ul{g}(b) =_{\ul X} x.
		\]%
	\item[IL11] \emph{Axioms for dual images}: For a morphism \(f: X \to Y\), a
		family of subobjects \(f_i: A_i \hookrightarrow X\), and \(g: B
		\hookrightarrow Y\) such that \(B \cong \forall_{f} (A_1 \Longrightarrow
		A_2)\), the sequents%
		\[%
			\exists b^{\ul B} \ul{g}(b) =_{\ul{Y}} y \dashv\vdash \forall x^{\ul{X}}
			\left(\left( \ul{f}(x) =_{\ul{Y}} y \wedge \exists a_1^{\ul{A_1}}
			\ul{f_1}(a_1) = x \right)\Longrightarrow \exists a_2^{\ul{A_2}}
			\ul{f_2}(a_1) = x \right).
		\]%
\end{enumerate}

\subsection{Limits and Colimits in the Syntactic Category as Sequents}
\label{subsection:sc-schemata}
\begin{proposition}\label{prop:sc-axiom-schemata}
	Let \(\C(T)\) be the syntactic category of a small coherent theory \(T\). The following relates 
	(co)limits of \(\C(T)\) to
	sequent provability in \(T\).%
\end{proposition}
{%
	\renewcommand\labelenumi{(SC\theenumi)}%
	\begin{enumerate}%
        \item \emph{(Identity)} \(1_{\fb{\phi}} : \fb{\phi} \to \fb{\phi}\) is an identity
                if and only if \(\phi(x) \vdash 1_{\fb{\phi}}(x,x)\).
        \item \emph{(Commutative diagrams)} \(h : \fb{\phi} \to \fb{\eta}\) is 
        \(\fb{\phi} \xrightarrow{f} \fb{\psi} \xrightarrow{g} \fb{\eta}\) if and only if
        \[\phi(x) \vdash \exists z \left(h(x,z) \wedge \exists y\left(f(x,y) \wedge g(y,z)\right)\right).\]
		\item \emph{(Monomorphism)} \(f : \fb{\phi} \to \fb{\psi}\) is a
			monomorphism if and only if \(f(x_1,y) \wedge f(x_2,y) \vdash x_1 =
			x_2\).
		
		\item \emph{(Terminal object)} \(\fb{\phi}\) is a terminal object if and
			only if \(\phi(x) \wedge \phi(y) \vdash x = y\) and \(\vdash \exists x \,
			\phi(x)\).
		
		\item \emph{(Equalizer)} A diagram \(\fb{\eta}
			\xrightarrow{e} \fb{\phi} \doublerightarrow{f}{g} 
                \fb{\psi}\) is an equalizer if and only if%
			\begin{align*}%
				e(x_1,y) \wedge e(x_2,y) &\vdash x_1 = x_2,\\
				\eta(x) &\vdash \exists z\left(\exists y_1 \left(e(x,y_1) \wedge
				f(y_1,z)\right) \wedge \exists y_2\left(e(x,y_2) \wedge
				g(y_2,z)\right)\right),\\
				f(y_1,z) \wedge g(y_1,z) &\vdash \exists x \, e(x,y).
			\end{align*}%
		
		\item \emph{(Product)} A diagram \(\fb{\phi} \xleftarrow{f} \fb{\eta}
			\xrightarrow{g} \fb{\psi}\) is a product if and only if%
			\begin{align*}%
					f(x_1,y) \wedge f(x_1,y) \wedge g(x_1,z) \wedge g(x_2,z) &\vdash x_1
					= x_2,\\
					\phi(y) \wedge \psi(z) &\vdash \exists x \left(f(x,y) \wedge
					g(x,z)\right).
			\end{align*}
		
		\item \emph{(Initial object)} \(\fb{\phi}\) is an initial object if and
			only if \(\phi(x) \vdash \bot\).
		
		\item \emph{(Unions)} \(g : \fb{\psi} \hookrightarrow \fb{\eta}\) is the
			union of \(f_i : \fb{\phi_i} \hookrightarrow \fb{\eta}\) if and only if%
			\[%
				\exists b \, g(b,x) \dashv\vdash \bigvee_i \exists a_i f_i(a_i,x).
			\]%
		
		\item \emph{(Images)} \(f : \fb{\phi} \to \fb{\psi}\) is a regular
			epimorphism if and only if \(\psi(y) \vdash \exists x f(x,y)\).
		
		\item \emph{(Intersections)} \(g : \fb{\psi} \hookrightarrow \fb{\eta}\) is
			the intersection of \(f_i : \fb{\phi_i} \hookrightarrow \fb{\eta}\) if
			and only if%
			\[%
				\exists b \, g(b,x) \dashv\vdash \bigwedge_i \exists a_i f_i(a_i,x).
			\]%
		
		\item \emph{(Dual images)} \(g : \fb{\psi} \hookrightarrow \fb{\eta}\) is
			\(\forall_f\left(\fb{\phi_1} \Longrightarrow \fb{\phi_2}\right)\) for \(f
			: \fb{\phi} \to \fb{\eta}\), \(f_i : \fb{\phi_i} \hookrightarrow
			\fb{\phi}\) if and only if%
			\[%
				\exists b \, g(b,y) \dashv\vdash \forall x\left(\left(f(x,y) \wedge
				\exists a_1 f_1(a_1,x)\right) \Longrightarrow \exists a_2
			f_2(a_2,x)\right).
			\]%
	\end{enumerate}
}%
\subsection{Reconstrual Properties}
\label{subsection:reconstrual-properties}

\begin{proof}%
        [Proof of Proposition~\ref{prop:e.p.-domain-formula}]
		Since \(F\) is a translation and the sequent \(s'=_\sigma t' \vdash t'=_\sigma t'\) holds by \(=\)\-/introduction, 
		we have \(E^F_\sigma(s,t) \equiv F\fb{s'=_\sigma t'}(s,t) \vdash F\fb{t'=_\sigma t'}(t,t) \equiv D^F_\sigma(t)\). 
		Since \(F\) is e.p., the forward direction is due to \(\wedge\)\-/introduction. For the converse direction, we 
		note that \(E^F_\sigma(t,t) \equiv D^F_\sigma(t)\) is obtained by Rule~\ref{rule:context-duplication}, allowing us 
		to conclude with \(=\)-elimination and the cut rule.
\end{proof}

\begin{proof}%
	[Proof of Proposition~\ref{prop:reconstrual-strictly-composes}]
		We induct on logical connectives. We need to show that the substitution
		class \(\pwrap{GF}^+\fb\phi\) is identical to the substitution class \(G^+
		F^+ \fb\phi\). The base case of the induction proof is the case where either
		\(\fb\phi \equiv \top\) or \(\bot\) or where \(\fb\phi\) is the substitution
		class of a relation \(R \in \Sigma_1\), i.e., \(\fb\phi \equiv
		\fb{R(\vec{x})}\). For \(\top\) and \(\bot\) the proposition is a consequence
		of Rule~\ref{rule:top-bot-preservation}. For the second case, we apply
		Rule~\ref{rule:relations}.%
		\[%
			G^+\pwrap{F^+\fb\phi} \equiv G^+\pwrap{F R} \equiv \pwrap{GF}(R) \equiv
			\pwrap{GF}^+\fb{R(\vec x)} \equiv \pwrap{GF}^+\fb\phi.
		\]%
		For the inductive step, we need to show that, given \(\Sigma_1\)-formulae
		\(\phi_1,\ldots,\phi_n\) satisfying the inductive hypothesis \(G^+\pwrap{F^+
		\fb{\phi_i}} \equiv \pwrap{GF}^+\fb{\phi_i}\), then any combination of
		\(\phi_1,\ldots,\phi_n\) using Rules~\ref{rule:conjunction-preservation}\-/\ref{rule:term-reduction} also
		satisfies the induction hypothesis.		

		We begin with Rule~\ref{rule:conjunction-preservation}. Given
		\(\Sigma_1\)-formulae \(\phi_1\) and \(\phi_2\) satisfying the inductive
		hypothesis, we have the following chain of identities:%
		\begin{gather*}%
			G^+\pwrap{F^+\fb{\phi_1(x_1) \wedge \phi_2(x_2)}}
			\equiv G^+ \fb{F^+\phi_1(x_1') \wedge F^+\phi_2(x_2') }\\
			\equiv \fb{G^+ F^+ \phi_1(x_1'') \wedge G^+ F^+ \phi_2(x_2'')}
			\equiv \fb{\pwrap{GF}^+ \phi_1(x_1'') \wedge \pwrap{GF}^+ \phi_2(x_2'')}\\
			\equiv \pwrap{GF}^+\fb{\phi_1(x_1) \wedge \phi_2(x_2)},
		\end{gather*}%
		where the third identity uses the inductive hypothesis for \(\phi_1\) and
		\(\phi_2\), and the fourth identity is
		Rule~\ref{rule:conjunction-preservation} for the reconstrual \(GF\).

		All the other reconstrual rules use the same argument, where the third
		identity uses the inductive step, and the fourth identity will always use
		the pertinent reconstrual rule. For example, for
		Rule~\ref{rule:term-reduction}, we consider four identities.
		\begin{gather*}
			G^+ F^+ \fb{\phi(x_1,f(y),x_2)}
			\equiv G^+ \fb{\exists t' F f(y',t') \wedge F^+\phi(x_1',t',x_2')}\\
			\equiv \fb{\exists t'' G^+ F f(y'',t'') \wedge G^+ F^+\phi(x_1'',t'',x_2'')}\\
			\equiv \fb{\exists t'' (GF)f(y'',t'') \wedge \pwrap{GF}^+\phi(x_1'',t'',x_2'')}
			\equiv \pwrap{GF}^+ \fb{\phi(x_1,f(y),x_2)}.
		\end{gather*}
		Since any formula is built inductively from logical connectives, this proves
		the inductive statement holds for any \(\Sigma_1\)-substitution class, as
		desired.
\end{proof}

\begin{proof}%
	[Proof of Proposition~\ref{prop:reconstrual-log-eq-to-identity}]
	We induct on logical connectives. For the case \(\fb\phi
	\equiv \fb{R(\vec x)}\), observe \(1^+_\Sigma \fb\phi \equiv 1_\Sigma R \equiv
	\fb R\), so \(1_\Sigma^+ \fb R \dashv\vdash \fb R\). In the case where \(\phi\) contains
	no function symbols, the reconstrual rules guarantee \(1_\Sigma^+ \fb\phi
	\equiv \fb\phi\); hence \(1_\Sigma^+ \fb\phi \dashv\vdash \fb\phi\). This
	leaves the case where \(\phi\) has function symbols, so we need to use
	Rule~\ref{rule:term-reduction}. Let \(\phi'\) be the \(\Sigma\)-formula such
	that \(\phi'(x_1,f(y),x_2) \equiv \phi(x_1,y,x_2)\). Assume \(\phi'\)
	satisfies the inductive hypothesis: \(1_\Sigma^+ \fb{\phi'(x_1,t,x_2)}
	\dashv\vdash \fb{\phi'(x_1,t,x_2)}\). Observe:
	\begin{gather*}
		1_\Sigma^+ \fb{\phi(x_1,y,x_2)} 
		\equiv 1_\Sigma^+ \fb{\phi'(x_1,f(y),x_2)}
		\equiv \fb{\exists t 1_\Sigma f(y,t) \wedge 1_\Sigma^+ \phi'(x_1,t,x_2)}\\
		\equiv \fb{\exists t f(y) = t \wedge 1_\Sigma^+\phi'(x_1,t,x_2)}
		\equiv \fb{\exists t f(y) = t \wedge \phi'(x_1,t,x_2)}\\
		\dashv\vdash \fb{\phi'(x_1,f(y),x_2)} \equiv \fb{\phi(x_1,y,x_2)},
	\end{gather*}
	where the fourth identity uses the inductive hypothesis for \(\phi'\). Like in
	the proof of Proposition~\ref{prop:reconstrual-strictly-composes}, this shows that the inductive hypothesis is
	satisfied by any \(\Sigma\)-substitution class, as desired.
\end{proof}

\subsection{Bicategory Proofs}\label{subsection:bicategory-proofs}

\begin{proof}%
	[Proof of Lemma~\ref{lemma:t-map-domain-lemma}]
	The key point is that \(x = x\) is tautological.
	\begin{prooftree}
			\Axiom$D_\sigma^{F_1}(x) \fCenter x =_{F_1 \sigma} x$
			\RightLabel{\(F_2\)}
			\UnaryInf$F_2 D_\sigma^{F_1}(x') \fCenter D_{F_1\sigma}^{F_2}(x')$
			\UnaryInf$D_\sigma^{F_2 F_1}(x') \fCenter D_{F_1 \sigma}^{F_2}(x')$
	\end{prooftree}
\end{proof}

\begin{proof}%
	[Proof of Lemma~\ref{lemma:t-map-equality-lemma}]
	This uses the fact that \(x = x\) is a tautology along with
	\(=\)-introduction.
	\begin{prooftree}
			\Axiom$ \fCenter x =_\sigma x$
			\RightLabel{\(F_1\)}
			\UnaryInf$ D_\sigma^{F_1}(x') \fCenter E_\sigma^{F_1}(x',x')$
			\RightLabel{\(=\)-intro.}
			\UnaryInf$ x' =_{F_1 \sigma} y' \wedge \left( D_\sigma^{F_1}(x') \vee
			D_\sigma^{F_1}(y') \right) \fCenter E_\sigma^{F_1}(x',y')$
			\RightLabel{\(F_2\)}
			\UnaryInf$E^{F_2}_{F_1\sigma}(x'',y'') \wedge \left( D_\sigma^{F_2
			F_1}(x'') \vee D_\sigma^{F_2 F_1}(y'') \right) \fCenter E_\sigma^{F_2
			F_1}(x'',y'')$
	\end{prooftree}
\end{proof}

\begin{proof}%
	[Proof of Lemma~\ref{lemma:t-map-Z-lemma}]
	We begin with the forward sequent \(\pwrap{\eta\circ\chi}_\sigma(s,t) \vdash
	Z_\sigma(s,t)\). Using \(\exists\)\-/introduction on \(w\) and \(t'\), it
	suffices to prove%
	\begin{align*}%
		F_2\chi_\sigma(s,w) \wedge &\eta_{G_1\sigma}(w,t') \wedge E^{G_2
		G_1}_\sigma(t',t) \\%
		&\vdash \exists {s'}^{F_2 F_1\sigma}E^{F_2 F_1}_\sigma(s,s') \wedge \exists
		v^{G_2 F_1 \sigma} \pwrap{%
			\eta_{F_1\sigma}(s',v) \wedge G_2\chi_\sigma(v,t)%
		}.
	\end{align*}%
	Applying \(F_2\) to \axiom[\chi]{TM1} implies \(F_2\chi_\sigma(s,w) \vdash
	D^{F_2 F_1}_\sigma(s)\). Lemma~\ref{lemma:t-map-domain-lemma} and the cut
	rule implies \(F_2\chi_\sigma(s,w) \vdash D^{F_2}_{F_1\sigma}(s)\).
	Therefore, \axiom[\eta]{TM3} and the cut rule arrives at the sequent \(F_2\chi_\sigma(s,w)
	\vdash \exists v^{G_2 F_1\sigma} \eta_{F_1\sigma}(s,v)\). Moreover,
	\axiom[\eta]{TM5} implies%
	\[%
		\eta_{F_1\sigma}(s,v) \wedge \eta_{G_1\sigma}(w,t') \wedge
		F_2\chi_\sigma(s,w) \vdash G_2\chi_\sigma(v,t').
	\]%
	We can apply \(G_2\) to \axiom[\chi]{TM2} to deduce \(E^{G_2
	G_1}_\sigma(t',t) \wedge G_2\chi_\sigma(v,t') \vdash G_2\chi_\sigma(v,t)\).
	Gathering together the above sequents, we deduce:%
	\begin{align*}
		F_2\chi_\sigma(s,w) \wedge \eta_{G_1\sigma}(w,t') \wedge
		&E^{G_2G_1}_\sigma(t',t)\\
		&\vdash \exists v^{G_2 F_1\sigma} \pwrap{%
			\eta_{F_1\sigma}(s,v) \wedge G_2\chi_\sigma(v,t)
		}\\%
		&\vdash E^{F_2 F_1}(s,s) \wedge \exists v \pwrap{%
			\eta_{F_1\sigma}(s,v) \wedge G_2\chi_\sigma(v,t)
		}\\%
		&\vdash \exists {s'}^{F_2 F_1\sigma} E^{F_2F_1}_\sigma(s,s') \wedge \exists
		v \pwrap{%
			\eta_{F_1\sigma}(s',v) \wedge G_2\chi_\sigma(v,t)
		}.%
	\end{align*}%
	We now prove the converse sequent. Like before, it suffices to show%
	\begin{align*}%
		E^{F_2 F_1}_\sigma(s,s') \wedge &\eta_{F_1\sigma}(s',v) \wedge
		G_2\chi_\sigma(v,t)\\%
		&\vdash \exists {t'}^{G_2 G_1\sigma} \exists w^{F_2 G_1\sigma}
		F_2\chi_\sigma(s,w) \wedge \eta_{G_1\sigma}(w,t') \wedge E^{G_2
		G_1}_\sigma(t',t).
	\end{align*}%
	Note \(E^{F_2 F_1}_\sigma(s,s') \vdash D^{F_2 F_1}_\sigma(s)\).
	Applying the translation \(F_2\) to \axiom[\chi]{TM3} implies \(D^{F_2
	F_1}_\sigma(s) \vdash \exists w^{F_2 G_1\sigma} F_2\chi_\sigma(s,w)\);
	therefore the cut rule implies \(E^{F_2 F_1}_\sigma(s,s') \vdash \exists
	w^{F_2 G_1 \sigma} F_2\chi_\sigma(s,w)\). On the other hand,
	\axiom[\chi]{TM1} and \(F_2\) imply \(F_2\chi_\sigma(s,w) \vdash D^{F_2
	G_1}_\sigma(w)\). By Lemma~\ref{lemma:t-map-domain-lemma} and the cut rule, \(F_2\chi_\sigma(s,w) \vdash D^{F_2}_{G_1\sigma}(w)\). Applying
	\axiom[\eta]{TM3} and the cut rule gives \(F_2\chi_\sigma(s,w) \vdash \exists
	{t'}^{G_2 G_1\sigma} \eta_{G_1\sigma}(w,t')\). Finally, we have by
	\axiom[\eta]{TM5} the following sequent.%
	\[%
		\eta_{F_1\sigma}(s',v) \wedge \eta_{G_1\sigma}(w,t') \wedge
		F_2\chi_\sigma(s',w) \vdash G_2\chi_\sigma(v,t').
	\]%
	We use the consequent and \(G_2\) applied to \axiom[\chi]{TM4} to deduce
	\(G_2 \chi_\sigma(v,t') \wedge G_2\chi_\sigma(v,t) \vdash E^{G_2
	G_1}_\sigma(t,t')\). Gathering together these sequents, with the cut rule, we
	deduce the desired sequent.
\end{proof}

\begin{proof}%
	[Proof of Proposition~\ref{prop:horizontal-composition-yields-tmap}]
	This proof consists of five parts, one for each property of a t-map. 
        
	\emph{Proving \axiom{TM1}:}%
	\[%
			(\eta\circ\chi)_\sigma(s,t) \vdash D_\sigma^{F_2 F_1}(s) \wedge
			D_\sigma^{G_2 G_1}(t).
	\]%
	Note \(E^{F_2 F_1}_\sigma(s,s') \vdash D^{F_2 F_1}_\sigma(s)\) and \(E^{G_2
	G_1}_\sigma(t',t) \vdash D^{G_2 G_1}_\sigma(t)\). Therefore \axiom{TM1} is a
	consequence of Lemma~\ref{lemma:t-map-Z-lemma}, since both \(E^{F_2
	F_1}_\sigma(s,s')\) and \(E^{G_2 G_1}_\sigma(t',t)\) appear in one of the two
	formulae presenting \(\pwrap{\eta\circ\chi}_\sigma\).

	\emph{Proving \axiom{TM2}:}%
	\[%
			E_\sigma^{F_2 F_1}(s_1,s_2) \wedge E_\sigma^{G_2 G_1}(t_1,t_2) \wedge
			(\eta\circ\chi)_\sigma(s_1,t_1) \vdash (\eta\circ\chi)_\sigma(s_2,t_2).    
	\]%
	Transitivity of equality implies, after applying the translation \(G_2
	G_1\),%
	\[%
		E^{G_2 G_1}_\sigma(t_1,t_2) \wedge E^{G_2 G_1}_\sigma(t_1,t') \vdash E^{G_2
		G_1}_\sigma(t_2,t').
	\]%
	By the cut rule, this implies \(E^{G_2 G_1}_\sigma(t_1,t_2) \wedge
	\pwrap{\eta\circ\chi}_\sigma(s_1,t_1) \vdash
	\pwrap{\eta\circ\chi}_\sigma(s_1,t_2)\). Therefore it suffices to show that
	the following sequent is provable in \(T_3\).%
	\[%
		E^{F_2 F_1}_\sigma(s_1,s_2) \wedge \pwrap{\eta\circ\chi}_\sigma(s_1,t_2)
		\vdash \pwrap{\eta\circ\chi}_\sigma(s_2,t_2)
	\]%
	This follows from the following proof tree.
	\begin{prooftree}
		\AxiomC{}
		\RightLabel{\axiom[\chi]{TM2}}
		\UI$ E^{F_1}_\sigma(s_1',s_2') \wedge \chi_\sigma(s_1',w') \fCenter
		\chi_\sigma(s_2',w') $
		\RightLabel{$F_2$}
		\UI$ E^{F_2 F_1}_\sigma(s_1,s_2) \wedge F_2\chi_\sigma(s_1,w) \fCenter
		F_2\chi_\sigma(s_2,w) $
		\RightLabel{$\wedge$-intro., $\exists$-intro.}
		\UI$ E^{F_2 F_1}_\sigma(s_1,s_2) \wedge \pwrap{\eta\circ\chi}_\sigma(s_1,t_2)
		\fCenter \pwrap{\eta\circ\chi}_\sigma(s_2,t_2) $
	\end{prooftree}

	\emph{Proving \axiom{TM3}:}%
	\[%
			D^{F_2 F_1}_\sigma(s) \vdash \exists t^{G_2 G_1 \sigma}
			\pwrap{\eta\circ\chi}_\sigma(s,t).
	\]%
	Begin with \axiom[\chi]{TM3} and apply \(F_2\) to deduce \(D^{F_2
	F_1}_\sigma(s) \vdash \exists w^{F_2 G_1\sigma} F_2\chi_\sigma(s,w)\). The
	cut rule, \axiom[\chi]{TM1}, \(F_2\), and
	Lemma~\ref{lemma:t-map-domain-lemma} imply \(F_2\chi_\sigma(s,w) \vdash
	D^{F_2}_{F_1\sigma}(w)\). \axiom[\eta]{TM3} and cut gives
	\(F_2\chi_\sigma(s,w) \vdash \exists t^{G_2 G_1\sigma} F_2\chi_\sigma(s,w)
	\wedge \eta_{G_1\sigma}(s,t)\). On the other hand, \(F_2\chi_\sigma(s,w)
	\vdash D^{F_2 G_1}_\sigma(w) \equiv E^{F_2 G_1}_\sigma(w,w)\), via \(F_2\)
	applied to \axiom[\chi]{TM1}. Furthermore, \axiom[\eta]{TM5} implies
	\(\eta_{G_1\sigma}(w,t) \wedge E^{F_2 G_1}_\sigma(w,w) \vdash E^{G_2
	G_1}_\sigma(t,t)\). By cut, this gives%
	\[%
		D^{F_2 F_1}_\sigma(s) \vdash \exists t^{G_2 G_1\sigma} \exists w^{F_2
		G_1\sigma} \pwrap{%
			F_2\chi_\sigma(s,w) \wedge \eta_{G_1\sigma}(w,t)
		} \wedge E^{G_2 G_1}_\sigma(t,t).%
	\]%
	Using \(\exists\)-introduction and elimination, we conclude%
	\begin{align*}%
		D^{F_2 F_1}_\sigma(s) \vdash &\exists t^{G_2 G_1\sigma}
			\exists {t'}^{G_2 G_1\sigma} \exists w^{F_2 G_1\sigma} \pwrap{%
				F_2\chi_\sigma(s,w) \wedge \eta_{G_1\sigma}(w,t')
			} \wedge E^{G_2 G_1}_\sigma(t',t).%
	\end{align*}%
	The consequent is \(\exists t^{G_2
	G_1\sigma}\pwrap{\eta\circ\chi}_\sigma(s,t)\), as desired.

	\emph{Proving \axiom{TM4}:}%
	\[%
			(\eta\circ\chi)_\sigma(s,t) \wedge (\eta\circ\chi)_\sigma(s,w) \fCenter
			E_\sigma^{G_2 G_1}(t,w).
	\]%
	Apply \(F_2\) to \axiom[\chi]{TM4} to get \(F_2\chi_\sigma(s,w_1) \wedge
	F_2\chi_\sigma(s,w_2) \vdash E^{F_2 G_1}_\sigma(w_1,w_2)\). By
	\axiom[\eta]{TM5}, \(\eta_{G_1\sigma}(w_1,t_1') \wedge
	\eta_{G_1\sigma}(w_2,t_2') \wedge E^{F_2 G_1}_\sigma(w_1,w_2) \vdash E^{G_2
	G_1}_\sigma(t_1',t_2')\). Since \(E^{G_2 G_1}_\sigma\) is transitive,%
	\[%
		E^{G_2 G_1}_\sigma(t_1',t_1) \wedge E^{G_2 G_1}_\sigma(t_1',t_2') \wedge
		E^{G_2 G_1}_\sigma(t_2',t_2) \vdash E^{G_2 G_1}_\sigma(t_1,t_2).
	\]%
	We deduce \axiom[\eta\circ\chi]{TM4} using the cut rule.

	\emph{Proving \axiom{TM5}:} Let \(\phi \hookrightarrow \vec\sigma\) be a
	\(T_1\)\-/formula, where \(\vec\sigma \equiv \sigma_1,\ldots,\sigma_n\). We
	need to show \(\pwrap{\eta\circ\chi}_{\vec\sigma}(\vec s, \vec t) \wedge F_2
	F_1 \phi(\vec s) \vdash G_2 G_1 \phi(\vec t)\).

	Apply \(F_2\) to \axiom[\chi]{TM5} to deduce \(\bigwedge_{i=1}^n
	F_2\chi_\sigma(s_i,w_i) \wedge F_2 F_1 \phi(\vec s) \vdash F_2 G_1
	\phi(\vec{w})\). Using \axiom[\eta]{TM5}, \(\bigwedge_{i=1}^n \eta_{G_1
	\sigma_i}(w_i, t_i') \wedge F_2 G_1 \phi(\vec w) \vdash G_2 G_1
	\phi(\vec{t'})\). Applying \(G_2 G_1\) to \(=\)\-/introduction yields
	\(E^{G_2 G_1}_\sigma(\vec{t_1'},\vec{t}) \vdash G_2 G_1\phi(\vec{t'}) \vdash
	G_2 G_1 \phi(\vec{t})\). Therefore, by the cut rule \(
	\pwrap{\eta\circ\chi}_{\vec\sigma}(\vec s, \vec t) \wedge F_2 F_1\phi(\vec s)
	\vdash G_2 G_1\phi(\vec t) \).
\end{proof}

\begin{proof}%
	[Proof of Theorem~\ref{theorem:godements-law}]
	We first expand the definitions of \(A \equivdef ((\eta_2 \cdot
	\chi_2)\circ(\eta_1\cdot\chi_1))_\sigma\) and
	\(B \equivdef ((\eta_2\circ\eta_1)\cdot(\chi_2\circ\chi_1))_\sigma\) below.%
	\begin{align*} 
			A_\sigma&(s,t) \equiv%
			\exists {t'}^{H_2 H_1\sigma} \exists a^{F_2 H_1\sigma}%
				F_2(\eta_1\cdot\chi_1)_\sigma(s,a) \wedge
				(\eta_2\cdot\chi_2)_{H_1\sigma}(a,t')
			\wedge E^{H_2 H_1}_\sigma(t',t) \\
			&\dashv\vdash \exists {t'}^{H_2 H_1\sigma} \exists a^{F_2 H_1 \sigma}
			\exists b^{F_2 G_1\sigma} \exists c^{G_2 H_1 \sigma}%
			F_2 (\chi_1)_{\sigma}(s,b)\\
			&\quad\wedge F_2(\eta_1)_\sigma(b,a)
			\wedge (\chi_2)_{H_1\sigma}(a,c) \wedge (\eta_2)_{H_1\sigma}(c,t') \wedge%
			E^{H_2 H_1}_\sigma(t',t)\\
			B_\sigma&(s,t) \equiv%
			\exists \alpha^{G_2 G_1\sigma} \pwrap{%
				(\chi_2\circ\chi_1)_\sigma(s,\alpha) \wedge
				(\eta_2\circ\eta_1)_\sigma(\alpha,t)
			}\\%
			&\dashv\vdash \exists \alpha^{G_2 G_1\sigma} \exists {\alpha'}^{G_2
			G_1\sigma}%
				\exists \beta^{F_2 G_1\sigma} %
					F_2 (\chi_1)_\sigma(s,\beta) \wedge (\chi_2)_{G_1\sigma}(\beta,\alpha')
			\wedge E^{G_2 G_1}_\sigma(\alpha',\alpha)\\ 
                &\quad\wedge \exists {t'}^{H_2
				H_1\sigma} \Big(%
					\exists \gamma^{G_2 H_1\sigma} \big(%
						G_2 (\eta_1)_\sigma(\alpha,\gamma)
						\wedge (\eta_2)_{H_1\sigma}(\gamma,t')
					\big) \wedge E^{H_2 H_1}_\sigma(t',t)%
			\Big)\\
			&\dashv\vdash \exists {t'}^{H_2 H_1\sigma} \exists \alpha^{G_2 G_1
			\sigma} \exists \beta^{F_2 G_1\sigma} \exists \gamma^{G_2 H_1\sigma} F_2
			(\chi_1)_{\sigma}(s,\beta)\\
			&\quad\wedge (\chi_2)_{G_1\sigma}(\beta,\alpha) \wedge
			G_2(\eta_1)_{\sigma}(\alpha,\gamma) \wedge (\eta_2)_{H_1\sigma}(\gamma,t')
			\wedge E^{H_2 H_1}_\sigma(t',t)
	\end{align*}%
	To keep the notation clean, we will omit the sort symbols from the
	t-map formulae. We do not lose any information by doing this since the domains
	of the variables indicate the appropriate domains for the t-maps.
	We begin with the forward implication. Both sides contain \(E^{H_2
	H_1}_\sigma(t',t)\), and when replacing \(\beta\) with \(b\) and \(\gamma\)
	with \(c\), both sides also contain some of the same conjuncts. Therefore it
	suffices to show that the following sequent is provable in \(T_3\).
	\[%
		F_2\chi_1(s,b) \wedge F_2\eta_1(b,a) \wedge \chi_2(a,c) \wedge \eta_2(c,t')
		\vdash \exists \alpha \pwrap{\chi_2(b,\alpha) \wedge G_2\eta_1(\alpha,c)}
	\]%
	We obtain the original sequent \(A \vdash
	B\) from the above sequent by
	reintroducing the duplicate conjuncts and quantifying over \(b,c\) and
	\(t'\).

	Akin to the proof of
	Proposition~\ref{prop:horizontal-composition-yields-tmap}, we have the
	following chain of sequents.%
	\[%
		F_2\chi_1(s,b) \vdash D^{F_2 G_1}_\sigma(b) \vdash D^{F_2}_{G_1\sigma}(b)
		\vdash \exists \alpha^{G_2 G_1 \sigma} \chi_2(b,\alpha)
	\]%
	Therefore it suffices to prove:%
	\[%
		F_2 \chi_1(s,b) \wedge F_2\eta_1(b,a) \wedge \chi_2(a,c) \wedge
		\eta_2(c,t') \wedge \chi_2(b,\alpha) \vdash G_2\eta_1(\alpha,c).
	\]%
	This is a consequence of \axiom[\chi_2]{TM5}:%
	\[%
		\chi_2(a,c) \wedge \chi_2(b,\alpha) \wedge F_2\eta_1(b,a) \vdash
		G_2\eta_1(\alpha,c).
	\]%
	We proceed to the converse implication. It suffices to prove:%
	\begin{gather*}
		\exists \alpha \exists \beta \exists \gamma \, F_2\chi_1(s,\beta) \wedge
		\chi_2(\beta,\alpha) \wedge G_2\eta_1(\alpha,\gamma) \wedge
		\eta_2(\gamma,t') \wedge E^{H_2 H_1}_\sigma(t',t)\\%
		\vdash\exists t'' \exists a \exists b \exists c \, F_2\chi_1(s,b) \wedge
		F_2\eta_1(b,a) \wedge \chi_2(a,c) \wedge \eta_2(c,t'') \wedge E^{H_2
		H_1}_\sigma(t'',t).
	\end{gather*}
	We claim that the following sequents are provable in \(T_3\).
	\begin{gather}
		F_2 \chi_1(s,\beta) \wedge \chi_2(\beta,\alpha) \wedge
		G_2\eta_1(\alpha,\gamma) \wedge \eta_2(\gamma,t')\label{sequent:godement-1}\\
		\vdash \exists a \exists c \, F_2\chi_1(s,\beta) \wedge F_2\eta_1(\beta,a)
		\wedge \chi_2(a,c) \notag\\
		\chi_2(\beta,\alpha) \wedge \chi_2(a,c) \wedge F_2\eta_1(\beta,a) \vdash
		G_2\eta_1(\alpha,c) \label{sequent:godement-2}\\
		G_2\eta_1(\alpha,\gamma) \wedge G_2\eta_1(\alpha,c) \vdash E^{G_2
		H_1}_\sigma(\gamma,c) \label{sequent:godement-3}\\
		\eta_2(\gamma,t') \wedge \eta_2(c,t'') \wedge E^{G_2 H_1}_\sigma(\gamma,c)
		\vdash E^{H_2 H_1}_\sigma(t',t'') \label{sequent:godement-4}\\
		E^{H_2 H_1}_\sigma(t',t) \wedge E^{H_2 H_1}_\sigma(t',t'') \vdash E^{H_2
		H_1}_\sigma(t'',t). \label{sequent:godement-5}
	\end{gather}
	Applying the cut rule to (\ref{sequent:godement-1}) through
	(\ref{sequent:godement-5}), along with \(\exists\)\-/introduction and
	elimination, we can deduce the sequent%
	\begin{gather*}
		F_2 \chi_1(s,\beta) \wedge \chi_2(\beta,\alpha) \wedge
		G_2\eta_1(\alpha,\gamma) \wedge \eta_2(\gamma,t') \wedge E^{H_2
		H_1}_\sigma(t',t)\\
		\vdash \exists t'' \exists a \exists c \, F_2\chi_1(s,\beta) \wedge
		F_2\eta_1(\beta,a) \wedge \chi_2(a,c) \wedge \eta_2(c,t'') \wedge E^{H_2
		H_1}_\sigma(t'',t).
	\end{gather*}
	Quantifying over \(\beta\) and renaming it to \(b\) yields%
	\begin{gather*}
		F_2\chi_1(s,\beta) \wedge \chi_2(\beta,\alpha) \wedge
		G_2\eta_1(\alpha,\gamma) \wedge \eta_2(\gamma,t') \wedge E^{H_2
		H_1}_\sigma(t',t)\\
		\vdash \exists t'' \exists a \exists b \exists c F_2\chi_1(s,b) \wedge
		F_2\eta_1(b,a) \wedge \chi_2(a,c) \wedge \eta_2(c,t'') \wedge E^{H_2
		H_1}_\sigma(t'',t).
	\end{gather*}
	The desired converse implication is now a consequence of
	\(\exists\)\-/introduction applied to \(\alpha,\beta,\gamma,\) and \(t'\).

	All that remains is proving (\ref{sequent:godement-1}) through
	(\ref{sequent:godement-5}). (\ref{sequent:godement-2}) is a consequence of
	\axiom[\chi_2]{TM5}. Applying \(G_2\) to \axiom[\eta_1]{TM4} yields
	(\ref{sequent:godement-3}). \axiom[\eta_2]{TM5} implies
	(\ref{sequent:godement-4}). Since \(E^{H_2 H_1}_\sigma\) is transitive and
	symmetric, we deduce (\ref{sequent:godement-5}). This leaves
	(\ref{sequent:godement-1}).
	\begin{prooftree}
		\AxiomC{}
		\RightLabel{\axiom{(1)}}
		\UI$F_2 \chi_1(s,\beta) \fCenter D^{F_2 G_1}_\sigma(\beta)$
		\RightLabel{(2)}
		\UI$F_2 \chi_1(s,\beta) \fCenter \exists a\, F_2\eta_1(\beta,a)$
		
		\AxiomC{}
		\RightLabel{\axiom{(3)}}
		\UI$F_2 \eta_1(\beta,a) \fCenter D^{F_2}_{H_1\sigma}(a)$
		\RightLabel{\axiom{(4)}}
		\UI$F_2 \eta_1(\beta,a) \fCenter \exists c\, \chi_2(a,c)$
		\BI$ F_2\chi_1(s,\beta) \fCenter \exists a \exists c \, F_2\eta_1(\beta,a)
		\wedge \chi_2(a,c)$
		\UI$ F_2\chi_1(s,\beta) \fCenter \exists a \exists c \, F_2\chi_1(s,\beta)
		\wedge F_2\eta_1(\beta,a) \wedge \chi_2(a,c)$
	\end{prooftree}
        where (1) is \(F_2\pwrap{\axiom[\chi_1]{TM1}}\), \(\wedge\)\-/elim.; 
	(2) is \(F_2\pwrap{\axiom[\eta_1]{TM3}}\), cut.; 
	(3) is \(F_2\pwrap{\axiom[\eta_1]{TM1}}\), \(\wedge\)\-/elim., 
	Lemma~\ref{lemma:t-map-domain-lemma}, cut.; and (4) is \axiom[\chi_2]{TM3}, cut.
	Finally, (\ref{sequent:godement-1}) follows from \(\wedge\)\-/introduction
	and the cut rule.
\end{proof}

\begin{proof}%
	[Proof of Proposition~\ref{prop:horizontal-composition-unital}]
	We wish to show \(\One^G \circ \One^F = \One^{GF}\). Unpacking the definition
	of horizontal composition, this means we need to prove%
	\[%
		\exists {t'}^{GF\sigma} \exists w^{GF\sigma} \pwrap{%
			E^{GF}_\sigma(s,w) \wedge E^G_{F\sigma}(w,t')
		} \wedge E^{GF}_\sigma(t',t) \dashv\vdash E^{GF}_\sigma(s,t).%
	\]%
	We begin with the converse direction. By
	Lemma~\ref{lemma:t-map-domain-lemma}, \(E^{GF}_\sigma(s,t) \vdash
	E^{GF}_\sigma(t,t) \equiv D^{GF}_\sigma(t) \vdash D^G_{F\sigma}(t) \equiv
	E^G_{F\sigma}(t,t)\). Therefore,
	\begin{align*}
		E^{GF}_\sigma(s,t) &\vdash E^{GF}_\sigma(s,t) \wedge E^G_{F\sigma}(t,t)
		\wedge E^{GF}_\sigma(t,t)\\
		&\vdash \exists {t'}^{GF\sigma} \exists w^{GF\sigma} \, \pwrap{%
			E^{GF}_\sigma(s,w) \wedge E^G_{F\sigma}(w,t')
		} \wedge E^{GF}_\sigma(t',t).%
	\end{align*}
	We proceed to the forward direction. Note that \(E^{GF}_\sigma(s,w) \vdash
	D^{GF}_\sigma(w)\), so by Lemma~\ref{lemma:t-map-equality-lemma} and the cut
	rule, \(E^{GF}_\sigma(s,w) \wedge E^G_{F\sigma}(w,t') \vdash
	E^{GF}_\sigma(s,t')\). Since \(E^{GF}_\sigma\) is transitive,
	\(E^{GF}_\sigma(s,t') \wedge E^{GF}_\sigma(t',t) \vdash E^{GF}_\sigma(s,t)\).
	By \(\exists\)\-/introduction and the cut rule, we conclude.
\end{proof}

\subsection{Coherence Proofs}\label{subsection:coherence-proofs}

\begin{proof}[Proof of Proposition~\ref{prop:C-trivial-compositor}]
	The key part of this proof is
	Proposition~\ref{prop:reconstrual-strictly-composes}. Given an object
	\(\fb\phi\) of \(\C(T_1)\), we have%
	\[%
		\C(GF)\fb\phi \equiv \pwrap{GF}^+ \fb\phi \equiv G^+ F^+ \fb\phi \equiv \C(G)\C(F)\fb\phi.
	\]%
	For a morphism \(\theta: \fb\phi \to \fb\psi\) in \(\C(T_1)\), the above
	argument applied to a substitution class presenting \(\theta\) shows that
	\(\C(GF)\theta = \C(G)\C(F)\theta\). In particular, \(\C(G)\C(F)\) and
	\(\C(GF)\) are identical functors. Therefore we define the compositor
	\(\C_{GF}: \C(G)\C(F) \Rightarrow \C(GF)\) to be the identity 2\-/cell
	\(\One^{\C(GF)}\).

	Naturality of the compositor is the equation \(\C_{G_2 G_1} \cdot
	(\C(\eta)\circ\C(\chi)) = \C(\eta\circ\chi) \cdot \C_{F_2 F_1}\) for any pair of\
	t\-/maps \(\chi: F_1 \Rightarrow G_1\) and \(\eta: F_2 \Rightarrow G_2\) in
	\(\CThEq\). Since the compositor's components are identity 2\-/cells, this
	reduces to \(\C(\eta)\circ\C(\chi) = \C(\eta\circ\chi)\). It suffices to 
        verify that \(\C(\eta\circ\chi)_{\fb\phi} \dashv\vdash (\C(\eta)\circ\C(\chi))_{\fb\phi}\)
	for any object \(\fb\phi \hookrightarrow \sigma\) of \(\C(\Dom F_1)\). From
	Remark~\ref{remark:tmap-horizontal-composition-ep} and applying \(\C\), we find%
	\[%
		\C(\eta\circ\chi)_{\fb\phi}(x,y) \equiv%
			\exists z^{F_2 G_1 \sigma} \pwrap{F_2 \chi_\sigma(x,z) \wedge \eta_{G_1
			\sigma}(z,y)} \wedge D_\sigma^{G_2 G_1}(y) \wedge F_2 F_1\phi(x).
	\]%
	Applying TM1(\(\eta\)) allows us to drop \(D_\sigma^{G_2 G_1}\):%
	\[%
		\C(\eta\circ\chi)_{\fb\phi}(x,y) \dashv\vdash \exists z^{F_2 G_1\sigma}
		\pwrap{ F_2\chi_\sigma(x,z) \wedge \eta_{G_1\sigma}(z,y) } \wedge F_2
		F_1\phi(x).
	\]%
	On the other hand, \(\pwrap{ \C(\eta) \circ \C(\chi) }_{\fb\phi} =
	\C(\eta)_{\C(G_1)\fb\phi} \cdot \C(F_2)\C(\chi)_{\fb\phi}\), so the
	expression \(\pwrap{\C(\eta)\circ\C(\chi)}_{\fb\phi}(x,y)\) is identical to%
	\[%
			\exists z^{F_2 G_1 \sigma} F_2\chi_\sigma(x,z) \wedge F_2 F_1 \phi(x)
			\wedge \eta_{G_1\sigma}(z,y) \wedge F_2 G_1\phi(z).
	\]%
	Applying TM5(\(\chi\)) and since \(F_2\) is a translation, we can
	drop \(F_2 G_1 \phi(z)\); hence%
	\[%
		\pwrap{\C(\eta)\circ\C(\chi)}_{\fb\phi}(x,y) \dashv\vdash \exists z^{F_2
		G_1\sigma} \pwrap{ F_2\chi_\sigma(x,z) \wedge \eta_{G_1 \sigma}(z,y) }
		\wedge F_2 F_1\phi(x).
	\]%
	As shown earlier, this is logically equivalent to the formula presenting
	\(\C(\eta\circ\chi)_{\fb\phi}\), so the components of both natural
	transformations are logically equivalent.
\end{proof}

\begin{proof}[Proof of Proposition~\ref{prop:C-identitor}]
	Let \(\theta_\phi(x,y) \equivdef \phi(x) \wedge x = y\). We need to verify
	that this formula presents a morphism \(\fb\phi \to \C(1_T)\fb\phi\) in
	\(\C(T)\). This means we need to show \axiom[\theta_\phi]{DM1} through
	\axiom[\theta_\phi]{DM3} are provable in \(T\). \axiom[\theta_\phi]{DM1} is
	the sequent \(\theta(x,y) \vdash \phi(x) \wedge 1_T^+\phi(y)\). This follows
	from \(\theta_\phi(x,y) \vdash \phi(x)\) and
	Proposition~\ref{prop:reconstrual-log-eq-to-identity}, along with
	\(=\)-elimination. \axiom[\theta_\phi]{DM2} is the sequent
	\(\theta_\phi(x,y_1) \wedge \theta_\phi(x,y_2) \vdash y_1 = y_2\), which
	follows from transitivity of equality. \axiom[\theta_\phi]{DM3} is the sequent
	\(\phi(x) \vdash \exists y \theta_\phi(x,y)\), which is an application of
	\(=\)-introduction.

	This shows that the components of the proposed identitor are indeed morphisms
	in \(\C(T)\). Furthermore, each morphism \(\theta_\phi\) has an inverse,
	presented by \(\phi(y) \wedge x = y\). What remains is showing that the
	proposed identitor is indeed a natural transformation. Suppose we have a
	morphism \(\eta: \fb\phi \to \fb\phi\) in \(\C(T)\). Naturality of the
	identitor is the equation \(\theta_\phi \circ \eta = \C(1_T)\eta \circ
	\theta_\phi\). This reduces to the sequent \(\eta(x,y) \wedge \psi(y)
	\dashv\vdash \phi(x) \wedge 1_T^+\eta(x,y)\).  Due to
	Proposition~\ref{prop:reconstrual-log-eq-to-identity} and \axiom[\eta]{DM1},
	this sequent is provable in \(\C(T)\). Thus the family of morphisms
	\(\theta_\phi\) indeed define the components of a homotopy \(1_{\C(T)}
	\Rightarrow \C(1_T)\).
\end{proof}

\begin{proof}[Proof of Proposition~\ref{prop:C-pseudofunctor}]
	We need to verify the hexagon and triangle identities for \axiom[\C]{PF5}.
	Since \(\C\) has a trivial compositor and \(\Coh\) and \(\CThEq\) have trivial
	associators, the hexagon identity degenerates into the following triangle.
	\[%
		\begin{tikzpicture}%
			\node[
				regular polygon,
				regular polygon sides=3,
				minimum width=30mm,
				rotate=90,
			] (TG) {}
			(TG.corner 1) node (TG1) {\(\C(H)\C(G)\C(F)\)}
			(TG.corner 2) node (TG2) {\(\C(H)\C(G)\C(F)\)}
			(TG.corner 3) node (TG3) {\(\C(H)\C(G)\C(F)\)}
			;
			\draw[->]%
				(TG1)--(TG3) node[midway, above left] {\scriptsize\(\C_{HG}\circ\One^{\C(F)}\)};
			\draw[->]%
				(TG3)--(TG2) node[midway, right] {\scriptsize\(\One\)};
			\draw[->]%
				(TG1)--(TG2) node[midway, below left] {\scriptsize\(\One^{\C(H)}\circ\C_{GF}\)};
		\end{tikzpicture}%
	\]%
	Since all the compositors are trivial, each side of this diagram is the
	2\-/cell \(\One^{\C(H)\C(G)\C(F)}\), so the diagram commutes. The two square
	identities degenerate into the equations \(\One^{\C(F)}\circ\C_{1_{T_1}} =
	\C(r_F)\) and \(\C_{1_{T_2}}\circ\One^{\C(F)} = \C(l_F)\). These identities
	follow from the observation that, for any object \(\fb\phi\) of \(\C(T_1)\),
	\(\C(r_F)_{\fb\phi}\) and \(\C(l_F)_{\fb\phi}\) are presented by the formula
	\(F\phi(x) \wedge x = y\).
\end{proof}

\begin{proof}[Proof of Proposition~\ref{prop:T-compositor}]
	We first need to show that \(\kappa_{\Gee\Eff}\) is a homotopy from
	\(\T(\Gee)\T(\Eff)\) to \(\T(\Gee\Eff)\). Since \(\T(\Gee)\), \(\T(\Eff)\),
	and \(\T(\Gee\Eff)\) are e.p.\ translations with trivial domain classes,
	\axiom[\kappa]{TM1} through \axiom[\kappa]{TM4}, \axiom[\kappa]{TM6}, and
	\axiom[\kappa]{TM7} are provable trivially. This leaves \axiom[\kappa]{TM5}
	and \axiom[\kappa]{TM8}; i.e., we must show that for any \(\T(C_1)\)-formula
	\(\phi\hookrightarrow\vec\sigma\), the following sequents are
	provable:%
	\begin{gather*}
		(\kappa_{\Gee\Eff})_{\vec\sigma}(x,y) \wedge \T(\Gee\Eff)\phi(x) \vdash
		\T(\Gee)\T(\Eff)\phi(y),\\
		(\kappa_{\Gee\Eff})_{\vec\sigma}(x,y) \wedge \T(\Gee)\T(\Eff)\phi(y) \vdash
		\T(\Gee\Eff)\phi(x).
	\end{gather*}%
	Since \((\kappa_{\Gee\Eff})_{\vec\sigma}(x,y) \equiv x = y\), this is equivalent
	to showing the sequents%
	\[%
		\T(\Gee\Eff)\phi(x) \dashv\vdash \T(\Gee)\T(\Eff)\phi(x).
	\]%
	As translations satisfy the same reconstrual laws, it suffices to
	verify two cases:%
	\begin{gather*}
		\T(\Gee\Eff)R(x) \dashv\vdash \T(\Gee)\T(\Eff)R(x),\\
		\T(\Gee\Eff)f(x,y) \dashv\vdash \T(\Gee)\T(\Eff)f(x,y),
	\end{gather*}%
	where \(R\) and \(f\) are an arbitrary relation symbol and function symbol
	respectively. The only relations in the internal language of a coherent
	category are equality relations. Indeed the above sequent holds for the case
	\(R(x_1,x_2) \equiv x_1 = x_2\) because \(E^{\T(\Eff)}(x_1,x_2) \equiv x_1 =
	x_2\) for any coherent functor \(\Eff\). This leaves the case of an arbitrary
	function symbol \(\ul{f}: \ul{X} \to \ul{Y}\) in \(\T(C_1)\). For this case,
	note that \(\T(\Gee\Eff)\ul{f}(x,y) \equiv \ul{\Gee\Eff f}(x) = y\). Then,%
	\begin{align*}%
		\T(\Gee)\T(\Eff)\ul{f}(x,y) &\equiv \T(\Gee)\fb{\ul{\Eff f}(x') = y'}(x,y)\\
		&\equiv \exists t\, \T(\Gee)\ul{\Eff f}(x,t) \wedge t = y\\
		&\equiv \exists t\, \ul{\Gee\Eff f}(x) = t \wedge t = y\\
		&\dashv\vdash \ul{\Gee\Eff f}(x) = y,
	\end{align*}
	so the desired sequent holds for \(\ul{f}\). Therefore \(\kappa_{\Gee\Eff}\)
	is a homotopy.

	Now we need to verify that \(\kappa\) is natural, i.e., for any pair of
	natural transformations \(\eta_1: \Eff_1 \Rightarrow \Gee_1\) and \(\eta_2:
	\Eff_2 \Rightarrow \Gee_2\) between coherent functors \(\Eff_1,\Gee_1: C_1 \to
	C_2\) and \(\Eff_2,\Gee_2: C_2 \to C_3\), we must have the equation
	\(\kappa_{\Gee_2\Gee_1} \cdot \pwrap{\T(\eta_2)\circ\T(\eta_1)} =
	\T(\eta_2\circ\eta_1) \cdot \kappa_{\Eff_2\Eff_1}\). The \(\ul{X}\) component
	of the left side is the following formula.%
	\[%
		\exists t_2 \kappa_{\Gee_2\Gee_1}(t_1,y) \wedge \exists t_1
		\T(\Eff_2)\T(\eta_1)_{\ul{X}}(x,t_1) \wedge \T(\eta_2)_{\ul{\Gee_1 X}}(t_1,t_2)
		\wedge D^{\T(\Gee_2)\T(\Gee_1)}_{\ul{X}}\!(y)
	\]%
	After expanding all terms to their definitions and applying \(=\)-introduction
	and elimination, the above formula is logically equivalent to%
	\[%
		\ul{(\eta_2)_{\Gee_1 X}}\pwrap{\ul{(\Eff_2\eta_1)_X}(x)} = y \dashv\vdash
		\exists t\, \ul{(\eta_2)_{\Gee_1 X}}\pwrap{\ul{(\Eff_2\eta_1)_{X}}(t)} = y
		\wedge x = t.
	\]%
	The latter is the \(\ul{X}\) component of the t\-/map
	\(\T(\eta_2\circ\eta_1)\cdot\kappa_{\Eff_2\Eff_1}\), so \(\kappa\) is natural.
\end{proof}

\begin{proof}[Proof of Proposition~\ref{prop:T-identitor}]
	We first verify that, for any coherent category \(C\), \(1_{\T(C)}\) and
	\(\T(1_C)\) are identical translations. For a sort \(\ul{X}\) of \(\T(C)\),
	\(\T(1_C)\ul{X} \equiv \ul{1_C X} \equiv \ul{X} \equiv 1_{\T(C)}\ul{X}\), and
	\(E^{\T(1_C)}_{\ul{X}}(x_1,x_2) \equiv x_1 =_{\T(1_C)\ul{X}} x_2 \equiv x_1
	=_{\ul X} x_2 \equiv E^{1_{\T(C)}}_{\ul{X}}(x_1,x_2)\). For a function symbol
	\(\ul{f}: \ul{X} \to \ul{Y}\) of \(\T(C)\), \(\T(1_C)\ul{f}(x,y) \equiv
	\ul{1_C f}(x) = y \equiv \ul{f}(x) = y \equiv 1_{\T(C)}\ul{f}(x,y)\). Thus the
	underlying reconstruals of \(\T(1_C)\) and \(1_{\T(C)}\) are identical, so
	the two translations are identical. This justifies setting the identitor
	\(\T_{1_C}\) to the identity t\-/map \(\One^{\T(1_C)} \equiv
	\One^{1_{\T(C)}}\). This is a homotopy.
\end{proof}

\begin{proof}[Proof of Proposition~\ref{prop:T-pseudofunctor}]
	We need to verify \axiom[\T]{PF5}. Since \(\Coh\) and \(\CThEq\) have trivial
	associators, the hexagon identity degenerates into a
	diamond.%
	\[%
		\begin{tikzpicture}%
			\node[
				diamond,
				minimum width=50mm,
				minimum height=20mm
			] (DM) {}
			(DM.north) node (DM1) {\(\T(\Eta\Gee)\T(\Eff)\)}
			(DM.east) node (DM2) {\(\T(\Eta\Gee\Eff)\)}
			(DM.south) node (DM3) {\(\T(\Eta)\T(\Gee\Eff)\)}
			(DM.west) node (DM4) {\(\T(\Eta)\T(\Gee)\T(\Eff)\)}
			;
			\draw[->]%
				(DM4)--(DM1) node[midway, above left] {\scriptsize\(\T_{\Eta\Gee}\circ\One^\Eff\)};
			\draw[->]%
				(DM1)--(DM2) node[midway, above right] {\scriptsize\(\T_{(\Eta\Gee)\Eff}\)};
			\draw[->]%
				(DM4)--(DM3) node[midway, below left] {\scriptsize\(\One^\Eta \circ\T_{\Gee\Eff}\)};
			\draw[->]%
				(DM3)--(DM2) node[midway, below right] {\scriptsize\(\T_{\Eta(\Gee \Eff)}\)};
		\end{tikzpicture}%
	\]%
	Since all t\-/maps in this diagram are presented by \(\fb{x=y}\), its
	compositions are also presented by \(\fb{x=y}\); thus the above diagram commutes. Since \(\Coh\)
	has trivial unitors and \(\CThEq\) has a trivial identitor, the two square
	identities of \axiom[\T]{PF5} degenerate into the equations \(r_{\T(\Eff)} =
	\T_{\Eff1_{C_1}}\) and \(l_{\T(\Eff)} = \T_{1_{C_2}\Eff}\). All
	t\-/maps in these equations are presented by \(\fb{x=y}\), so the identities
	hold.
\end{proof}

\subsection{Biequivalence Proofs}
\label{subsection:biequivalence-proofs}

\begin{proof}[Proof of Lemma~\ref{lemma:ep-presentation}]
	We induct on the complexity of formulae.
	Using \axiom{IL1} and \axiom{IL2} axioms we can reduce any composition of
	function symbols to a single function symbol. Therefore the base case is when
	\(\phi\) is an atomic formula of the form \(R(f(\vec x))\), where \(f:
	\vec\sigma \to \tau\) is an \(n\)\-/ary function symbol and \(R\) a relation
	symbol in \(\Sigma\). In \(\C(T)\) we have a pullback square%
	\[%
		\begin{tikzcd}
			\fb{R(f(\vec x))} \arrow[r,"\minisub{\dom}{{\fb{R(f(
			\vec{x}))}}}"{yshift=.25em}] \arrow[d,"\psi_f",swap] &
			\fb{\vec\sigma} \arrow[d,"\theta_f"]\\
			\fb{R} \arrow[r,"\dom_{[R]}"] & \fb{\tau},
		\end{tikzcd}
	\]%
	where \(\psi_f\) is presented by \(\psi_f(\vec{x},y') \dashv\vdash
	\theta_f(\vec x, y') \wedge R\pwrap{f\pwrap{\vec x}}\). For ease of reference,
	let \(\Phi_f\) denote the induced morphism on the product
	\(\dom_{\fb{R(f(\vec{x}))}} \times \psi_f: \fb{R(f(\vec x))} \to
	\fb{\vec\sigma} \times \fb{R}\).  Since we have a pullback, for appropriate
	projections \(\pi_1\) and \(\pi_2\), the following diagram is an equalizer.%
	\[%
		\begin{tikzcd}%
			\fb{R(f(\vec x))} \arrow[r,"\Phi_f"] &%
			\fb{\vec\sigma} \times \fb{R} \arrow[r,"\theta_f \pi_1",shift left]%
			\arrow[r,"\minisub{\dom}{[R]} \pi_2",shift right,swap] &%
			\fb{\tau}
		\end{tikzcd}
	\]%
	By Rule~\ref{rule:term-reduction}, we know that%
	\[%
		\ep_T \fb{R(f(\vec x))} \equiv \fb{%
			\exists t^{\ul{\fb\tau}}\pwrap{%
				\ep_T f\pwrap{\vec x, t} \wedge \exists y^{\ul{\fb R}}\pwrap{%
					\ul{\dom_{\fb R}}(y) = t
				}%
			}%
		}.%
	\]%
	This is logically equivalent to%
	\[%
		\fb{%
			\exists z^{\ul{\fb{\vec\sigma}}}\pwrap{%
				\exists y^{\ul{\fb R}} \pwrap{%
					\bigwedge_{i=1}^n \ul{\pi^{\vec\sigma}_{\sigma_i}}(z) = x_i \wedge
					\ul{\theta_f}(z) = \ul{\dom_{\fb R}}(y)
				}%
			}%
		}.%
	\]%
	With \axiom[\theta_f\pi_1]{IL2} and \axiom[\dom_{\fb R}\pi_2]{IL2},
	\axiom{IL5} for the aforementioned equalizer implies%
	\[%
		\ul{\theta_f}\pwrap{\ul{\pi_1}(x)} = \ul{\dom_{\fb R}}\pwrap{\ul{\pi_2}(x)}
		\vdash \exists w^{\ul{\fb{R(f(\vec x))}}} \ul{\Phi_f}(w) = x.
	\]%
	We can further unpack the right side. \axiom{IL6} for the product
	\(\fb{\vec\sigma} \times \fb R\) allows us to relate the projections of
	\(\ul{\Phi_f}(w)\) and \(x\):%
	\[%
		\ul{\Phi_f}(w) = x \dashv\vdash \ul{\pi_1}\pwrap{\ul{\Phi_f}(w)} = \ul{\pi_1}(x)
		\wedge \ul{\pi_2}\pwrap{\ul{\Phi_f}(w)} = \ul{\pi_2}(x).
	\]%
	The universal property of the product \(\fb{\vec\sigma} \times \fb{R}\)
	implies \(\dom_{\fb{R(f(\vec x))}} = \pi_1 \Phi_f\) and \(\psi_f = \pi_2 \Phi_f\).
	We can use \axiom[\pi_1\Phi_f]{IL2} and \axiom[\pi_2\Phi_f]{IL2} to deduce%
	\[%
		\ul{\Phi_f}(w) = x \dashv\vdash \ul{\dom_{\fb{R(f(\vec x))}}}(w) =
		\ul{\pi_1}(x) \wedge \ul{\psi_f}(w) = \ul{\pi_2}(x).
	\]%
	We replace \(\ul{\Phi_f}(w) = x\) in an earlier sequent with the right side to
	deduce that the following sequent is provable in \(\T\C(T)\).
	\[%
		\ul{\theta_f}(z) = \ul{\dom_{\fb R}}(y) \vdash \exists
		w^{\ul{\fb{R(f(\vec{x}))}}} \pwrap{%
			\ul{\dom_{\fb{R(f(\vec x))}}}(w) = z \wedge \ul{\psi_f}(w) = y
		}%
	\]%
	Recalling \(\ep_T\fb{R(f(\vec x))}\), and expanding \(\ep_T f\),%
	\[%
		\ep_T \fb{R(f(\vec{x}))} \vdash \fb{\exists w^{\ul{\fb{R(f(\vec{x}))}}}
		\bigwedge_{i=1}^n \ul{\pi^{\vec{\sigma}}_{\sigma_i}}
		\left(\ul{\dom_{[R(f(\vec{x}))]}}(w)\right) = x_i}.
	\]%
	The \axiom{IL2} axiom for the commutative square of the aforementioned
	pullback along with the \axiom{IL6} axioms for the product \(\fb{\vec\sigma}
	\to \fb{\sigma_i}\) yield the converse sequent.

	We now show the inductive step. Suppose first that \(\phi(x)\) and \(\psi(y)\)
	are \(T\)\-/formulae such that \(\ep_T\fb\phi\) and \(\ep_T\fb\psi\) satisfy
	the lemma. We need to deduce%
	\[%
		\ep_T\fb{\phi(x) \wedge \psi(y)} \dashv\vdash%
		\fb{%
			\exists z^{\ul{\fb{\phi(x) \wedge \psi(y)}}} \bigwedge_{i=1}^n%
			\ul{\pi^{\vec\sigma}_{\sigma_i}} \pwrap{%
				\ul{\dom_{\fb{\phi(x) \wedge \psi(y)}}}(z)%
			} = x_i%
		}.%
	\]%
	This follows from Rule~\ref{rule:conjunction-preservation} and the \axiom{IL10}
	axioms for the intersection \(\fb\phi \wedge \fb\psi\) as subobjects of
	\(\fb{\Dom\phi} \times \fb{\Dom\psi}\). The disjunction (resp.\
	existential quantifier) case is due to Rule~\ref{rule:disjunction-preservation} and
        \axiom{IL8} axioms (resp.\ Rule~\ref{rule:quantifier-preservation} and
	\axiom{IL9} axioms).
\end{proof}
\begin{proof}[Proof of Proposition~\ref{prop:ep_F-definable-iso}]
	We prove that \((\ep_F)_\sigma : D^{\T(\C(F))\ep_{T_1}}_\sigma \to
	D^{\ep_{T_2} F}_\sigma\) and the map \((\ep^{-1}_F)_\sigma:
	D_\sigma^{\ep_{T_2}F} \to D_\sigma^{\T(\C(F)) \ep_{T_1}} \) (presented by
	\((\ep^{-1}_F)_\sigma(\vec{y}, x) \equivdef (\ep_F)_\sigma(x,\vec{y})\))
	satisfy the sequents of a definable map.

	We begin with \(\pwrap{\ep_F}_\sigma\). \axiom{DM2} is a consequence of the
	transitivity of the relation \(=\). \axiom{DM3} is an application of the
	introduction rules for \(=\) and \(\exists\). To prove \axiom{DM1}, we first
	show that \(D^{\T(\C(F))\ep_{T_1}}_\sigma\) is tautological. By
	Rule~\ref{rule:context-duplication}, \(D^{\ep_{T_1}}_\sigma \equiv
	\fb{E^{\ep_{T_1}}_\sigma(x',x')}\), which is logically equivalent to \(\fb{x'
	= x'}\) by Proposition~\ref{prop:epsilon-translation}. Thus,
	\(D_\sigma^{\T(\C(F)) \ep_{T_1}}(x) \equiv \T(\C(F)) D_\sigma^{\ep_{T_1}}(x)
	\dashv\vdash \T(\C(F)) [x' = x'](x) \dashv\vdash x = x\) by
	Definition~\ref{def:T-1-cells}. Lastly, apply
	Lemma~\ref{lemma:ep-presentation} to \(D^{\ep_{T_2}F}_\sigma(\vec{y})\) and
	observe that it is provable from the tautology \((\ep_F)_\sigma(x,\vec{y})
	\vdash (\ep_F)_\sigma(x,\vec{y})\).

	We now turn to \((\ep_F^{-1})_\sigma\). By symmetry, the proof of DM1 is
	analogous. \axiom{DM2} follows by the first \axiom{IL6} axiom applied to each
	conjunct, \axiom{IL3} applied to the monomorphism \(\dom_{D^F_\sigma}\), and
	the cut rule. \axiom{DM3} can be deduced from the
	tautology \(D^{\ep_{T_2} F}_\sigma(\vec{x}) \vdash D^{\ep_{T_2}
	F}_\sigma(\vec{x})\) and expanding \(D^{\ep_{T_2}F}_\sigma \equiv \ep_{T_2}
	D^F_\sigma\) using Lemma~\ref{lemma:ep-presentation}. 
\end{proof}
\begin{proof}[Proof of Proposition~\ref{prop:ep_F-homotopy}]
	By Propositions~\ref{prop:ep_F-definable-iso}, all the axioms of a homotopy
	are satisfied except for \axiom{TM5} and \axiom{TM8}. We show these below.

	Let \(\phi \hookrightarrow \vec\sigma\) be an \(n\)-ary \(T_1\)-formula.  Then
	for each \(\sigma_i\)  in \(\vec\sigma\), declare \(F\sigma_i\) as a list of
	\(T_2\)-sorts \(\vec{\tau}_i \equivdef \tau_{i1},\ldots,\tau_{im_i}\) for
	\(m_i \in \mathbb{Z}_+\). Hence \(F\phi \hookrightarrow \vec\tau\) is an
	\(m\)-ary \(T_2\)-substitution class for \(m = m_1 + \ldots + m_n\). We begin
	by showing \axiom{TM5}, i.e.,
	\(\pwrap{\ep_F}_{\vec\sigma}\pwrap{\vec{x},\vec{y}} \wedge
	\T\pwrap{\C(F)}\ep_{T_1} \phi(\vec{x}) \vdash \ep_{T_2} F \phi(\vec{y})\),
	where (by convention),
	\[%
		(\ep_{F})_{\vec{\sigma}}(\vec{x},\vec{y}) \equiv \bigwedge_{i=1}^n (\ep_{F})
		_{\sigma_i}(x_i, \vec{y_i}) \equiv \bigwedge_{i=1}^n \bigwedge_{j=1}^{m_i} 
		\ul{\pi^{\vec{\tau}_i}_{\tau_{ij}}}\left(\ul{\dom_{D^F_{\sigma_i}}}(x_i)
		\right) = y_{ij}.
	\]%
	By Lemma~\ref{lemma:ep-presentation} and applying the translation
	\(\T(\C(F))\) to the result,%
	\begin{equation}\label{eq:formula-TCF-epT1}%
		\T\pwrap{\C(F)} \ep_{T_1} \phi(\vec{x}) \dashv\vdash%
		\exists z^{\ul{F\phi}} \bigwedge_{i=1}^n \ul{\C(F)
		\pi^{\vec\sigma}_{\sigma_i}}\pwrap{\ul{\C(F)\dom_{\fb\phi}}(z)} = x_i.
	\end{equation}%

	On the other hand, by Lemma~\ref{lemma:ep-presentation} and \axiom{IL2}
	applied to \(\pi^{\vec{\tau}}_ {\tau_{ij}} = \pi^{\vec{\tau}_i}_{\tau_{ij}}
	\pi^{\vec{\tau}}_{\vec{\tau}_i} \) (coming from the universal property of
	products in \(\C(T_2)\)), we have%
	\[%
		\ep_{T_2} F \phi(\vec{y}_1,\dots,\vec{y}_n) \dashv\vdash \exists z^{\ul
		{F\phi}} \bigwedge_{i=1}^n\bigwedge_{j=1}^{m_i}\ul{\pi^{\vec{\tau}_i}_{\tau_
		{ij}}}\left(\ul{\pi^{\vec{\tau}}_{\vec{\tau}_i}}\left(\ul{\dom_{F\phi}}(z)
		\right)\right) = y_{ij}.
	\]%

	Moreover, from the definition of \(\Eff = \C(F)\), note
	\(\dom_{D^F_{\vec{\sigma}}}\Eff\dom_{[\phi]} = \dom_{F\phi}\) in \(\C(T_2)\).
	Thus we can apply the \axiom{IL2} axiom and first \axiom{IL6} axiom to yield
	the sequent
	\[%
		(\ep_F)_{\vec{\sigma}} (\vec x,\vec y) \wedge%
		\exists z^{\ul{F\phi}} \bigwedge_{i=1}^n
		\ul{\pi^{\vec{\tau}}_{\vec{\tau}_i}}\pwrap{%
			\ul{\dom_{D^F_\sigma}}\pwrap{ \ul{\Eff\dom_{[\phi]}}(z) }%
		} = \ul{\dom_{D^F_{\sigma_i}}}(x_i) \vdash \ep_{T_2} F\phi(\vec{y}).
	\]%
	Since \(D^F_{\vec{\sigma}}(x_1',\dots,x_n') \equiv \bigwedge_{i=1}^n 
	D^F_{\sigma_i}(x_i')\), the following diagram commutes.%
	\[%
		\begin{tikzcd}%
			F\phi \arrow[r,"\Eff\dom_{[\phi]}"{yshift=.1em}] &%
			D^F_{\vec{\sigma}} \arrow [r,"\Eff\pi^{\vec{\sigma}}_{\sigma_i}"]
			\arrow[d,"\dom_{D^F_{\vec{\sigma}}} ",swap] &%
			D_{\sigma_i}^F \arrow[d,"\dom_{D^F_{\sigma_i}}"]\\
			& \fb{\vec{\tau}} \arrow[r,"\pi^{\vec{\tau}}_{\vec{\tau}_i}"] &%
			\fb{\vec {\tau}_i}
		\end{tikzcd}%
	\]%
	Applying \axiom{IL2} axioms to this diagram, the left side of the previous
	sequent is logically equivalent in \(\T\C(T_2)\) to%
	\begin{equation}\label{eq:formula-TM5-post-diagram}%
		\pwrap{\ep_F}_{\vec\sigma}\pwrap{\vec{x},\vec{y}} \wedge%
		\exists z^{\ul{F\phi}}\bigwedge_{i=1}^n \ul{\dom_{D^F_{\sigma_i}}} 
		\pwrap{%
			\ul{\Eff\pi^{\vec{\sigma}}_{\sigma_i}}\pwrap{%
				\ul{\Eff\dom_{[\phi]}}(z)%
			}%
		} = \ul{\dom_{D^F_{\sigma_i}}}(x_i).
	\end{equation}%
	The second conjunct of Formula~\ref{eq:formula-TM5-post-diagram} is the
	result after applying the function symbol \(\ul{\dom_{D^F_{\sigma_i}}}\) to
	all terms in the right side of Sequent~\ref{eq:formula-TCF-epT1}. Therefore%
	\[%
		\pwrap{\ep_F}_{\vec\sigma}\pwrap{\vec{x},\vec{y}} \wedge%
		\T\pwrap{\C(F)}\ep_{T_1}\phi(\vec{x}) \vdash \text{A.2}(\vec{x},\vec{y})
		\vdash \ep_{T_2} F \phi(\vec{y}).
	\]%
	By the cut rule, this proves \axiom[\ep_F]{TM5}. As for \axiom{TM8}, it
	follows a similar argument, where the IL2 axioms are applied in reverse order
	so as to obtain the same function symbols as the ones in the expression for
	\(\ep_{T_2}F\phi\).
\end{proof}
\begin{proof}[Proof of Proposition~\ref{prop:ep-pseudonatural}]
	We need to prove conditions \axiom{PNT1} through \axiom{PNT4} of
	Definition~\ref{def:PNTi}. 

	(\axiom{PNT1}) We need to show that \(\ep_T: T \to \T\C(T)\) is an e.p.\
	translation for every coherent theory \(T\). This was proven in
	Proposition~\ref{prop:epsilon-translation}.

	(\axiom{PNT2}) Consider a t\-/map \(\chi: F \Rightarrow G\) in \(\CThEq\),
	where \(F,G: T_1 \to T_2\). We need to prove the equation%
	\[%
		\ep_G \cdot \pwrap{ \T\C(\chi) \circ \One^{\ep_{T_1}} } =%
		\pwrap{ \One^{\ep_{T_2}} \circ \chi } \cdot \ep_F.
	\]%
	Let \(\sigma\) be a \(T_1\)\-/sort. Proving the above equation amounts to
	showing that the \(\sigma\) components of both sides are presented by
	logically equivalent formulae. Using \axiom[\ep_G]{TM1}, the proof of
	Proposition~\ref{prop:ep_F-definable-iso}, and
	Definition~\ref{def:T-2-cells}, the left side is presented by the formula
	\(
        \pwrap{\ep_G}_\sigma\pwrap{\ul{\C(\chi)_{\fb\sigma}}(s),t}.
	\)
	Using \axiom[\chi]{TM1}, the assumption that \(G\) is e.p.,
	\axiom[\One^{\ep_{T_2}}]{TM3}, and the translation \(\ep_{T_2}\), the right
	side is presented by the formula
	\(
		\exists z^{\ep_{T_2} F\sigma} \pwrap{%
			\ep_{T_2} \chi_\sigma(z,t) \wedge \pwrap{\ep_F}_\sigma(s,z)
		}.
	\)
	Therefore \axiom{PNT2} reduces to proving%
	\begin{equation}\label{eq:PNT2-ep}%
		\pwrap{\ep_G}_\sigma\pwrap{\ul{\C(\chi)_{\fb\sigma}}(s),t}%
		\dashv\vdash%
		\exists z^{\ep_{T_2} F\sigma} \pwrap{%
			\ep_{T_2} \chi_\sigma(z,t) \wedge \pwrap{\ep_F}_\sigma(s,z)
		}.%
	\end{equation}%
	Let \(\vec\omega \equivdef \tau_1,\dots,\tau_n,\upsilon_1,\dots, \upsilon_m\)
	where \(\vec{\tau} \equivdef F\sigma\) and \(\vec{\upsilon} \equivdef
	G\sigma\). Let \(u \equivdef z,t\). Note that \(\Dom \ep_{T_2}\chi_\sigma
	\equiv \ul{\fb{\tau_1}},\dots, \ul{\fb
	{\tau_n}},\ul{\fb{\upsilon_1}},\dots,\ul{\fb{\upsilon_m}}\).
	Lemma~\ref{lemma:ep-presentation} allows us to expand the right side to a
	conjunction of equations involving products
	\(\ul{\pi^{\vec\omega}_{\omega_i}}\) and \(\ul{\dom_{\chi_\sigma}}\). The
	definition of \(\pwrap{\ep_F}_\sigma\) is similar. By the universal property
	of products in \(\C(T_2)\), we have \(\pi^{\vec\omega}_{\omega_i} =
	\pi^{\vec\tau}_{\tau_i} \pi^{\vec\omega}_{\vec\tau}\) for \(1 \leq i \leq n\)
	and \(\pi^{\vec\omega}_{\omega_i} = \pi^{\vec\upsilon}_{\upsilon_{i-n}}
	\pi^{\vec\omega}_{\vec\upsilon}\) for \(n+1 \leq i \leq n + m\). By using
	\axiom{IL2} for these projections, we can relate the projections coming from
	\(\ep_{T_2}\chi_\sigma\) with the projections in the definition of
	\(\pwrap{\ep_F}_\sigma\). With this in mind, we can use \axiom{IL6} to deduce
	that the right side is logically equivalent to the formula%
	\begin{equation}\label{eq:PNT2-LHS}%
		\exists y^{\ul{\chi_\sigma}} \pwrap{%
			\ul{\pi^{\vec\omega}_{\vec\tau}}\pwrap{\ul{\dom_{\chi_\sigma}}(y)} =%
			\ul{\dom_{D^F_\sigma}}(s) \wedge \bigwedge_{j=1}^m
			\ul{\pi^{\vec\upsilon}_{\upsilon_j}}\pwrap{%
				\ul{\pi^{\vec\omega}_{\vec\upsilon}}\pwrap{\ul{\dom_{\chi_\sigma}}(y)}
			} = t_j%
		}.%
	\end{equation}%
	As for the left side, it expands to the formula%
	\begin{equation}\label{eq:PNT2-RHS}%
		\bigwedge_{j=1}^m \ul{\pi^{\vec\upsilon}_{\upsilon_j}}\pwrap{%
			\ul{\dom_{D^G_\sigma}}\pwrap{\ul{\C\pwrap{\chi}_{\fb\sigma}}(s)}
		} = t_j.%
	\end{equation}%
	To show a logical equivalence between these formulae, we first prove the
	 sequent%
	\begin{equation}\label{eq:PNT2-simplification}%
		\ul{\pi^{\vec\omega}_{\vec\tau}}\pwrap{\ul{\dom_{\chi_\sigma}}(y)} =%
		\ul{\dom_{D^F_\sigma}}(s) \vdash%
		\ul{\dom_{D^G_\sigma}}\pwrap{\ul{\C\pwrap{\chi}_{\fb\sigma}}(s)} =%
		\ul{\pi^{\vec\omega}_{\vec\upsilon}}\pwrap{\ul{\dom_{\chi_\sigma}}(y)}.
	\end{equation}%
	This new sequent applied to the right side of Sequent~\ref{eq:PNT2-ep} allows
	us to replace the term
	\(\ul{\pi^{\vec\omega}_{\vec\upsilon}}\pwrap{\ul{\dom_{\chi_\sigma}}(y)}\)
	with \(\ul{\dom_{D^G_\sigma}}\pwrap{\ul{\C\pwrap{\chi}_{\fb\sigma}}(s)}\).
	This establishes the converse of Sequent~\ref{eq:PNT2-ep}.

	Let \(f: \chi_\sigma \twoheadrightarrow D^F_\sigma\) be the morphism in
	\(\C(T_2)\) presented by \(f(x_1,x_2,y) \equivdef \chi_\sigma(x_1,x_2) \wedge
	x_1 = y\). \axiom[\chi]{TM1} implies \(\dom_{D^F_\sigma} f =
	\pi^{\vec\omega}_{\vec\tau} \dom_{\chi_\sigma}\). The \axiom{IL2} axioms for
	this equation yield the sequent%
	\[%
		\ul{\pi^{\vec\omega}_{\vec\tau}}\pwrap{\ul{\dom_{\chi_\sigma}}(y)} =%
		\ul{\dom_{D^F_\sigma}}(s) \vdash \ul{\dom_{D^F_\sigma}}(\ul{f}(y)) =%
		\ul{\dom_{D^F_\sigma}}(s).
	\]%
	Since \(\dom_{D^F_\sigma}\) is monic, \axiom[\dom_{D^F_\sigma}]{IL3} shows
	that the left side of the above sequent entails the formula \(\ul{f}(y) =
	s\). The definition of \(\C(\chi)_{\fb\sigma}\) implies \(\dom_{D^G_\sigma}
	\C(\chi)_{\fb\sigma} f\) equals the morphism
	\(\pi^{\vec\omega}_{\vec\upsilon} \dom_{\chi_\sigma}\). From \(\ul{f}(y) =
	s\), this shows that Sequent~\ref{eq:PNT2-simplification} is provable.

	Returning to the logical equivalence, Sequent~\ref{eq:PNT2-simplification}
	shows that Formula~\ref{eq:PNT2-RHS} entails the formula%
	\[%
		\exists y^{\ul{\chi_\sigma}}\,\,%
		\ul{\dom_{D^G_\sigma}}\pwrap{\ul{\C\pwrap{\chi}_{\fb\sigma}}(s)} =%
		\ul{\pi^{\vec\omega}_{\vec\upsilon}}\pwrap{\ul{\dom_{\chi_\sigma}}(y)}%
		\wedge \bigwedge_{j=1}^m \ul{\pi^{\vec\upsilon}_{\upsilon_j}}\pwrap{%
			\ul{\pi^{\vec\omega}_{\vec\upsilon}}\pwrap{\ul{\dom_{\chi_\sigma}}(y)}
		} = t_j.%
	\]%
	The converse sequent \(\axiom{\ref{eq:PNT2-RHS}} \vdash
	\axiom{\ref{eq:PNT2-LHS}}\) follows by eliminating the variable \(y\) using
	the term involving \(\C\pwrap{\chi}_{\fb\sigma}(s)\). This leaves the forward
	sequent \(\axiom{\ref{eq:PNT2-LHS}} \vdash \axiom{\ref{eq:PNT2-RHS}}\).
	\axiom[\chi]{TM3} implies \(f\) is a regular epimorphism; therefore
	\axiom[f]{IL9}, namely \(\vdash \exists y^{\ul{\chi_\sigma}} \ul{f}(y) = s\),
	is an axiom of \(\T\C(T_2)\). Combining this with Formula~\ref{eq:PNT2-LHS}
	shows that \ref{eq:PNT2-LHS} entails%
	\[%
		\exists y^{\ul{\chi_\sigma}}\,\,%
		\ul{\dom_{D^F_\sigma}}\pwrap{\ul{f}(y)} = \ul{\dom_{D^F_\sigma}}(s)%
		\wedge \bigwedge_{j=1}^m \ul{\pi^{\vec\upsilon}_{\upsilon_j}}\pwrap{%
			\ul{\C\pwrap{\chi}_{\fb\sigma}}\pwrap{\ul{f}(y)}
		} = t_j.%
	\]%
	Combining this with the \axiom{IL2} axioms for the equations
	\(\dom_{D^F_\sigma} f = \pi^{\vec\omega}_{\vec\tau} \dom_{\chi_\sigma}\) and
	\(\dom_{D^G_\sigma}\C\pwrap{\chi}_{\fb\sigma} f =
	\pi^{\vec\omega}_{\vec\upsilon} \dom_{\chi_\sigma}\) shows that the above
	formula entails Formula~\ref{eq:PNT2-RHS}, establishing the forward sequent.
	This completes the proof that~\ref{eq:PNT2-LHS} and~\ref{eq:PNT2-RHS} are
	logically equivalent, which completes the proof of \axiom{PNT2}.

	(\axiom{PNT3}) We need to show that the following equation holds for any
	coherent theory \(T\).%
	\[%
		\pwrap{\One^{\ep_T} \circ \pwrap{\id_\CThEq}_{1_T}} \cdot r^{-1}_{\ep_T}
		\cdot l_{\ep_T} = \ep_{1_T} \cdot \pwrap{ \pwrap{\T\C}_{1_T} \circ
		\One^{\ep_T} }
	\]%
	Given a sort \(\sigma\) of \(T\), expanding the \(\sigma\) components of both
	sides of this equation yields expressions involving only \(E^{\ep_T}\) and
	the equality relation \(=_{\ul{\fb\sigma}}\). By
	Proposition~\ref{prop:epsilon-translation}, both sides are presented by \(s =_{\ul{\fb\sigma}} t\) and thus equal as t\-/maps.

	(\axiom{PNT4}) Given a pair of e.p.\ translations \(T_1 \xrightarrow{F}
	T_2 \xrightarrow{G} T_3\), we must show that the following equation of
	t\-/maps holds.%
	\begin{gather*}%
		\pwrap{\One^{\ep_{T_3}} \circ \pwrap{\id_\CThEq}_{GF}} \cdot%
		a_{\ep_{T_3} GF} \cdot \pwrap{\ep_G \circ \One^F} \cdot%
		a^{-1}_{\T\C(G)\ep_{T_2} F} \cdot\\%
		\pwrap{ \One^{\T\C(G)} \circ \ep_F } \cdot a_{\T\C(G)\T\C(F) \ep_{T_1}} =%
		\ep_{GF} \cdot \pwrap{ \pwrap{\T\C}_{GF} \circ \One^{\ep_{T_1}} }
	\end{gather*}%
	Let \(\sigma\) be a \(T_1\)\-/sort. Since all translations involved are e.p.,
	we can simplify the \(\sigma\) component of both sides significantly. Using
	Lemma~\ref{lemma:t-map-lemma-2} we can collect all domain formulae into one
	formula. Thus the \(\sigma\) component of the left side is presented by the
	formula%
	\[%
		\exists y^{\ul{D^G_{F\sigma}}}\pwrap{%
			D^{\ep_{T_2}GF}_\sigma(t) \wedge%
			\T\pwrap{\C(G)}\pwrap{\ep_F}_\sigma(s,t) \wedge%
			\pwrap{\ep_G}_{F\sigma}(y,t)
		}.%
	\]%
	For the same reason, the \(\sigma\) component of the right side is presented
	by%
	\[%
		\exists w^{\ul{D^{GF}_\sigma}}\pwrap{%
			\pwrap{\pwrap{\T\C}_{GF}}_{\ul{\fb\sigma}}(s,w) \wedge%
			\pwrap{\ep_{GF}}_\sigma(w,t)
		}.%
	\]%
	Since the compositor of \(\C\) is trivial and the identitor of \(\T\) is
	trivial, \(\T(\C_{GF}) = \T(\One^{\C(G)\C(F)}) = \One^{\T(\C(G)\C(F))}\), so we arrive at the result that
        \[%
        (\T\C)_{GF} = \T(\C_{GF}) \cdot \T_{\C(G)\C(F)} = \One^{\T(\C(G)\C(F))} \cdot \T_{\C(G)\C(F)} = \T_{\C(G)\C(F)}.
        \]%
	Thus \(((\T\C)_{GF})_{\ul{\fb{\sigma}}}(s,w) \equiv s = w\), so the right
	side of the original equation reduces to
	\(\left(\ep_{GF}\right)_\sigma\!(s,t)\). 

	Now that we have simplified the left and right \(\sigma\) components, we
	prove that they are logically equivalent. We will omit the domain formula on
	the left side, since it can be recovered from \axiom[\ep_{GF}]{TM1}. Let
	\(\vec\tau \equivdef F\sigma\) and \(G\tau_i \equivdef \vec{\alpha_i}
	\equivdef \alpha_{i 1},\ldots,\alpha_{i m_i}\) so that \(GF\sigma \equiv
	\vec\alpha\). Invoking Rule~\ref{rule:term-reduction} twice, the left side is
	logically equivalent to the following pair of formulae.%
	\begin{gather*}%
		\exists y^{\ul{D^G_{\vec\tau}}} \pwrap{%
			\bigwedge_{i=1}^n \ul{\C(G)\pi^{\vec\tau}_{\tau_i}}\pwrap{%
				\ul{\C(G)\dom_{D^F_\sigma}}(s)
			} = y_i \wedge \bigwedge_{i=1}^n \pwrap{\ep_G}_{\tau_i}(y_i,t_i)%
		}\\%
		\dashv\vdash \bigwedge_{i=1}^n \bigwedge_{j=1}^{m_i}%
		\ul{\pi^{\vec{\alpha_i}}_{\alpha_{i j}}}\pwrap{%
			\ul{\dom_{D^G_{\tau_i}}}\pwrap{%
				\ul{\C(G)\pi^{\vec\tau}_{\tau_i}}\pwrap{%
					\ul{\C(G)\dom_{D^F_\sigma}}(s)
				}%
			}%
		} = t_{i j}%
	\end{gather*}%
	Meanwhile, the right side is presented by \(\pwrap{\ep_{GF}}_\sigma(s,t)\),
	which expands to%
	\[%
		\bigwedge_{i=1}^n \bigwedge_{j=1}^{m_i}%
		\ul{\pi^{\vec\alpha}_{\alpha_{i j}}}\pwrap{%
			\ul{\dom_{D^{GF}_\sigma}}(s)
		} = t_{i j}.%
	\]%
	Since both sides are a conjunction with the same number of terms, it suffices
	to show that each pair of conjuncts of the same index are logically
	equivalent, i.e.,%
	\[%
		\ul{\pi^{\vec{\alpha_i}}_{\alpha_{i j}}}\pwrap{%
			\ul{\dom_{D^G_{\tau_i}}}\pwrap{%
				\ul{\C(G)\pi^{\vec\tau}_{\tau_i}}\pwrap{%
					\ul{\C(G)\dom_{D^F_\sigma}}(s)
				}%
			}%
		} = t_{i j} \dashv\vdash%
		\ul{\pi^{\vec\alpha}_{\alpha_{i j}}}\pwrap{%
			\ul{\dom_{D^{GF}_\sigma}}(s)
		} = t_{i j}.%
	\]%
	Since both sides contain the variable \(t_{i j}\) it suffices to prove that the
	following equation holds in \(\C(T_3)\)%
	\[%
		\dom_{D^G_{\tau_i}} \pwrap{\C(G)\pi^{\vec\tau}_{\tau_i}}
		\pwrap{\C(G)\dom_{D^F_\sigma}} = \dom_{D^{GF}_\sigma}.
	\]%
	Unpacking both sides of this equation to defining formulae, this equation is
	equivalent to showing the following pair of formulae are logically
	equivalent:%
	\[%
		D^{GF}_\sigma(\vec x) \wedge D^G_{\tau_i}(x_i) \dashv\vdash
		D^{GF}_\sigma(x).
	\]%
	This is a consequence of Lemma~\ref{lemma:t-map-lemma-1}. Thus \axiom{PNT4}
	is proven, so \(\ep\) is pseudonatural.
\end{proof}

\begin{proof}[Proof of Proposition~\ref{prop:delta-pseudonatural}]
	As for \(\ep\), we show that \(\delta\) is a pseudonatural transformation by
	proving conditions \axiom{PNT1} through \axiom{PNT4} of
	Definition~\ref{def:pseudofunctor}.

	(\axiom{PNT1}) The \axiom{IL1} and \axiom{IL2} axioms ensure that
	\(\delta_C: C \to \C\T(C)\) is a functor for any coherent category \(C\). To
	see that this functor is coherent, it suffices to show that every diagram
	mentioned by the \axiom{IL} axiom schemata is preserved by \(\delta_C\). This
	follows almost immediately from Proposition~\ref{prop:sc-axiom-schemata}.
	
	(\axiom{PNT2}) For any natural transformation \(\chi : \Eff \Rightarrow \Gee\) with
	\(\Eff,\Gee : C_1 \to C_2\), we need the equation%
	\[%
		\delta_{\Gee} \cdot \pwrap{%
			\C\T\pwrap\chi \circ \One^{\delta_{C_1}}
		} = \pwrap{%
			\One^{\delta_{C_2}} \circ \chi
		} \cdot \delta_{\Eff}.
	\]%
	By Proposition~\ref{prop:delta-F-one}, it suffices to show \(\C\T\pwrap\chi
	\circ \One^{\delta_{C_1}} = \One^{\delta_{C_2}} \circ \chi\). We conclude,
	for they are maps from \(\delta_{C_2}\Eff\) to \(\delta_{C_2}\Gee\) with the
	same components:%
	\begin{gather*}
		\pwrap{\C\T\pwrap\chi \circ \One^{\delta_{C_1}}}_A =%
		\pwrap{\C\T\pwrap\chi}_{\fb{\ul A}} \cdot 1_{\fb{\ul A}} =%
		\pwrap{\C\T\pwrap\chi}_{\fb{\ul A}} = \theta_{\ul{\chi_A}}\\
		\pwrap{\One^{\delta_{C_2}}\circ\chi}_A =%
		1_{\fb{\ul{\Gee A}}} \cdot \delta_{C_2}\pwrap{\chi_A} =%
		\delta_{C_2}\pwrap{\chi_A} = \theta_{\ul{\chi_A}}.
	\end{gather*}
	
	(\axiom{PNT3}) For any coherent category \(C\), we must have the equation%
	\[%
		\pwrap{\One^{\delta_C} \circ \pwrap{\id_{\Coh}}_{1_C}} \cdot
		r^{-1}_{\delta_C} \cdot l_{\delta_C} =%
		\delta_{1_C} \cdot \pwrap{\pwrap{\C\T}_{1_C} \circ \One^{\delta_C}}.
	\]%
	By Proposition~\ref{prop:delta-F-one} and the triviality of the unitors and
	associators of \(\Coh\), it suffices to prove \(\One^{\delta_C} =
	\pwrap{\C\T}_{1_C} \circ \One^{\delta_C}\).
	Proposition~\ref{prop:delta-F-one} confirms that \(\pwrap{\C\T}_{1_C} \circ
	\One^{\delta_C}\) is also a map from \(\delta_C \Rightarrow \delta_C\). Thus,
	we show that the two maps agree along the component of an object \(A\). The
	first expression yields \(\One^{\delta_C}_A = 1_{\delta_C A}\).
	For the second, we note that%
	\[%
		\pwrap{\C\T}_{1_C} = \C\pwrap{\One^{\T(1_C)}} \circ \C_{1_{\T(C)}} =%
		\One^{\C\T(1_C)} \circ \C_{1_{\T(C)}} = \One^{\delta_C} \circ
		\C_{1_{\T(C)}},
	\]%
	using Definition~\ref{def:pseudofunctor} for the identitor, Proposition
	\ref{prop:T-identitor}, \axiom[\C]{PF2}, and
	Proposition~\ref{prop:delta-F-one}. By Proposition~\ref{prop:C-identitor},
	\(\pwrap{\C_{1_{\T(C)}}}_ {\fb{\ul{A}}}= 1_{\delta_C A}\), so by
	Proposition~\ref{prop:horizontal-composition-unital} we conclude that
	\(\pwrap{\pwrap{\C\T}_{1_C} \circ \One^{\delta_C}}_A = \One^{\delta_C}_A
	\cdot \pwrap{\C_{1_{\T(C)}}}_{\fb{\ul A}} \cdot \One^{\delta_C}_A =
	1_{\delta_C A}\), as desired.

	(\axiom{PNT4}) Consider a pair of coherent functors \(C_1 \xrightarrow{\Eff}
	C_2 \xrightarrow{\Gee} C_3\). Since associators are trivial in \(\Coh\) and
	Proposition~\ref{prop:delta-F-one} shows that \(\delta_{\Eff}\) and
	\(\delta_\Gee\) are trivial, proving the commutative diagram for
	\axiom[\delta]{PNT4} simplifies to showing the equation
	\(\One^{\delta_{C_3}\Gee\Eff} =
	\C\pwrap{\T_{\Gee\Eff}}\circ\One^{\delta_{C_1}}\)
	(Indeed, Proposition~\ref{prop:delta-F-one} shows that both sides of this equation have the same
	domain and codomain.) It suffices to show that
	both sides have equal components along an object \(A\); hence it suffices to
	show that both sides are presented by logically equivalent formulae. The left
	side is presented by \(\One^{\delta_{C_3}\Gee\Eff}_A(s,t) \dashv\vdash s =
	t\), where \(s\) and \(t\) belong to the sort \(\ul{\Gee\Eff A}\). Since the
	compositor of \(\C\) is trivial, we can expand the \(A\) component of the
	right side using Definition~\ref{def:pseudofunctor-compose}
	and~\axiom[\C]{PF2} to deduce%
	\[%
		\pwrap{\pwrap{\C\T}_{\Gee\Eff}\circ \One^{\delta_{C_1}}}_A =%
		\C\pwrap{\T_{\Gee\Eff}}_{\delta_{C_1} A} \cdot
		\C\pwrap{\T\pwrap{\Gee\Eff}}\pwrap{1_{\delta_{C_1}} A} =%
		\C\pwrap{\T_{\Gee\Eff}}_{\fb{\ul A}}.
	\]%
	The latter is also presented by \(s = t\), as desired.
\end{proof}

\subsection{Properties of Exact Completions}%
\label{subsection:exact-completion-properties}

\begin{proof}[Proof of Proposition~\ref{prop:exact-completion-2-cells}]
	We need to define the components of \(\chi^\ex\). Recall there are two classes of
	objects in the exact completion \(C^\ex\).%
	\begin{enumerate}
		\item Objects of the form \(I A\), where \(A\) is an object of \(C\) (and
		\(I: C \to C^\ex\) is the inclusion functor in
		Proposition~\ref{prop:exact-completion}).
		\item Quotients \(I(A) / I(R)\), where \(R \hookrightarrow A \times A\) is a
		congruence in \(C\).
	\end{enumerate}
	For objects in the first class, we extend \(\chi\) by setting \(\chi^\ex_{I A}
	\eqdef \chi_A\). The triangle in \axiom{EC3} of
	Proposition~\ref{prop:exact-completion} ensures that this makes sense. For
	objects in the second class, we need to diagram chase. Consider a quotient
	\(I(A)/I(R)\) in \(C^\ex\), and abbreviate this object to \(A/R\). Since \(\Eff\)
	and \(\Gee\) are coherent functors, \(\Eff R\) and \(\Gee R\) are congruences over \(\Eff A\)
	and \(\Gee A\) respectively. Since \(D\) is Barr\-/exact, this means that there
	exist objects \(Q^\Eff_A\) and \(Q^\Gee_A\) of \(D\) and quotient morphisms \(q^\Eff_A:
	\Eff A \to Q^\Eff_A\) and \(q^\Gee_A: \Gee A \to Q^\Gee_A\) which coequalize \(\Eff R\) and \(\Gee R\)
	respectively. Coherent functors preserve quotients (since coherent functors
	preserve images and pullbacks), so \(\Eff^\ex(A/R)\) is isomorphic to \(Q^\Eff_A\)
	and there exists a bijective correspondence between diagrams in \(D\)%
	\[%
		\begin{array}{cc}%
			\begin{tikzcd}
				\Eff^\ex(A/R) \arrow[d,dashed] \\
				\ast
			\end{tikzcd} \mapsto &
			\begin{tikzcd}
				\Eff R \arrow[r] \arrow[d] & \Eff A \arrow[d,dashed] \\
				\Eff A\arrow[r,dashed] & \ast
			\end{tikzcd}
		\end{array}%
	\]%
	between a morphism from \(\Eff^\ex(A/R)\) and a morphism which coequalizes
	the projections of the congruence \(\Eff R \rightrightarrows \Eff A\). The
	analogous correspondence is true for \(\Gee^\ex(A/R)\). Using the morphisms
	\(\chi_R\) and \(\chi_A\), we stitch three squares together
	(Figure~\ref{fig:naturality-cube}, left). 
	\begin{figure}[h]
		\tdplotsetmaincoords{70}{60}
		\begin{tikzpicture}%
			[scale=2,tdplot_main_coords, every node/.style={scale=1}]
			\coordinate [label=above left:$\Eff R$] (FR) at
				(0,0,1);
			\coordinate [label=above:$\Eff A$] (FPhi1) at (0,1,1);
			\coordinate [label=above:$\Gee R$] (GR) at (1,0,1);
			\coordinate [label=above right:$\Gee A$] (GPhi1) at (1,1,1);
			\coordinate [label=below left:$\Eff A$] (FPhi2) at (0,0,0);
			\coordinate [label=below:$\Gee A$] (GPhi2) at (1,0,0);
			\coordinate [label=below right:$\Gee^\ex(A/R)$]
				(Gquot) at (1,1,0);

			\begin{scope}%
				[decoration={markings,mark=at position 0.5 with {\arrow{to}}}]
				\draw[postaction={decorate}] (FR) -- (FPhi1);
				\draw[postaction={decorate}] (FR) -- (GR) node[midway, below]
					{\scriptsize$\chi_R$};
				\draw[postaction={decorate}] (FPhi1) -- (GPhi1) node[midway, above
					right] {\scriptsize$\chi_A$};
				\draw[postaction={decorate}] (GR) -- (GPhi1);

				\draw[postaction={decorate}] (FR) -- (FPhi2);
				\draw[postaction={decorate}] (GR) -- (GPhi2);
				\draw[postaction={decorate}] (FPhi2) -- (GPhi2) node[midway,below left]
					{\scriptsize$\chi_A$};

				\draw[postaction={decorate}] (GPhi1) -- (Gquot) node[midway,right]
					{\scriptsize$q^\Gee_A$};
				\draw[postaction={decorate}] (GPhi2) -- (Gquot) node[midway,below]
					{\scriptsize$q^\Gee_A$};
			\end{scope}%
		\end{tikzpicture}
		\begin{tikzpicture}%
			[scale=2,tdplot_main_coords, every node/.style={scale=1}]
			\coordinate [label=above left:$\Eff R$] (FR) at
				(0,0,1);
			\coordinate [label=above:$\Eff A$] (FPhi1) at (0,1,1);
			\coordinate [label= above:$\Gee R$] (GR) at (1,0,1);
			\coordinate [label= right:$\Gee A$] (GPhi1) at (1,1,1);
			\coordinate [label=below left:$\Eff A$] (FPhi2) at (0,0,0);
			\coordinate [label=below:$\Gee A$] (GPhi2) at (1,0,0);
			\coordinate [label=below right:$\Gee^\ex (A/R)$] (Gquot) at (1,1,0);
			\coordinate (Fquot) at (0,1,0);

			\begin{scope}%
				[decoration={markings,mark=at position 0.5 with {\arrow{to}}}]
				\draw[postaction={decorate}] (FR) -- (FPhi1);
				\draw[postaction={decorate}] (FR) -- (GR) node[midway, below]
					{\scriptsize$\chi_R$};
				\draw[postaction={decorate}] (FPhi1) -- (GPhi1) node[midway, above right]
                {\scriptsize$\chi_A$};
				\draw[postaction={decorate}] (GR) -- (GPhi1);

				\draw[postaction={decorate}] (FR) -- (FPhi2);
				\draw[postaction={decorate}] (GR) -- (GPhi2);
				\draw[postaction={decorate}] (FPhi2) -- (GPhi2) node[midway,below left]
					{\scriptsize$\chi_A$};

				\draw[postaction={decorate}] (GPhi1) -- (Gquot) node[midway,right]
					{\scriptsize$q^\Gee_A$};
				\draw[postaction={decorate}] (GPhi2) -- (Gquot) node[midway,below]
					{\scriptsize$q^\Gee_A$};

				\draw[postaction={decorate}] (FPhi1) -- (Fquot) node[midway,below
					right] {\scriptsize$ q^\Eff_A$};
				\draw[postaction={decorate}] (FPhi2) -- (Fquot) node[midway,above left]
					{\scriptsize$q^\Eff_A$};
				\draw[postaction={decorate},dashed] (Fquot) -- (Gquot);
			\end{scope}%
		\end{tikzpicture}
		\caption{%
			\label{fig:naturality-cube}%
			Naturality and Quotient Diagrams
		}%
	\end{figure}
	All the faces commute because
	\(\chi\) is a natural transformation. Therefore the outer hexagon gives a
	fourth commutative square.%
	\[%
		\begin{tikzcd}
			\Eff R \arrow[r] \arrow[d] & \Eff A \arrow[d,"q^\Gee_A\circ\chi_A"]\\
			\Eff A \arrow[r,"q^\Gee_A\circ\chi_A",swap] & \Gee^\ex(A/R)
		\end{tikzcd}
	\]%
	We invoke the universal property of \(q^\Eff_A\) to make a commutative cube
	(Figure~\ref{fig:naturality-cube}, right). Set \(\chi^\ex_{A/R}: \Eff^\ex
	(A/R) \to \Gee^\ex(A/R)\) to be the morphism represented by the dashed line.
	This completes the definition of \(\chi^\ex\). The fact that \(\chi^\ex\) is
	a natural transformation can be proven via diagram chasing using the
	naturality of \(\chi\) and the bijective correspondence mentioned earlier. By
	juxtaposing cubes in front of the other, we see that \(\chi \mapsto
	\chi^\ex\) preserves vertical composition. In the case that \(\chi =
	\One^\Eff\), the induced map \(\chi^\ex_{A/R}\) is just the identity, per the
	bijective correspondence mentioned earlier. Thus \(\chi \mapsto \chi^\ex\) is
	a functor.
\end{proof}

\end{document}